\documentclass[a4paper,11pt]{article}

\usepackage{epsfig,amsmath,amsfonts,amssymb,graphics,xypic}

\def\dst{\displaystyle}
\def\tst{\textstyle}
\def\sst{\scriptstyle}
\def\ens{\enspace}

\def\al{\alpha}
\def\be{\beta}
\def\ga{\gamma}
\def\Ga{{\Gamma}}
\def\de{\delta}

\def\De{\Delta}
\def\eps{{\varepsilon}}

\def\la{\lambda}
\def\La{\Lambda}
\def\om{\omega}
\newcommand{\dom}{{\overset{\raisebox{-.23ex}{$\scriptscriptstyle\bullet$}}{\om}}}
\def\Om{\Omega}

\def\sig{{\sigma}}
\def\Sig{{\Sigma}}

\def\th{{\theta}}
\def\Th{\Theta}
\newcommand{\ph}{\varphi}

\def\ze{{\zeta}}
\newcommand\dze{\makebox[.6em]{
         \raisebox{0ex}[1.4ex][0ex]{$
            \begin{array}[b]{c}
                {\scriptscriptstyle\bullet}\\ [-1.1ex] \ze
            \end{array}
         $}
         }}
\newcommand\sdze{\makebox[.4em]{
         \raisebox{0ex}[1.4ex][0.3ex]{$
            \begin{array}[b]{c}
                {\scriptscriptstyle\bullet}\\ [-1.7ex] {\scriptstyle \ze}
            \end{array}
         $}
         }}
\newcommand\dla{\makebox[.6em]{
         \raisebox{0ex}[1.4ex][0ex]{$
            \begin{array}[b]{c}
                {\scriptscriptstyle\bullet}\\ [-1.1ex] \la
            \end{array}
         $}
         }}

\newcommand{\bino}[2]{\mbox{$
\begin{pmatrix}#1\\#2\end{pmatrix}%
                $}}

\newcommand{\demi}{\frac{1}{2}}
\newcommand{\dem}{\tfrac{1}{2}}

\newcommand{\ao}{\{\,}          
\newcommand{\lao}{\left\{\,}          
\newcommand{\af}{\,\}}          
\newcommand{\laf}{\,\right\}}          

\def\ov{\overline}
\def\wh{\widehat}

\newcommand{\const}{\mathop{\hbox{{\rm const}}}\nolimits}

\newcommand{\dist}{\operatorname{dist}}

\newcommand{\card}{\operatorname{card}}

\newcommand{\Sing}{\operatorname{sing}}
\newcommand{\cont}{\operatorname{cont}}

\def\ID{\mathop{\hbox{{\rm Id}}}\nolimits}

\newcommand{\I}{{\mathrm i}}
\newcommand{\dd}{{\mathrm d}}

\newcommand{\ee}{\mathrm e}

\def\pa{\partial}

\def\ii{^{-1}}
\newcommand{\ti}{\tilde}

\def\IM{\mathop{\Im m}\nolimits}
\def\RE{\mathop{\Re e}\nolimits}

\newcommand{\proof}{\noindent{\em Proof.}\ }
\newcommand{\idproof}{\noindent{\em Idea of the proof.}\ }
\newcommand{\eopf}{
\hfill\begin{picture}(.24,.24)\thinlines
\put(0,0){\line(1,0){.24}}
\put(.24,0){\line(0,1){.24}}
\put(.24,.24){\line(-1,0){.24}}
\put(0,.24){\line(0,-1){.24}}
\end{picture}\bigskip }

\def\ie{{\it i.e.}\ }
\def\cf{{\it cf.}\ }
\def\eg{{\it e.g.}\ }
\def\resp{{resp.}\ }
\def\wrt{{with respect to}}
\def\lhs{{left-hand side}}
\def\rhs{{right-hand side}}

\def\dst{\displaystyle}
\def\tst{\textstyle}
\def\sst{\scriptstyle}

\newcommand{\refeq}[1]{{\rm (\ref{#1})}}
\newcommand{\refeqa}[2]{{\rm (\ref{#1}#2)}}

\newcommand{\C}{\mathbb{C}}
\newcommand{\Clog}{{\underset{\raisebox{.4ex}{$\scriptscriptstyle\bullet$}}{\C}}}

\newcommand{\G}{\mathbb{G}}
\newcommand{\N}{\mathbb{N}}

\newcommand{\R}{\mathbb{R}}

\newcommand{\Z}{\mathbb{Z}}

\def\cB{\mathcal{B}}

\def\cD{\mathcal{D}}
\def\cE{\mathcal{E}}
\def\cF{\mathcal{F}}
\def\cG{\mathcal{G}}
\renewcommand{\cH}{\mathcal{H}}

\def\cL{\mathcal{L}}

\def\cO{\mathcal{O}}

\def\cR{\mathcal{R}}

\def\cV{\mathcal{V}}
\def\cW{\mathcal{W}}

\newtheorem{thm}{Theorem}

\newtheorem{lemma}{Lemma}

\newtheorem{prop}{Proposition}

\newtheorem{cor}{Corollary}

\newtheorem{Def}{Definition}

\newcommand{\subsub}[1]{%
\subsubsection*{{\em #1}}%
\addcontentsline{toc}{subsubsection}{{\em #1}}}

\newcounter{parag}[subsection]

\newcounter{parage}

\newcounter{paraga}

\def\ss#1{^{(#1)}}
\def\sst#1{^{[#1]}}

\def\abs#1{\lvert#1\rvert}

\def\om{\omega}
\def\cS{{\cal S}}

\newcommand{\SM}{^{\text{sm}}}
\newcommand{\QU}{^{\text{qu}}}

\newcommand{\dDe}[1]{
  \raisebox{.23ex}{
  {$\stackrel{\raisebox{-.23ex}{$\scriptscriptstyle\bullet$}}\Delta_{
  \raisebox{-.23ex}{$\scriptstyle#1$}}$}}
}

\newcommand{\dDep}[1]{
  \raisebox{.23ex}{
  {$\stackrel{\raisebox{-.23ex}{$\scriptscriptstyle\bullet$}}\Delta_{
  \raisebox{-.23ex}{$\scriptstyle#1$}}^{
  \raisebox{-.85ex}{$\scriptstyle +$}}$}}
}

\newcommand{\dDeom}{
  \raisebox{.23ex}{
  {$\stackrel{\raisebox{-.23ex}{$\scriptscriptstyle\bullet$}}\Delta_{
  \raisebox{-.23ex}{$\scriptstyle\om$}}$}}
}

\newcommand{\hta}[1]{{\stackrel{\raisebox{-.23ex}{$\scriptscriptstyle\wedge$}}{#1}}}

\newcommand{\tr}[1]{{\stackrel{\raisebox{-.23ex}{$\scriptscriptstyle\triangledown$}}{#1}}}

\newcommand{\htb}[1]{\raisebox{-.23ex}{${\stackrel{
            \raisebox{-.18ex}{$\scriptscriptstyle\wedge$}
          }{#1}
     }$}}
\newcommand{\shtb}[1]{
        \raisebox{0.15ex}[2ex][0ex]{${\stackrel{
            \raisebox{-.18ex}{$\scriptscriptstyle\wedge$}
          }{\scriptstyle#1}
     }$}}
\newcommand{\chb}[1]{\raisebox{-.23ex}{${\stackrel{
            \raisebox{-.23ex}{$\scriptscriptstyle\vee$}
          }{#1}
     }$}}
\newcommand{\chbph}[1]{ \raisebox{-.23ex}{${\stackrel{
            \raisebox{-.23ex}{$\scriptscriptstyle\vee$}
          }{\ph}
     }$}^{\raisebox{-.63ex}{$\scriptstyle#1$}} }
\newcommand{\trb}[1]{\raisebox{-.23ex}{${\stackrel{
            \raisebox{-.23ex}{$\scriptscriptstyle\triangledown$}
          }{#1}
     }$}}

\newcommand{\htn}[1]{\raisebox{.23ex}{${\stackrel{
            \raisebox{-.23ex}{$\scriptscriptstyle\wedge$}
          }{#1}
     }$}}
\newcommand{\chn}[1]{\raisebox{.23ex}{${\stackrel{
            \raisebox{-.23ex}{$\scriptscriptstyle\vee$}
          }{#1}
     }$}}
\newcommand{\trn}[1]{\raisebox{.23ex}{${\stackrel{
            \raisebox{-.23ex}{$\scriptscriptstyle\triangledown$}
          }{#1}
     }$}}

\newcommand{\htx}[1]{\raisebox{-.09ex}{${\stackrel{
            \raisebox{-.23ex}{$\scriptscriptstyle\wedge$}
          }{#1}
     }$}}
\newcommand{\chx}[1]{\raisebox{-.09ex}{${\stackrel{
            \raisebox{-.23ex}{$\scriptscriptstyle\vee$}
          }{#1}
     }$}}
\newcommand{\trx}[1]{\raisebox{-.09ex}{${\stackrel{
            \raisebox{-.23ex}{$\scriptscriptstyle\triangledown$}
          }{#1}
     }$}}

\newcommand{\chcR}[1]{\raisebox{.23ex}{${\stackrel{
            \raisebox{-.23ex}{$\scriptscriptstyle\vee$}
          }{\cR}
     }$}^{\raisebox{-.85ex}{$\scriptstyle #1$}}}

\newcommand{\ChcR}{\raisebox{.23ex}{${\stackrel{
            \raisebox{-.23ex}{$\scriptscriptstyle\vee$}
          }{\cR}
     }$}}

\newcommand{\trRES}{\raisebox{.01ex}{${\stackrel{
            \raisebox{-.23ex}{$\scriptscriptstyle\bigtriangledown$}
          }{\operatorname{RES}}
     }$}}
\newcommand{\trRESint}[1]{\raisebox{.01ex}{${\stackrel{
            \raisebox{-.23ex}{$\scriptscriptstyle\bigtriangledown$}
          }{\operatorname{RES}}}$}\,^{\mathrm{int}}_{#1}}
\newcommand{\trRESsimp}[1]{\raisebox{.01ex}{${\stackrel{
            \raisebox{-.23ex}{$\scriptscriptstyle\bigtriangledown$}
          }{\operatorname{RES}}}$}\,^{\mathrm{simp}}_{#1}}

\newcommand{\whRES}{\operatorname{\widehat{RES}}}

\newcommand{\Rsimp}{\operatorname{RES^{\mathrm{simp}}}}

\newcommand{\trRsramOm}{\raisebox{.01ex}{${\stackrel{
            \raisebox{-.23ex}{$\scriptscriptstyle\bigtriangledown$}
          }{\operatorname{RES}}}$}\,\!^{\raisebox{.3ex}{\scriptsize{s.ram.}}}_{\Om}}

\newcommand{\trRsramZ}{\raisebox{.01ex}{${\stackrel{
            \raisebox{-.23ex}{$\scriptscriptstyle\bigtriangledown$}
          }{\operatorname{RES}}}$}\,\!^{\raisebox{.3ex}{\scriptsize{s.ram.}}}_{2\pi\I\Z}}

\newcommand{\wtRsimp}{\operatorname{\widetilde{RES}\,\!^{\mathrm{simp}}}}

\newcommand{\cRexp}{\operatorname{\hat\cR_{\text{exp}}}}
\newcommand{\cRexpom}{\operatorname{\hat\cR_{\text{exp}}^\om}}

\newcommand{\fracC}{\C(\!(z\ii)\!)_1}

\newcommand{\SING}{\operatorname{SING}}
\newcommand{\ANA}{\operatorname{ANA}}
\newcommand{\integ}{^{\operatorname{int}}}
\newcommand{\simp}{^{\mathrm{simp}}}
\newcommand{\sram}{^{\mathrm{s.ram.}}}

\newcommand{\bem}{\mbox{}^\flat\hspace{-.8pt}}
\newcommand{\sing}{\operatorname{sing}}
\newcommand{\var}{\operatorname{var}}

\def\Bbibitem#1#2{\bibitem[#1]{#2}}

\newcommand{\med}{\medskip}

\newcommand{\sm}{\smallskip}

\begin{document}


\thispagestyle{empty}

\begin{center}

{\bf \huge
Resurgent functions and splitting problems}\\[6ex]

{\large David Sauzin (CNRS--IMCCE, Paris)}\\[6ex]

{\small April 2006}

\end{center}


\vspace{1cm}

\begin{abstract}

The present text is an introduction to \'Ecalle's theory of resurgent functions
and alien calculus, in connection with problems of exponentially small
separatrix splitting.
An outline of the resurgent treatment of Abel's equation for resonant dynamics
in one complex variable is included.
Some proofs and details are omitted.
The emphasis is on examples of nonlinear difference equations, as a simple and
natural way of introducing the concepts.

\end{abstract}

\vspace{1.5cm}


\tableofcontents

\pagebreak


\suppressfloats

\section{The algebra of resurgent functions}

Our first purpose is to present a part of \'Ecalle's theory of resurgent functions and
alien calculus in a self-contained way. Our main sources are the series of books
\cite{Eca81} (mainly the first two volumes), a course taught by Jean \'Ecalle at
Paris-Sud university (Orsay) in 1996 and the book~\cite{CNP}.


\subsection{Formal Borel transform}	\label{secFBT}


A resurgent function can be viewed as a special kind of power series, the radius
of convergence of which is zero, but which can be given an analytical meaning
through Borel-Laplace summation.
It is convenient to deal with formal series ``at infinity'', \ie with elements
of~$\C[[z\ii]]$. We denote by~$z\ii\C[[z\ii]]$ the subset of formal series
without constant term.

\begin{Def}
The formal Borel transform is the linear operator
\begin{equation}	\label{eqdefBorel}
\cB\,:\;
\ti\ph(z) = \sum_{n\ge0} c_n z^{-n-1} \in z\ii\C[[z\ii]]
\ens\mapsto\ens
\hat\ph(\ze) = \sum_{n\ge0} c_n \frac{\ze^n}{n!} \in \C[[\ze]].
\end{equation}
\end{Def}

Observe that if $\ti\ph(z)$ has nonzero radius of convergence, say if $\ti\ph(z)$
converges for $|z\ii|<\rho$, then $\hat\ph(\ze)$ defines an entire function,
of exponential type in every direction: if $\tau>\rho\ii$, then
$|\hat\ph(\ze)|\le\const\,\ee^{\tau|\ze|}$ for all $\ze\in\C$.

\begin{figure}

\begin{center}

\epsfig{file=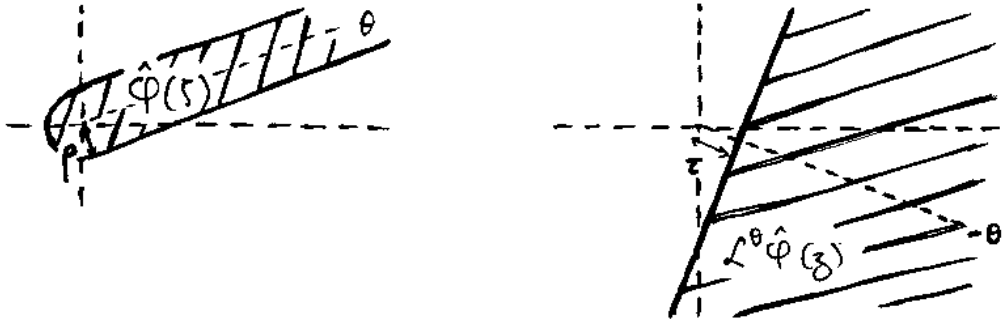,height=5cm,angle = 0}

\end{center}

\vspace{-.5cm}

\caption{\label{figstriphalf} Laplace integral in the direction~$\th$ gives rise
to functions analytic in the half-plane $\RE (z\,\ee^{\I\th})>\tau$.}

\end{figure}

\begin{Def}
For any $\th\in\R$, we define the Laplace transform in the direction~$\th$ as
the linear operator~$\cL^\th$,
\begin{equation}	\label{eqdefLapl}
\cL^\th \hat\ph(z) = \int_0^{\ee^{\I\th}\infty} \hat\ph(\ze)\,\ee^{-z\ze}\,\dd\ze.
\end{equation} 
Here, $\hat\ph$ is assumed to be a function such that
$r\mapsto\hat\ph(r\,\ee^{\I\th})$ is analytic on~$\R^+$ and 
$|\hat\ph(r\,\ee^{\I\th})|\le\const\,\ee^{\tau r}$.
The function~$\cL^\th\hat\ph$ is thus analytic in the half-plane 
$\RE (z\,\ee^{\I\th})>\tau$ (see Figure~\ref{figstriphalf}).
\end{Def}

Recall that $z^{-n-1} =
\int_0^{+\infty}\frac{\ze^n}{n!}\,\ee^{-z\ze}\,\dd\ze$ for $\RE z>0$,
thus
\begin{equation}	\label{eqLaplmon}
z^{-n-1} = \cL^\th\left(\frac{\ze^n}{n!}\right), \qquad \RE (z\,\ee^{\I\th})>0.
\end{equation}
(For that reason, $\cB$ is sometimes called ``formal inverse Laplace transform''.)
As a consequence, if $\hat\ph$ is an entire function of exponential type in
every direction, that is if $\hat\ph=\cB\ti\ph$ with $\ti\ph(z)\in
z\ii\C\{z\ii\}$, we recover~$\ti\ph$ from~$\hat\ph$ by applying the Laplace
transform:
it can be shown\footnote{	\label{footomitpf}
Here, as sometimes in this text, we omit the details of the proof.
See \eg~\cite{Mal} for the properties of the Laplace and Borel transforms.
}
that $\cL^\th\hat\ph(z)=\ti\ph(z)$ for all~$z$ and~$\th$ such
that $\RE (z\,\ee^{\I\th})$ is large enough.

\subsub{Fine Borel-Laplace summation}

Suppose now that $\cB\ti\ph = \hat\ph \in \C\{\ze\}$ but $\hat\ph$ is not
entire, \ie $\hat\ph$ has finite radius of convergence. The radius of
convergence of~$\ti\ph$ is then zero.
Still, it may happen that $\hat\ph(\ze)$ extends analytically to a half-strip 
$\ao \ze\in\C \mid \dist(\ze,\ee^{\I\th}\,\R^+)\le \rho \af$,
with exponential type less than a $\tau\in\R$.
In such a case, formula~\refeq{eqdefLapl} makes sense and the formal
series~$\ti\ph$ appears as the asymptotic expansion of~$\cL^\th\hat\ph$ in the
half-plane $\ao\RE (z\,\ee^{\I\th})>\max(\tau,0)\af$
(as can be deduced from~\refeq{eqLaplmon})\footnote{
See footnote~\ref{footomitpf}.
}.
This is more or less the classical definition of a ``Borel-summable'' formal
series~$\ti\ph$.
One can consider the function~$\cL^\th\cB\ti\ph$ as a ``sum'' of~$\ti\ph$,
associated with the direction~$\th$.
This summation is called ``fine'' when~$\hat\ph$ is only known to extend to a
half-strip in the direction~$\th$, which is sufficient for recovering~$\ti\ph$
as asymptotic expansion of~$\cL^\th\hat\ph$; more often, Borel-Laplace sums are
associated with sectors.

\med

\noindent
{\bf Note:} From the inversion of the Fourier transform, one can deduce a formula for
the {\em integral Borel transform} which allows one to recover~$\hat\ph(\ze)$
from~$\cL^\th\hat\ph(z)$.
For instance, $\hat\ph(\ze) =
\frac{1}{2\pi\I}\int_{\rho-\I\infty}^{\rho+\I\infty} \cL^0\hat\ph(z) \,
\ee^{z\ze} \, \dd z$ for small $\ze\ge0$, with suitable $\rho>0$.

\subsub{Sectorial sums}

\begin{figure}

\begin{center}

\epsfig{file=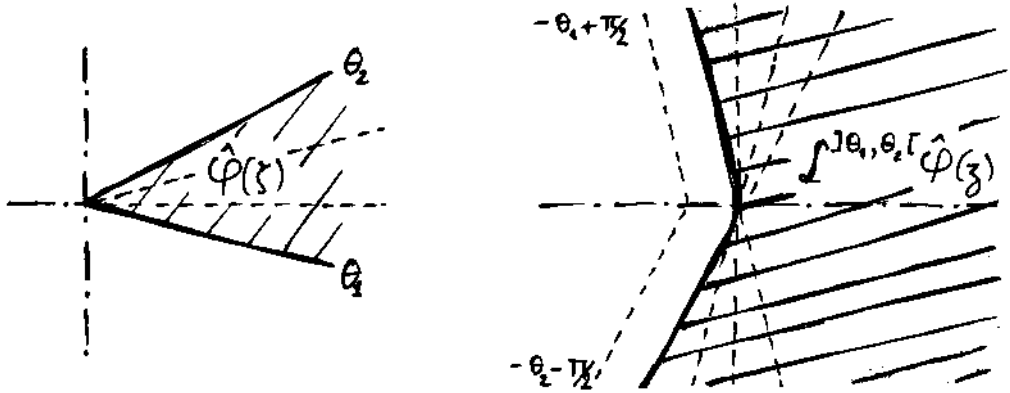,height=6cm,angle = 0}

\end{center}


\caption{\label{figsectors} Sectorial sums.}

\end{figure}

Suppose that $\hat\ph(\ze)$ converges near the origin and extends analytically
to a sector $\ao \ze\in\C \mid \th_1 < \arg\ze < \th_2 \af$ (where
$\th_1,\th_2\in\R$, $|\th_2-\th_1|<2\pi$), with exponential type less
than~$\tau$, then we can move the direction of integration~$\th$
inside~$\left]\th_1,\th_2\right[$.
According to
the Cauchy theorem, $\cL^{\th'}\hat\ph$ is the analytic continuation
of~$\cL^\th\hat\ph$ when $|\th'-\th|<\pi$, we can thus glue together these
holomorphic functions and obtain a function
$\cL^{\left]\th_1,\th_2\right[}\hat\ph$ analytic in the union of the half-planes 
$\ao\RE (z\,\ee^{\I\th})>\tau\af$, which is a sectorial neighbourhood of infinity
contained in $\ao -\th_2-\pi/2 < \arg z < -\th_1+\pi/2 \af$
(see Figure~\ref{figsectors}).
Notice however that, if $\th_2-\th_1>\pi$, the resulting function may be
multivalued, \ie one must consider the variable~$z$ as moving on the Riemann
surface of the logarithm.

\begin{figure}

\begin{center}

\epsfig{file=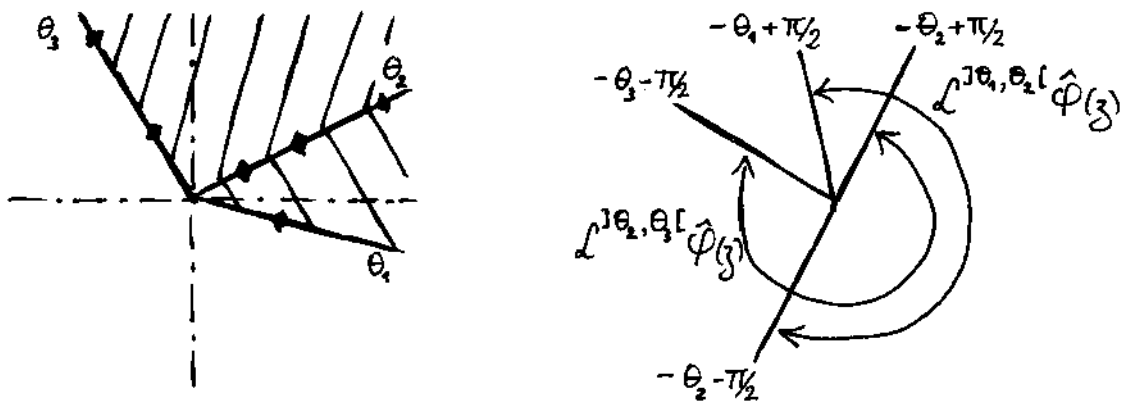,height=6cm,angle = 0}

\end{center}


\caption{\label{figfrequ} Several Borel-Laplace sums, analytic in different
domains, may be attached to a single divergent series.}

\end{figure}

A frequent situation is the following:
$\hat\ph=\cB\ti\ph$ converges and extends analytically to several infinite
sectors, with bounded exponential type, but also has singularities at finite
distance (in particular $\hat\ph$ has finite radius of convergence and~$\ti\ph$
is divergent).
Then several ``Borel-Laplace sums'' are available on various domains, but are
not the analytic continuations one of the other: the presence of singularities,
which separate the sectors one from the other, prevents one from applying the
Cauchy theorem.
On the other hand, all these ``sums'' share the same asymptotic expansion: the
mutual differences are exponentially small in the intersection of their domains
of definition (see Figure~\ref{figfrequ}).
%

\subsub{Resurgent functions}

It is interesting to ``measure'' the singularities in the~$\ze$-plane, since
they can be considered as responsible for the divergence of the common
asymptotic expansion~$\ti\ph(z)$ and for the exponentially small differences
between the various Borel-Laplace sums.
The resurgent functions can be defined as a class of formal series~$\ti\ph$
such that the analytic continuation of the formal Borel transform~$\hat\ph$
satisfies a certain condition regarding the possible singularities, which makes
it possible to develop a kind of singularity calculus (named ``alien
calculus''). These notions were introduced in the late 70s by J.~\'Ecalle, who
proved their relevance in a number of analytic problems~\cite{Eca81,Mal84}.
We shall not try to expound the theory in its full generality, but shall rather
content ourselves with explaining how it works in the case of certain difference
equations.

\med

\noindent
{\bf Note:} The formal Borel transform of a series~$\ti\ph(z)$ has positive
radius of convergence if and only if $\ti\ph(z)$ satisfies a ``Gevrey-$1$''
condition: $\hat\ph(\ze) \in \C\{\ze\} \Leftrightarrow \ti\ph(z) \in
z\ii{\C[[z\ii]]}_1$, where by definition
$$
\tst
z\ii{\C[[z\ii]]}_1 = \bigl\{\, \sum_{n\ge0} c_n z^{-n-1} \mid
\exists\, \rho>0
\text{\ such that\ } |c_n| = \cO(n!\rho^n) \,\bigr\}.
$$


\subsection{Linear and nonlinear difference equations}	\label{secdiffeq}


We shall be interested in formal series~$\ti\ph$ solutions of certain equations
involving the first-order difference operator $\ti\ph(z)\mapsto\ti\ph(z+1)-\ti\ph(z)$ (or
second-order differences). 
This operator is well defined in~$\C[[z\ii]]$, \eg by way of the Taylor formula
\begin{equation}	\label{eqdefdiff}
\ti\ph(z+1)-\ti\ph(z) = \pa\ti\ph(z) + \frac{1}{2!}\pa^2\ti\ph(z) 
  + \frac{1}{3!}\pa^3\ti\ph(z) + \cdots,
\end{equation}
where $\pa = \frac{\dd}{\dd z}$ and the series is formally convergent
because of increasing valuations
(we say that the series $\sum \frac{1}{r!}\pa^r\ti\ph$ is formally convergent
because the \rhs\ of~\refeq{eqdefdiff} is a well-defined formal series, each
coefficient of which is given by a finite sum of terms;
this is the notion of sequential convergence associated with the so-called Krull
topology).

It is elementary to compute the counterpart of the differential and difference
operators by~$\cB$:
$$
\cB\,:\; \pa\ti\ph(z) \mapsto -\ze\, \hat\ph(\ze),\quad\ens
\ti\ph(z+1) \mapsto \ee^{-\ze}\,\hat\ph(\ze).
$$
When~$\ti\ph(z)$ is obtained by solving an equation, a natural strategy is thus
to study~$\hat\ph(\ze)$ as solution of a transformed equation.
If a Laplace transform~$\cL^\th$ can be applied to~$\hat\ph$, one then recovers
an analytic solution of the original equation, because $\cL^\th\circ\cB$ commutes
with the differential and difference operators.

\subsub{Two linear equations}	\label{seclineardiffeq}

Let us illustrate this on two simple equations:
\begin{equation}	\label{eqnph}
\ti\ph(z+1) - \ti\ph(z) = a(z), \qquad 
\text{$a(z)\in z^{-2}\C\{z\ii\}$ given,}
\end{equation}
\begin{equation}	\label{eqnpsi}
\ti\psi(z+1) - 2\ti\psi(z) + \ti\psi(z-1) = b(z), \qquad 
\text{$b(z)\in z^{-3}\C\{z\ii\}$ given.}
\end{equation}
The corresponding equations for the formal Borel transforms are
$$
(\ee^{-\ze}-1)\hat\ph(\ze) = \hat a(\ze), \qquad 
\left(4\sinh^2\frac{\ze}{2}\right)\hat\psi(\ze) = \hat b(\ze).
$$
Here the power series $\hat a(\ze)$ and~$\hat b(\ze)$ converge to entire functions of
bounded exponential type in every direction, vanishing at the origin; moreover
$\hat b'(0)=0$.
We thus get in~$\C[[\ze]]$ unique solutions
$\hat\ph(\ze)=\hat a(\ze)/(\ee^{-\ze}-1)$ and
$\hat\psi(\ze)=\hat b(\ze)/\left(4\sinh^2\frac{\ze}{2}\right)$, 
which converge near the origin and define meromorphic functions, 
the possible poles being located in~$2\pi\I\,\Z^*$.

The original equations thus admit unique solutions $\ti\ph=\cB\ii\hat\ph$ and
$\ti\psi=\cB\ii\hat\psi$ in~$z\ii\C[[z\ii]]$. For each of them, Borel-Laplace
summation is possible and we get two natural sums, associated with two sectors:
$$
\ph^+(z) = \cL^\th\hat\ph (z), \quad 
\th\in\left]-\tfrac{\pi}{2},\tfrac{\pi}{2}\right[,
\qquad
\ph^-(z) = \cL^{\th'}\hat\ph (z), \quad 
\th'\in\left]\tfrac{\pi}{2},\tfrac{3\pi}{2}\right[,
$$
and similarly $\hat\psi(\ze)$ gives rise to $\psi^+(z)$ and~$\psi^-(z)$.

\begin{figure}

\begin{center}

\epsfig{file=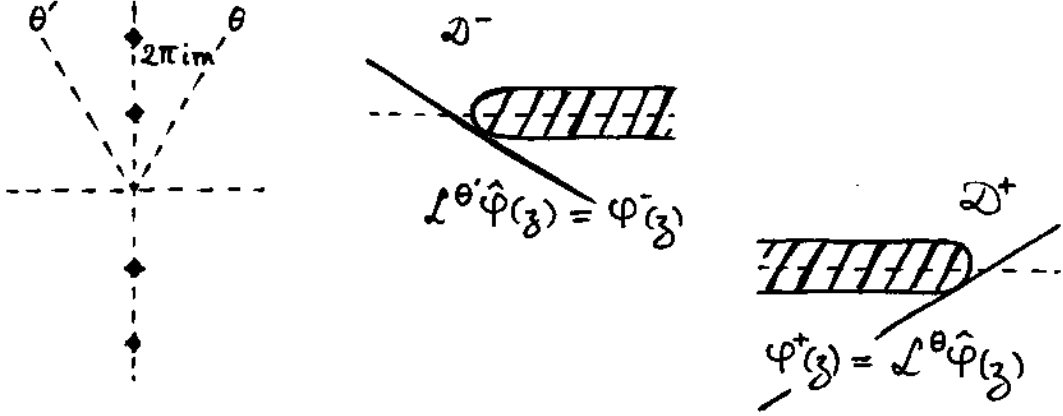,height=7cm,angle = 0}

\end{center}

\vspace{-.5cm}

\caption{\label{figlindiff} Borel-Laplace summation for the difference equation~\refeq{eqnph}.}

\end{figure}

The functions~$\ph^+$ and~$\psi^+$ are solutions of~\refeq{eqnph}
and~\refeq{eqnpsi}, analytic in a domain of the form 
$\cD^+ = \C \setminus\ao \dist(z,\R^-)\le\tau \af$.
The solutions~$\ph^-$ and~$\psi^-$ are defined in a symmetric domain~$\cD^-$
(see Figure~\ref{figlindiff}).
The intersection $\cD^+\cap\cD^-$ has two connected components,
$\ao \IM z< -\tau \af$ and $\ao \IM z> \tau \af$.
In the case of equation~\refeq{eqnph} for instance, the exponentially small
difference $\ph^+-\ph^-$ in the lower component is related to the singularities
of~$\hat\ph$ in~$2\pi\I\,\N^*$; it can be exactly computed by the resiuduum formula:
the singularity at $\om=2\pi\I m$ yields a contribution 
$$
A_\om\,\ee^{-\om z}, \qquad \text{with \;} A_\om = - 2\pi\I\, \hat a(\om)
$$
(the modulus of which is~$|A_\om|\ee^{2\pi m\IM z}$, which is exponentially small for
$\IM z\to-\infty$);
the difference
$(\ph^+-\ph^-)(z) = \int_{\ee^{\I\th'}\infty}^{\ee^{\I\th}\infty} 
\hat\ph(\ze)\,\ee^{-z\ze}\,\dd\ze$
is simply the sum of these contributions:
\begin{equation}	\label{eqlinearStokes}
\ph^+(z)-\ph^-(z) = \sum_{\om\in2\pi\I\,\N^*} A_\om\,\ee^{-\om z},
\qquad \IM z<-\tau
\end{equation}
as is easily seen by deforming the
contour of integration (choose $\th$ and~$\th'$ close enough to~$\pi/2$
according to the precise location of~$z$, and push the contour of integration upwards).
Symmetrically, the difference in the upper
component can be computed from the singularities in~$-2\pi\I\,\N^*$.

\sm

\noindent
{\bf Note:} $\ph^+(z)$ is the unique solution of~\refeq{eqnph} which tends to~$0$
when $\RE z\to+\infty$ and can be written as $-\sum_{k\ge0} a(z+k)$, and
$\ph^-(z) =\sum_{k\ge1} a(z-k)$ is the unique solution which tends to~$0$
when $\RE z\to-\infty$; the difference defines two $1$-periodic functions, the
Fourier coefficients of which can be expressed in term of the Fourier transform
of~$a(\pm\I\rho+z)$ (take $\rho>0$ large enough). One recovers the previous
formula for the difference by using the integral representation for the Borel
transform to compute the numbers~$\hat a(\om)$.

\subsub{Nonlinear equations} 

In the present text we shall show how one can deal with nonlinear difference equations
like
\begin{equation}	\label{eqniter}
\ti\ph(z+1) - \ti\ph(z) = a\bigl(z+\ti\ph(z)\bigr), \qquad 
\text{$a(z)\in z^{-2}\C\{z\ii\}$ given,}
\end{equation}
which is related to Abel's equation and the classification of holomorphic germs
in one complex variable, or
\begin{equation}	\label{eqnsecond}
\ti\psi(z+1) - 2\ti\psi(z) + \ti\psi(z-1) = b\bigl(\ti\psi(z),\ti\psi(z-1)\bigr),
\end{equation}
with certain $b(x,y)\in\C\{x,y\}$, which is related to splitting problems in two
complex variables.

\med

Dealing with nonlinear equations will require the study of {\em convolution},
which is the subject of sections~\ref{secAnContConv} and~\ref{secformconvmod}.
The Borel transforms~$\hat\ph(\ze)$ and~$\hat\psi(\ze)$ will still be
holomorphic at the origin but no longer meromorphic in~$\C$, as will be shown
later; their analytic continuations have more complicated singularities than
mere first- or second-order poles.
We shall introduce {\em alien calculus} in Section~\ref{secAlCalAbEq} and a more
general version of it in Section~\ref{secGenRes} to deal with these
singularities.


\subsection{The Riemann surface~$\cR$ and the analytic continuation of convolution}
\label{secAnContConv}


The first nonlinear operation to be studied is the multiplication of formal
series.

\begin{lemma}	\label{lemcounterpmult}
Let $\hat\ph$ and~$\hat\psi$ denote the formal Borel transforms of
$\ti\ph,\ti\psi\in z\ii\C[[z\ii]]$ and consider the product series
$\ti\chi=\ti\ph\ti\psi$. 
Then its formal Borel transform is given by the ``convolution''
\begin{equation}	\label{eqdefconvol}
(\cB\ti\chi)(\ze) =
(\hat\ph*\hat\psi)(\ze) = \int_0^\ze \hat\ph(\ze_1)\hat\psi(\ze-\ze_1)\,\dd\ze_1.
\end{equation}
\end{lemma}

The above formula must be interpreted termwise:
$\int_0^\ze \frac{\ze_1^n}{n!} \, \frac{(\ze-\ze_1)^m}{m!}\,\dd\ze_1
= \frac{\ze^{n+m+1}}{(n+m+1)!}$ (as can be checked \eg by induction on~$n$,
which is sufficient to prove the lemma).

\subsub{The problem of analytic continuation}

The formula can be given an analytic meaning in the case of Gevrey-$1$ formal
series: if $\hat\ph,\hat\psi\in\C\{\ze\}$, their convolution is convergent in
the intersection of the discs of convergence of~$\hat\ph$ and~$\hat\psi$ and
defines a new holomorphic germ~$\hat\ph*\hat\psi$ at the origin; 
formula~\refeq{eqdefconvol} then holds as a relation between holomorphic
functions, but only for~$|\ze|$ small enough (smaller than the radii of
convergence of~$\hat\ph$ and~$\hat\psi$).
What about the analytic continuation of~$\hat\ph*\hat\psi$ when~$\hat\ph$
and~$\hat\psi$ themselves admit an analytic continuation beyond their discs of
convergence?
What about the case when~$\hat\ph$ and~$\hat\psi$ extend to meromorphic
functions for instance?

\med

A preliminary answer is that~$\hat\ph*\hat\psi$ always admit an analytic
continuation in the intersection of the ``holomorphic stars'' of~$\hat\ph$
and~$\hat\psi$.
We define the holomorphic star of a germ as the union of
all the open sets~$U$ containing the origin in which it admits analytic
continuation and which are star-shaped \wrt\ the origin (\ie $\forall\ze\in U$,
$[0,\ze]\subset U$).
And it is indeed clear that if~$\hat\ph$ and~$\hat\psi$ are holomorphic in such
a~$U$, formula~\refeq{eqdefconvol} makes sense for all $\ze\in U$ and provides
the analytic continuation of~$\hat\ph*\hat\psi$.
With a view to further use we notice that, if $|\hat\ph(\ze)|\le\Phi(|\ze|)$ and
$|\hat\psi(\ze)|\le\Psi(|\ze|)$ for all $\ze\in U$, then
\begin{equation}	\label{ineqmajconv}
|\hat\ph*\hat\psi(\ze)| \le \Phi*\Psi(|\ze|), \qquad
\ze\in U.
\end{equation}

\med 

The next step is to study what happens on singular rays, behind singular points.
The idea is that convolution of poles generates ramification
(``multivaluedness'') but is easy to continue analytically.
For example, since
$$
1*\hat\ph(\ze) = \int_0^\ze \hat\ph(\ze_1)\,\dd\ze_1,
$$
we see that when~$\hat\ph$ is a meromorphic function with poles in a set
$\Om\subset\C^*$, $1*\hat\ph$ admits an analytic continuation along any path
issuing from the origin and avoiding~$\Om$; in other words, it defines a
function holomorphic on the universal cover\footnote{
\label{footClog}
Here it is understood that the base-point is at the origin.
If $\Om$ is a closed subset of~$\C$ with $\C\setminus\Om$ connected and
$\ze_0\in\C\setminus\Om$, the universal cover of~$\C\setminus\Om$ with
base-point~$\ze_0$ can be defined as the set of homotopy classes of paths
issuing from~$\ze_0$ and lying in~$\C\setminus\Om$ (equivalence classes for
homotopy with fixed extremities).
We denote it $\widetilde{(\C\setminus\Om,\ze_0)}$.
There is a covering map
$\pi:\widetilde{(\C\setminus\Om,\ze_0)}\to\C\setminus\Om$, which associates with
any class~$c$ the extremity~$\ga(1)$ of any path
$\ga:[0,1]\to\C\setminus\Om$ which represents~$c$,
and which allows one to define a Riemann surface structure
on~$\widetilde{(\C\setminus\Om,\ze_0)}$ by pulling back the complex structure
of~$\C\setminus\Om$ (see~\cite[pp.~81--89 and 105--112]{CNP}). 
For example, the Riemann surface of the logarithm is
$\widetilde{(\C\setminus\{0\},1)}$,
the points of which can be written ``$r\,\ee^{\I\th}$'' with $r>0$ and $\th\in\R$.
We often use the letter~$\ze$ for points of a universal cover, and then denote
by $\dze=\pi(\ze)$ their projection.
}
of~$\C\setminus\Om$, with
logarithmic singularities at the poles of~$\hat\ph$.

\med

But convolution may also create new singular points.
For instance, if $\hat\ph(\ze)=\frac{1}{\ze-\om'}$ and $\hat\psi(\ze) =
\frac{1}{\ze-\om''}$ with $\om',\om''\in\C^*$, one gets
$$
\hat\ph*\hat\psi(\ze) = \frac{1}{\ze-\om}\left(
\int_0^\ze \frac{\dd\ze_1}{\ze_1-\om'} + \int_0^\ze\frac{\dd\ze_1}{\ze_1-\om''} \right),
\qquad  \om=\om'+\om''.
$$
We thus have logarithmic singularities at~$\om'$ and~$\om''$, but also a pole
at~$\om$, the residuum of which is an integer multiple of~$2\pi\I$ which
depends on the path chosen to approach~$\om$. In other words,
$\hat\ph*\hat\psi$ extends meromorphically to the universal cover
of~$\C\setminus\{\om',\om''\}$, with a pole lying over~$\om$ (the residuum of which
depends on the sheet\footnote{
Again we can take the base-point at the origin to define the universal cover
of~$\C\setminus\Om$, here with $\Om=\{\om',\om''\}$. The word ``sheets'' usually
refers to the various lifts in the cover of an open subset~$U$ of the base space
which is star-shaped \wrt\ one of its points,
\ie to the various connected components of~$\pi\ii(U)$.
}
of the Riemann surface which is considered; in particular it
vanishes for the principal sheet\footnote{
In the case of a universal cover $\widetilde{(\C\setminus\Om,\ze_0)}$, the
``principal sheet''~$\ti U$ is obtained by considering the maximal open subset~$U$
of~$\C\setminus\Om$ which is star-shaped \wrt~$\ze_0$ and lifting it by means of
rectilinear segments: $\ti U$ is the set of all the classes of
segments~$[\ze_0,\ze]$, $\ze\in U$.
}
if $\arg\om'\not=\arg\om''$, which is consistent with
what was previously said on the holomorphic star).

\subsub{The Riemann surface~$\cR$}

With a view to the difference equations we are interested in and to the expected
behaviour of the Borel transforms, we define a Riemann surface which is obtained
by adding a point to the universal cover of $\C\setminus2\pi\I\,\Z$.

\begin{Def} 
Let~$\cR$ be the set of all homotopy classes of paths issuing from the origin and lying
inside~$\C\setminus 2\pi\I\,\Z$ (except for their initial point),
and let $\pi:\cR\to\C\setminus 2\pi\I\,\Z^*$ be the covering map, which
associates with any class~$c$ the extremity~$\ga(1)$ of any path
$\ga:[0,1]\to\C$ which represents~$c$.
We consider~$\cR$ as a Riemann surface by pulling back by~$\pi$ the complex structure
of~$\C\setminus 2\pi\I\,\Z^*$.
\end{Def}

Observe that~$\pi\ii(0)$ consists of only one point (the homotopy class of the
constant path), which we may call the origin of~$\cR$.
Let~$U$ be the complex plane deprived from the half-lines
$+2\pi\I\left[1,+\infty\right[$ and $-2\pi\I\left[1,+\infty\right[$.
We define the ``principal sheet'' of~$\cR$ as the set of all the classes of
segments~$[0,\ze]$, $\ze\in U$; equivalently, it is the connected component
of~$\pi\ii(U)$ which contains the origin.
We define the ``half-sheets'' of~$\cR$ as the various connected components
of~$\pi\ii(\{\RE\ze>0\})$ or of~$\pi\ii(\{\RE\ze<0\})$.

A holomorphic function of~$\cR$ can be viewed as a germ of holomorphic
function at the origin of~$\C$ which admits analytic continuation along any
path avoiding~$2\pi\I\,\Z$;
we then say that this germ ``extends holomorphically to~$\cR$''.
This definition a priori does not authorize analytic continuation along
a path which leads to the origin, unless this path stays in the principal
sheet\footnote{
That is, unless it lies in $U=\C\setminus\pm2\pi\I\left[1,+\infty\right[$.
We shall often identify the paths issuing from~$0$
in~$\C\setminus2\pi\I\,\Z$ and their lifts starting at the origin of~$\cR$.
Sometimes, we shall even identify a point of~$\cR$ with its projection by~$\pi$
(the path which leads to this point being understood), which amounts to treating
a holomorphic function of~$\cR$ as a multivalued function
on~$\C\setminus2\pi\I\,\Z$.
}.
More precisely, one can prove 
\begin{lemma}	\label{lemPhiph}
If $\Phi$ is holomorphic in~$\cR$, then its restriction to the principal sheet
defines a holomorphic function~$\ph$ of~$U$ which extends analytically along any
path~$\ga$ issuing from~$0$ and lying in~$\C\setminus2\pi\I\,\Z$. The analytic
continuation is given by $\ph(\ga(t))=\Phi(\Ga(t))$, where $\Ga$ is the lift
of~$\ga$ which starts at the origin of~$\cR$.

Conversely, given $\ph\in\C\{\ze\}$, if any $c\in\cR$ can be represented by a
path of analytic continuation for~$\ph$, then the value of~$\ph$ at the
extremity~$\ga(1)$ of this path depends only on~$c$ and the formula
$\Phi(c)=\ph(\ga(1))$ defines a holomorphic function of~$\cR$.
\end{lemma}

The absence of singularity at the origin on the principal sheet is the only
difference between~$\cR$ and the universal cover of $\C\setminus2\pi\I\,\Z$
with base-point at~$1$.
For instance, among the two series 
$$
\dst \sum_{m\in\Z^*} \frac{1}{\ze} \,\ee^{-|m|} \int_1^\ze 
\frac{\dd\ze_1}{\ze_1-2\pi\I m}, \quad
\dst \sum_{m\in\Z^*} \frac{1}{\ze} \,\ee^{-|m|} \int_0^\ze 
\frac{\dd\ze_1}{\ze_1-2\pi\I m},
$$ 
the first one defines a function which is holomorphic in the universal cover of
$\C\setminus2\pi\I\,\Z$ but not in~$\cR$, whereas the second one defines a
holomorphic function of~$\cR$.

\subsub{Analytic continuation of convolution in~$\cR$}

The main result of this section is

\begin{thm}	\label{thmstbconv}
If two germs at the origin extend holomorphically to~$\cR$, so does their
convolution product.
\end{thm}

\idproof
Let $\hat\ph$ and~$\hat\psi$ be holomorphic germs at the origin of~$\C$ which
admit analytic continuation along any path avoiding~$2\pi\I\,\Z$; we denote by
the same symbols the corresponding holomorphic functions of~$\cR$.
One could be tempted to think that, if a point~$\ze$ of~$\cR$ is defined by a
path~$\ga$, the integral
\begin{equation}	\label{eqtempt}
\hat\chi(\ze) = \int_{\ga} \hat\ph(\ze') \hat\psi(\ze-\ze') \,\dd\ze'
\end{equation}
would give the value of the analytic continuation of $\hat\ph*\hat\psi$
at~$\ze$. However, this formula does not always make sense, since
one must worry about the path~$\ga'$ followed by~$\ze-\ze'$ when~$\ze'$
follows~$\ga$: is~$\hat\psi$ defined on this path? 
In fact, even if~$\ga'$ lies in~$\C\setminus2\pi\I\,\Z$ (and thus
$\hat\psi(\ze-\ze')$ makes sense), even if~$\ga'$ coincides with~$\ga$,
it may happen that this integral does not give the analytic continuation of
$\hat\ph*\hat\psi$ at~$\ze$ (usually, the value of this integral does not depend
only on~$\ze$ but also on the path~$\ga$).\footnote{
\label{footmixedconv}
However, if~$\hat\psi$ is entire, it is true that the integral~\refeq{eqtempt}
does provide the analytic continuation of~$\hat\ph*\hat\psi$ along~$\ga$.
}

The construction of the desired analytic continuation relies on the idea of
``symmetrically contractile'' paths.
A path~$\ga$ issuing from~$0$ is said to be {\em $\cR$-symmetric} if it lies
in~$\C\setminus2\pi\I\,\Z$ (except for its starting point) and is symmetric
\wrt\ its midpoint: the paths $t\in[0,1]\mapsto\ga(1)-\ga(t)$ and
$t\in[0,1]\mapsto\ga(1-t)$ coincide up to reparametrisation.
A path is said to be {\em $\cR$-symmetrically contractile} if it is $\cR$-symmetric
and can be continuously deformed and shrunk to~$\{0\}$ within the class of
$\cR$-symmetric paths.
The main point is that any point of~$\cR$ can be defined by an
$\cR$-symmetrically contractile path. More precisely:
\begin{lemma}	\label{lemRsymcontr}
Let $\ga$ be a path issuing from~$0$ and lying in~$\C\setminus2\pi\I\,\Z$
(except for its starting point).
Then there exists an $\cR$-symmetrically contractile path~$\Ga$ which is
homotopic to~$\ga$.
Moreover, one can construct~$\Ga$ so
that there is a continuous map $(s,t)\mapsto H(s,t) = H_s(t)$ satisfying
\begin{enumerate}
\item $H_0(t)\equiv 0$ and $H_1(t)\equiv \Ga(t)$,
\item each~$H_s$ is an $\cR$-symmetric path with $H_s(0)=0$ and $H_s(1)=\ga(s)$.
\end{enumerate}
\end{lemma}

\begin{figure}

\begin{center}

\epsfig{file=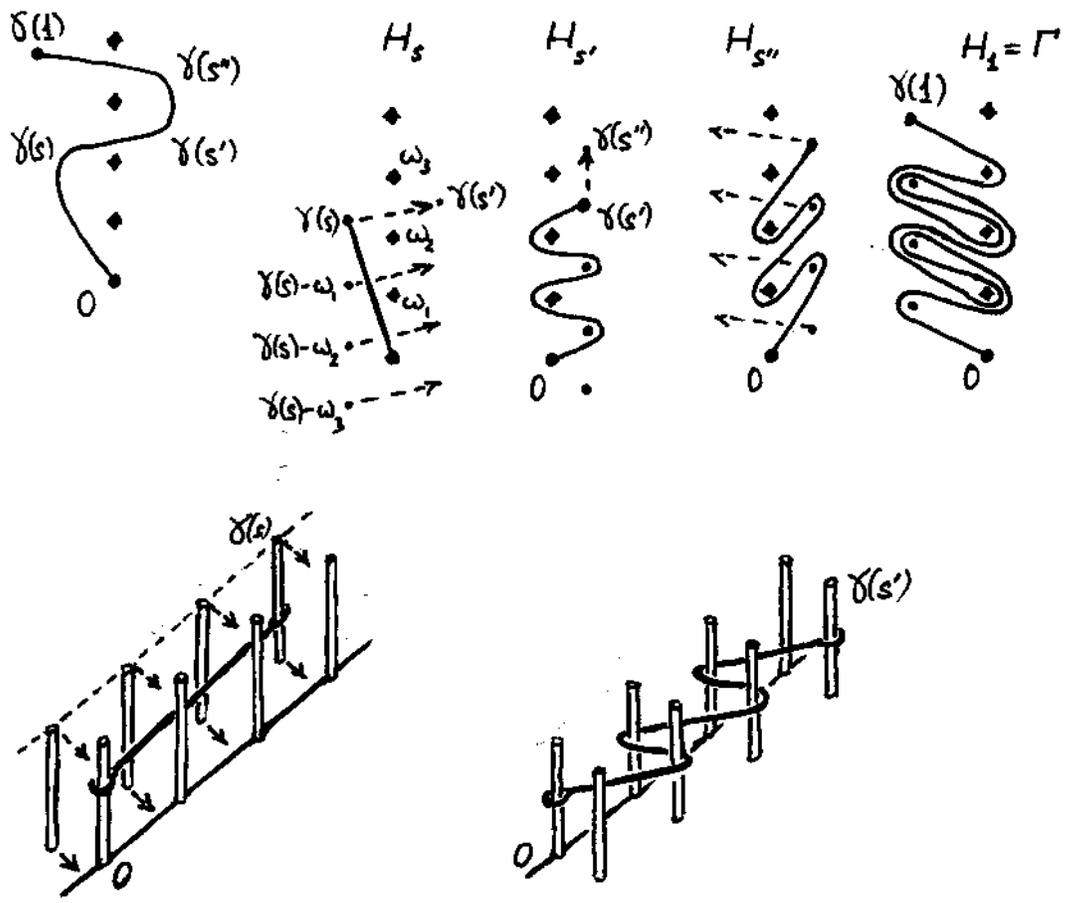,height=12cm,angle = 0}

\end{center}


\caption{\label{fignails} Construction of an $\cR$-symmetrically contractile path~$\Ga$
homotopic to~$\ga$.}

\end{figure}

We shall not try to write a formal proof of this lemma, but it is easy to
visualize a way of constructing~$H$.
Let a point $\ze=\ga(s)$ move along~$\ga$ (as $s$ varies from~$0$ to~$1$) and
remain connected to~$0$ by an extensible thread, 
with moving nails pointing downwards at each point of~$\ze-2\pi\I\,\Z$,
while fixed nails point upwards at each point of~$2\pi\I\,\Z$
(imagine for instance that the first nails are fastened to a moving rule and the
last ones to a fixed rule).
As $s$ varies, the thread is progressively stretched but it has to meander
between the nails.
The path~$\Ga$ is given by the thread in its final form, when~$\ze$ has reached
the extremity of~$\ga$; the paths~$H_s$ correspond to the thread at intermediary
stages\footnote{
Note that the mere existence of a continuous~$H$ satisfying conditions~{\em i)}
and~{\em ii)} implies that~$\ga$ and~$\Ga$ are homotopic, as is visually clear
(the formula 
$$
h_\la(t) = H\Bigl(\la+(1-\la)t, \frac{t}{\la+(1-\la)t}\Bigr),
\qquad 0\le\la\le1
$$
yields an explicit homotopy).
}
(see Figure~\ref{fignails}).

It is now easy to end the proof of Theorem~\ref{thmstbconv}.
Given $\hat\ph$, $\hat\psi$ as above and $\ga$ a path of~$\cR$ along which
we wish to follow the analytic continuation of~$\hat\ph*\hat\psi$,
we take~$H$ as in Lemma~\ref{lemRsymcontr} and let the reader convince himself
that the formula
\begin{equation}	\label{eqancontconv}
\hat\chi(\ze) = \int_{H_s} \hat\ph(\ze') \hat\psi(\ze-\ze') \,\dd\ze',
\qquad \ze = \ga(s),
\end{equation}
defines the analytic continuation~$\hat\chi$ of~$\hat\ph*\hat\psi$
along~$\ga$
(in this formula, $\ze'$ and $\ze-\ze'$ move on the same path~$H_s$ which avoids
$2\pi\I\,\Z$, by $\cR$-symmetry).
See \cite[Vol.~1]{Eca81}, \cite{CNP}, \cite{sauzin01} for more details.
\eopf

\sm

Of course, if the path~$\ga$ mentioned in the last part of the proof stays in
the principal sheet of~$\cR$, the analytic continuation is simply given by
formula~\refeq{eqdefconvol}. 
In particular, if~$\hat\ph$ and~$\hat\psi$ have bounded exponential type in a
direction $\arg\ze=\th$, $\th\notin\frac{\pi}{2}+\pi\Z$, it follows from
inequality~\refeq{ineqmajconv} that $\hat\ph*\hat\psi$ has the same property.


\subsection{Formal and convolutive models of the algebra of resurgent functions,
$\ti{\protect\cH}$ and~$\wh{\protect\cH}(\cR)$}
\label{secformconvmod}


In view of Theorem~\ref{thmstbconv}, the convolution of germs induces an
internal law on the space of holomorphic functions of~$\cR$, which is
commutative and associative (being the counterpart of multiplication of formal
series, by Lemma~\ref{lemcounterpmult}).
In fact, we have a commutative algebra (without unit), which can be viewed as a
subalgebra of the convolution algebra~$\C\{\ze\}$, and which corresponds
via~$\cB$ to a subalgebra (for the ordinary product of formal series)
of~$z\ii\C[[z\ii]]$.

\begin{Def}
The space~$\wh\cH(\cR)$ of all holomorphic functions of~$\cR$, equipped with the
convolution product, is an algebra called the convolutive model of the algebra
of resurgent functions.
The subalgebra $\ti\cH=\cB\ii\bigl(\wh\cH(\cR)\bigr)$ of~$z\ii\C[[z\ii]]$ is
called the multiplicative model of the algebra of resurgent functions.
\end{Def}

The formal series in~$\ti\cH$ (most of which have zero radius of convergence)
are called ``resurgent functions''. These definitions will in fact be extended
to more general objects in the following (see Section~\ref{secSingul} on
``singularities'').

\med

There is no unit for the convolution in~$\wh\cH(\cR)$. Introducing a new symbol
$\de=\cB 1$, we extend the formal Borel transform:
$$
\cB\,:\;
\ti\chi(z) = c_0 + \sum_{n\ge0} c_n z^{-n-1} \in \C[[z\ii]]
\ens\mapsto\ens
\hat\chi(\ze) = c_0\,\de + \sum_{n\ge0} c_n \frac{\ze^n}{n!} \in \C\,\de\oplus\C[[\ze]],
$$
and also extend convolution from~$\C[[\ze]]$ to~$\C\,\de\oplus\C[[\ze]]$
linearly, by treating~$\de$ as a unit (\ie so as to keep~$\cB$ a morphism of algebras).
This way, $\C\,\de\oplus\wh\cH(\cR)$ is an algebra for the convolution, which is
isomorphic via~$\cB$ to the algebra $\C\oplus\ti\cH$.
Observe that $$\C\{z\ii\} \subset \C\oplus\ti\cH \subset {\C[[z\ii]]}_1.$$

\med

Having dealt with multiplication of formal series, we can study {\em composition}
and its image in $\C\,\de\oplus\wh\cH(\cR)$:
\begin{prop}	\label{propcompos}
Let $\ti\chi\in\C\oplus\ti\cH$. Then composition by $z\mapsto z+\ti\chi(z)$ defines a
linear operator of~$\C\oplus\ti\cH$ into itself, and for any $\ti\psi\in\ti\cH$ the
Borel transform of
$\ti\al(z) = \ti\psi(z+\ti\chi(z))= \sum_{r\ge0} \frac{1}{r!}
\pa^r\ti\psi(z) \ti\chi^r(z)$ is given by the series of functions
\begin{equation}	\label{eqcomposconvol}
\hat\al(\ze) = \sum_{r\ge0} \frac{1}{r!}
\left((-\ze)^r\hat\psi(\ze)\right) * \hat\chi^{*r}(\ze)
\end{equation}
(where $\hat\chi=\cB\ti\chi$ and $\hat\psi=\cB\ti\psi$),
which is uniformly convergent in every compact subset of~$\cR$.
\end{prop}

The convergence of the series stems from the regularizing character
of convolution (the convergence in the principal sheet of~$\cR$ can be proved by use
of~\refeq{ineqmajconv}; see~\cite[Vol.~1]{Eca81} or~\cite{CNP} for the
convergence in the whole Riemann surface).

\med

The notation $\hat\al = \hat\psi \circledast (\de'+\hat\chi)$ and the name
``composition-convolution'' are used in~\cite{Mal},
with a symbol $\de'=\cB z$ which must be considered as the derivative of~$\de$.
%
The symbols~$\de$ and~$\de'$ will be interpreted as elementary singularities in
Section~\ref{secSingul}.

\med

In Proposition~\ref{propcompos}, the operator of composition by $z\mapsto
z+\ti\chi(z)$ is invertible;
in fact, $\ID+\ti\chi$ has a well-defined inverse for composition in
$\ID+\C[[z\ii]]$, which turns out to be also resurgent:

\begin{prop}	\label{propInversion}
If $\ti\chi\in\C\oplus\ti\cH$, the formal transformation $\ID+\ti\chi$ has an
inverse (for composition) of the form $\ID+\ti\ph$ with $\ti\ph\in\ti\cH$.
\end{prop}

This can be proven by the same arguments as Proposition~\ref{propcompos},
since the Lagrange inversion formula allows one to write
\begin{equation}	\label{eqLagr}
\ti\ph = \sum_{k\ge1} \frac{(-1)^k}{k!} \pa^{k-1}\bigl(\ti\chi^k\bigr),
\quad\text{hence}\ens
\hat\ph = - \sum_{k\ge1} \frac{\ze^{k-1}}{k!}\hat\chi^{*k}.
\end{equation}
One can thus think of $z\mapsto z+\ti\chi(z)$ as of a ``resurgent change of variable''.

\med

Similarly, {\em substitution} of a resurgent function without constant term into a
convergent series is possible:
\begin{prop}	\label{propsubs}
If $C(w)=\sum_{n\ge0}C_n w^n \in\C\{w\}$ and $\ti\psi\in\ti\cH$,
then the formal series
$
C\circ\ti\psi(z) = \sum_{n\ge0} C_n \ti\psi^n(z)
$
belongs to $\C\oplus\ti\cH$.
\end{prop}
The proof consists in verifying the convergence of the series
$\cB(C\circ\ti\psi) = \sum_{n\ge0} C_n \hat\psi^{*n}$.

\med

As a consequence, any resurgent function with nonzero constant term has a
resurgent multiplicative inverse:
$1/(c+\ti\psi) = \sum_{n\ge0}(-1)^n c^{-n-1}\ti\psi^n \in\C\oplus\ti\cH$.
The exponential of a resurgent function~$\ti\psi$ is also a resurgent function, the Borel
transform of which is the {\em convolutive exponential}
$$
\exp_*(\hat\psi) = \de + \hat\psi + \frac{1}{2!}\hat\psi*\hat\psi +
\frac{1}{3!}\hat\psi*\hat\psi*\hat\psi + \ldots
$$
(in this case the substitution is well-defined even if $\ti\psi(z)$ has a
constant term).

\med

We end this section by remarking that the role of the lattice~$2\pi\I\,\Z$ in
the definition of~$\cR$ is not essential in the theory of resurgent functions.
See Section~\ref{secGenRes} for a more general definition of the space of
resurgent functions (in which the location of singular points is not a priori
restricted to~$2\pi\I\,\Z$), with a property of stability by convolution as in
Theorem~\ref{thmstbconv}, and with alien derivations more general than the ones
to be defined in Section~\ref{secAlienDer}.


\section{Alien calculus and Abel's equation}	\label{secAlCalAbEq}


We now turn to the resurgent treatment of the nonlinear first-order difference
equation~\refeq{eqniter}, beginning with a few words of motivation.

\subsection{Abel's equation and tangent-to-identity holomorphic germs
of~$(\C,0)$}
\label{secAbEqtgtId}


One of the origins of \'Ecalle's work on Resurgence theory is the problem of the
classification of holomorphic germs~$F$ of~$(\C,0)$ in the ``resonant'' case.
This is the question, important for one-dimensional complex dynamics, of
describing the conjugacy classes of the group~$\G$ of local analytic
transformations $w\mapsto F(w)$ which are locally invertible, \ie of the form
$F(w)=\la w+\cO(w^2)\in\C\{w\}$ with $\la\in\C^*$.
It is well-known that, if the multiplier $\la=F'(0)$ has modulus $\neq1$, then
$F$ is holomorphically linearizable: there exists $H\in \G$ such that $H\ii\circ
F\circ H(w)=\la w$.
Resurgence comes into play when we consider the resonant case, \ie when~$F'(0)$
is a root of unity (the so-called ``small divisor problems'', which appear
when~$F'(0)$ has modulus~$1$ but is not a root of unity, are of different
nature---see S.~Marmi's lecture in this volume).

\sm

The references for this part of the text are: \cite[Vol.~2]{Eca81}, \cite{Eca84},
\cite{Mal84} (and Example~1 of~[Eca05] p.~235).
For non-resurgent approaches of the same problem, see \cite{MarRam},
\cite{DouadyHubb}, \cite{Shi}, \cite{Shi00}, \cite{Milnor}, \cite{Loray}.

\subsub{Non-degenerate parabolic germs}

Here, for simplicity, we limit ourselves to $F'(0)=1$, \ie to germs~$F$ which
are tangent to identity, with the further requirement that $F''(0)\neq0$, a
condition which is easily seen to be invariant by conjugacy.
Rescaling the variable~$w$ if necessary, one can suppose $F''(0)=2$.
It will be more convenient to work ``near infinity'', \ie to use the variable $z=-1/w$.

\begin{Def}
We call ``non-degenerate parabolic germ at the origin'' any $F(w)\in\C\{w\}$ of
the form
$$
F(w) = w + w^2 + \cO(w^3).
$$
We call ``non-degenerate parabolic germ at infinity'' a transformation $z\mapsto
f(z)$ which is conjugated by $z=-1/w$ to a non-degenerate parabolic germ~$F$ at the
origin:
$$
f(z) = -1/\bigl(F(-1/z)\bigr),
$$
\ie any $f(z) = z + 1 + a(z)$ with $a(z) \in z\ii\C\{z\ii\}$.
\end{Def}

Let $\G_1$ denote the subgroup of tangent-to-identity germs.
One can easily check that, if $F,G\in\G_1$ and $H\in\G$,
then $G=H\ii\circ F\circ H$ implies $G''(0)=H'(0)F''(0)$.
In order to work with non-degenerate parabolic germs only, we can thus restrict
ourselves to tangent-to-identity conjugating transformations~$H$, \ie we can
content ourselves with studying the adjoint action of~$\G_1$.

\sm

It turns out that formal transformations also play a role. 
Let $\ti\G_1$ denote the group (for composition) of formal series of the form
$\ti H(w)= w+\cO(w^2) \in \C[[w]]$.
It may happen that two parabolic germs~$F$ and~$G$ are conjugated by such a
formal series~$\ti H$, \ie $G=H\ii\circ F\circ H$ in~$\ti\G_1$, without being
conjugated by any convergent series: the $\G_1$-conjugacy classes we are interested in form
a finer partition than the ``formal conjugacy classes''.

\sm

In fact, the formal conjugacy classes are easy to describe.
One can check that, for any two non-degenerate parabolic germs
$F(w),G(w)=w+w^2+\cO(w^3)$, there exists $\ti H\in\ti\G_1$ such that
$G=H\ii\circ F\circ H$ if and only if the coefficient of~$w^3$ is the same
in~$F(w)$ and~$G(w)$.
In the following, this coefficient will usually be denoted~$\al$.

\sm

Let us rephrase the problem at infinity, using the variable $z=-1/w$, and thus
dealing with transformations belonging to $\ID+\C[[z\ii]]$.
The formula $h(z)=-1/H(-1/z)$ puts in correspondence the conjugating
transformations~$H$ of~$\G_1$ or~$\ti\G_1$ and the series of the form
\begin{equation}	\label{eqhb}
h(z) = z + b(z), \qquad b(z) \in \C\{z\ii\} \text{ or } b(z) \in \C[[z\ii]].
\end{equation}
Given a non-degenerate parabolic germ at infinity $f(z) =
-1/\bigl(F(-1/z)\bigr)$,  the coefficient~$\al$ of~$w^3$ in~$F(w)$ shows up in
the coefficient of~$z^{-1}$ in~$f(z)$:
\begin{equation}	\label{eqfa}
f(z) = z + 1 + a(z), \qquad a(z) = (1-\al)z\ii + \cO(z^{-2}) \in \C\{z\ii\}.
\end{equation}
The coefficient $\rho=\al-1$ is called ``r\'esidu it\'eratif'' in \'Ecalle's work,
or ``resiter'' for short.
Thus {\em any two germs of the form~\refeq{eqfa} are conjugated by a formal
transformation of the form~\refeq{eqhb} if and only if they have the same
resiter.}

\subsub{The related difference equations}

The simplest formal conjugacy class is the one corresponding to $\rho=0$. Any
non-degenerate parabolic germ $f$ or~$F$ with vanishing resiter is conjugated by
a formal~$h$ or~$H$ to $z\mapsto f_0(z) = z+1$ or $w\mapsto F_0(w) = \frac{w}{1-w}$.
We can be slightly more specific:

\begin{prop}	\label{propformalconj} 
Let $f(z) = z + 1 + a(z)$ be a non-degenerate parabolic germ at infinity with
vanishing resiter, \ie $a(z)\in z^{-2}\C\{z\ii\}$.
Then there exists a unique $\ti\ph(z) \in z\ii\C[[z\ii]]$ such that the formal
transformation $\ti u=\ID+\ti\ph$ is solution of
\begin{equation}	\label{eqnAbelInv}
u\ii\circ f\circ u(z) = z+1.
\end{equation}
The inverse formal transformation $\ti v=\ti u\ii$ is the unique transformation
of the form $\ti v=z+\ti\psi(z)$, with $\ti\psi(z) \in z\ii\C[[z\ii]]$, solution of
\begin{equation}	\label{eqnAbel}
v\bigl(f(z)\bigr) = v(z)+1.
\end{equation}
All the other formal solutions of equations~\refeq{eqnAbelInv}
and~\refeq{eqnAbel} can be deduced from~$\ti u$ and~$\ti v$: they are the series
\begin{equation}	\label{eqnuvphpsi}
u(z) = z + c + \ti\ph(z +c),
\quad
v(z) = z - c +\ti\psi(z),
\end{equation}
where $c$ is an arbitrary complex number.
\end{prop}

We omit the proof of this proposition, which can be done by substitution of an
indetermined series $u\in\ID+\C[[z\ii]]$ in~\refeq{eqnAbelInv}.
Setting $v=u\ii$, the $u$-equation then translates into
equation~\refeq{eqnAbel},
 as illustrated on the commutative diagram
$$
\xymatrix{
{\rule[-.6ex]{0ex}{0ex} \makebox[1em]{$z$}} \ar@<-.5ex>[d]_*+{u} \ar[rr] & 
& {\rule[-.6ex]{0ex}{0ex} \makebox[3em]{$z+1$}} \ar@<-.5ex>[d]_*+{u}  \\
{\rule[-1ex]{0ex}{3ex} \makebox[1em]{$z$}} \ar@<-.5ex>[u]_*+{v} 
%
\ar[dr];[]^*+{ \raisebox{-.5ex}[0ex][0ex]{$z=-1/w$} } \ar[rr] & 
& {\rule[-1ex]{0ex}{3ex} \makebox[2em]{$f(z)$}} \ar[dr];[] \ar@<-.5ex>[u]_*+{v} \\
& {w} \ar[rr] & & {F(w)}.
}
$$

\med

Notice that, under the change of unknown $u(z)=z+\ph(z)$, the conjugacy
equation~\refeq{eqnAbelInv} is equivalent to the equation
$$
\ph(z) + a(z+\ph(z)) = \ph(z+1),
$$
\ie to the difference equation~\refeq{eqniter} with $a(z)=f(z)-z-1$.
Equation~\refeq{eqnAbel} is called {\em Abel's equation}\footnote{
In fact, this name usually refers to the equation $V\circ F=V+1$, for
$V(w)=v(-1/w)$, which expresses the conjugacy by $w\mapsto V(w)=-1/w+\cO(w)$
between the given germ~$F$ at the origin and the normal form~$f_0$ at infinity. 
}. 

\sm

The formal solutions~$\ti u$ and~$\ti v$ mentioned in
Proposition~\ref{propformalconj} are generically divergent.
It turns out that they are always resurgent.
Before trying to explain this, let us mention that the case of a general
resiter~$\rho$ can be handled by studying the same equations~\refeq{eqnAbelInv}
and~\refeq{eqnAbel}: if $\rho\neq0$, there is no solution in~$\ID+\C[[z\ii]]$,
but one finds a unique formal solution of Abel's equation of the form 
$$
\ti v(z) = z + \ti\psi(z), \qquad
\ti\psi(z) = \rho\log z + \sum_{n\ge1} c_n z^{-n},
$$ 
the inverse of which is of the form
$$
\ti u(z) = z + \ti\ph(z), \qquad
\ti\ph(z) = -\rho\log z + \sum_{\substack{n,m\ge0\\ n+m\ge1}} C_{n,m} z^{-n} (z\ii\log z)^m,
$$
and these series~$\ti\psi$ and~$\ti\ph$ can be treated by Resurgence almost as
easily as the corresponding series in the case $\rho=0$.
In \'Ecalle's work, the formal solution~$\ti v$ of Abel's equation is called the
{\em iterator} (it\'erateur, in French) of~$f$ and its inverse~$\ti u$ is called
the {\em inverse iterator} because of their role in iteration theory (which we shall
not develop in this text---see however footnote~\ref{footiter}).

\subsub{Resurgence in the case $\rho=0$}

From now on we focus on the case $\rho=0$,  thus with
``formal normal forms'' $f_0(z)=z+1$ at infinity or 
$F_0(w)=\frac{w}{1-w}$ at the origin.
We do not intend to give the complete resurgent solution of the classification
problem, but only to convey some of the ideas used in \'Ecalle's approach.

\begin{thm}	\label{thmresur}
In the case $\rho=0$ (vanishing resiter), the formal series
$\ti\ph(z),\ti\psi(z)\in z\ii\C[[z\ii]]$ in terms of which the solutions of
equations~\refeq{eqnAbelInv} and~\refeq{eqnAbel} can be expressed as
in~\refeq{eqnuvphpsi} have formal Borel transforms~$\hat\ph(\ze)$
and~$\hat\psi(\ze)$ which converge near the origin and extend holomorphically
to~$\cR$, with at most exponential growth in the directions $\arg\ze=\th$,
$\th\notin\frac{\pi}{2}+\pi\Z$
(for every $\th_0\in\left]0,\frac{\pi}{2}\right[$, there exists $\tau>0$ such
that $\hat\ph$ and~$\hat\psi$ have exponential type $\le\tau$ in the sectors
$\ao -\th_0+n\pi \le \arg\ze \le \th_0+n\pi \af$, $n=0$ or~$1$).
\end{thm}

In other words, Abel's equation gives rise to resurgent functions and it is
possible to apply the Borel-Laplace summation process to~$\ti\ph$ and~$\ti\psi$.

\med

\noindent {\em Idea of the proof.}
Equation~\refeq{eqnAbel} for $v(z)=z+\psi(z)$ translates into
\begin{equation}	\label{eqnAbelpsi}
\psi\bigl(z+1+a(z)\bigr) - \psi(z) = - a(z).
\end{equation}
The proof indicated in~\cite[Vol.~2]{Eca81} or~\cite{Mal84} relies on the
expression of the unique solution in~$z\ii\C[[z\ii]]$ as an infinite sum of
iterated operators applied to~$a(z)$; the formal Borel transform then yields a
sum of holomorphic functions which is uniformly convergent on every compact
subset of~$\cR$.
One can prove in this way that $\hat\psi\in\wh\cH(\cR)$ with at most exponential
growth at infinity, and then deduce from Proposition~\ref{propInversion} and
formula~\refeq{eqLagr} that~$\hat\ph$ has the same property.

\sm

Let us outline an alternative proof, which makes use of
equation~\refeq{eqnAbelInv} to prove that $\hat\ph\in\wh\cH(\cR)$.
As already mentioned, the change of unkwnown $u=\ID+\ph$ leads to
equation~\refeq{eqniter} with $a(z)\in z^{-2}\C\{z\ii\}$, which we now treat as
a perturbation of equation~\refeq{eqnph}: we write it as
$$
\ph(z+1)-\ph(z) = a_0(z) + \sum_{r\ge1} a_r(z) \ph^r(z),
$$
with $a_r = \frac{1}{r!} \pa^r a$.
The unique formal solution without constant term,~$\ti\ph$, has a formal Borel
transform~$\hat\ph$ which thus satisfies
\begin{equation}	\label{eqEqFBTph}
\hat\ph = E\,\hat a_0 + E\,\sum_{r\ge1} \hat a_r*\hat\ph^{*r},
\end{equation}
where $\dst E(\ze) = \frac{1}{\ee^{-\ze}-1}$ and 
$\hat a_r(\ze) = \frac{1}{r!} (-\ze)^r \hat a(\ze)$,
$\hat a = \cB\, a$.

The convergence of~$\hat\ph$ and its analytic extension to the principal sheet
of~$\cR$ are easily obtained:
we have $\hat\ph = \sum_{n\ge1} \hat\ph_n$ with
\begin{equation}	\label{eqseriesfcns}
\hat\ph_1 = E\,\hat a_0, \qquad
\hat\ph_n = E\sum_{ \genfrac{}{}{0pt}{}{r\ge1}{n_1+\cdots+n_r=n-1} }
\hat a_r * \hat\ph_{n_1} * \cdots * \hat\ph_{n_r},
\qquad n\ge2
\end{equation}
(more generally $\ti u(z) = z + \sum_{n\ge1}\eps^n\ti\ph_n$ is the solution
corresponding to $f(z) = z + 1 + \eps a(z)$).
Observe that this series is well-defined and formally convergent, because $\hat
a \in \ze\C[[\ze]]$ and $E\in\ze\ii\C[[\ze]]$ imply $\hat\ph_n \in
\ze^{2(n-1)}\C[[\ze]]$,
and that each~$\hat\ph_n$ is convergent and extends holomorphically to~$\cR$ (by
virtue of Theorem~\ref{thmstbconv}, because $\hat a$ converges to an entire
function and $E$ is meromorphic with poles in~$2\pi\I\,\Z$);
we shall check that the series of functions $\sum \hat\ph_n$ is uniformly
convergent in every compact subset of the principal sheet.
Since $a(z)\in z^{-2}\C\{z\ii\}$, we can find positive constants~$C$
and~$\tau$ such that
$$
\abs{ \hat a(\ze) } \le C \min\bigl(1,\abs\ze\bigr)\,\ee^{\tau|\ze|},
\qquad \ze\in\C.
$$
Identifying the principal sheet of~$\cR$ with the cut plane
$\C\setminus\pm2\pi\I\left[1,+\infty\right[$, we can write it as the union over
$c>0$ of the sets 
$\cR_c\ss0 = \ao\ze\in\C \mid \dist([0,\ze],\pm2\pi\I)\ge c\af$ 
(with $c<2\pi$).
For each~$c>0$, we can find $\la=\la(c)>0$ such that 
$$
\abs{E(\ze)} \le \la \bigl( 1 + \abs\ze\ii \bigr), 
\qquad \ze\in\cR_c\ss0.
$$
We deduce that $\abs{\hat\ph_1(\ze)} \le 2\la C\,\ee^{\tau\abs\ze}$
in~$\cR_c\ss0$, and the fact that~$\cR_c\ss0$ is star-shaped \wrt\ the origin
allows us to construct majorants by inductive use of~\refeq{ineqmajconv}:
$$
\abs{\hat\ph_n(\ze)} \le \hat\Phi_n\bigl(\abs\ze\bigr)\,\ee^{\tau\abs\ze},
\qquad \ze\in\cR_c\ss0,
$$
with 
$$
\hat\Phi_1(\xi) = 2\la C, \quad
\hat\Phi_n = 2\la C \sum_{ \genfrac{}{}{0pt}{}{r\ge1}{n_1+\cdots+n_r=n-1} }
\frac{\xi^r}{r!} * \hat\Phi_{n_1} * \cdots * \hat\Phi_{n_r}
$$
(we also used the fact that
$\abs{\hat\al(\ze)}\le A\bigl(\abs\ze\bigr)\,\ee^{a\abs\ze}$ and
$\abs{\hat\be(\ze)}\le B\bigl(\abs\ze\bigr)\,\ee^{b\abs\ze}$ imply
$\abs{\hat\al*\hat\be(\ze)}\le (A*B)\bigl(\abs\ze\bigr) \, 
\ee^{\max(a,b)\abs\ze}$, and that
$\frac{1}{\xi}\bigl( (\xi A)*B \bigr) \le A*B$ for $\xi\ge0$).
The generating series $\hat\Phi=\sum\eps^n\hat\Phi_n$ is the formal Borel
transform of the solution $\ti\Phi=\sum\eps^n\ti\Phi_n$ of the equation
$\ti\Phi = 2\eps\la C z\ii + 2\eps\la C \sum z^{-r-1}\ti\Phi^r$.
We get $\ti\Phi = \frac{ 1 - (1-8\eps\la C z^{-2})^{1/2} }{2 z\ii}$
by solving this algebraic equation of degree~$2$,
hence $\ti\Phi_n(z) = \ga_n z^{-2n+1}$ with $0<\ga_n\le \Ga^n$ (with an explicit
$\Ga>0$ depending on~$\la C$),
and finally $\hat\Phi_n(\abs\ze) \le 
\Ga^n\frac{ \abs\ze^{2(n-1)} }{ (2(n-1))! }$.
Therefore the series $\sum\hat\ph_n$ converges in~$\cR_c\ss0$ and~$\hat\ph$
extends to the principal sheet of~$\cR$ with at most exponential growth.

\sm

The analytic continuation to the rest of~$\cR$ is more difficult. A natural idea
would be to try to extend the previous method of majorants to every half-sheet
of~$\cR$, but the problem is to find a suitable generalisation of
inequality~\refeq{ineqmajconv}. 
As shown in~\cite{sauzin01} or~\cite{OSS}, this can be done in the
union~$\cR\ss1$ of the half-sheets which are contiguous to the principal sheet, \ie
the ones which are reached after crossing the imaginary axis exactly once;
indeed, the symmetrically contractile paths~$\Ga$ constructed in
Lemma~\ref{lemRsymcontr} can be described quite simply for the points~$\ze$
belonging to these half-sheets and it is possible to define an
analogue~$\cR_c\ss1$ of~$\cR_c\ss0$.
But it is not so for the general half-sheets of~$\cR$, because of the complexity
of the symmetrically contractile paths which are needed.
The remedy employed in~\cite{sauzin01} and~\cite{OSS} consists in performing the
resurgent analysis, \ie describing the action of the alien derivations~$\De_\om$
to be defined in Section~\ref{secsimpl}, gradually: 
the possibility of following the analytic continuation of~$\hat\ph$ in~$\cR\ss1$
is sufficient to define $\De_{2\pi\I}\ti\ph$ and $\De_{-2\pi\I}\ti\ph$, which
amounts to computing the difference between the principal branch of~$\hat\ph$
at a given point~$\ze$ and the branch~$\hat\ph^\pm(\ze)$ of~$\hat\ph$
obtained by turning once around~$\pm2\pi\I$ and coming back at~$\ze$;
one then discovers that this difference $\hat\ph^\pm(\ze)-\hat\ph(\ze)$ is
proportional to~$\hat\ph(\ze\mp2\pi\I)$ (we shall try to make clear the reason
of this phenomenon in Section~\ref{secParBridge});
$\hat\ph^\pm$ is thus a function continuable along paths which cross the
imaginary axis once (as the sum of such functions), which means that~$\hat\ph$
is continuable to a set~$\cR\ss2$ defined by paths which are authorized to cross
two times (provided the first time is between~$2\pi\I$ and~$4\pi\I$ or
between~$-2\pi\I$ and~$-4\pi\I$).
The Riemann surface~$\cR$ can then be explored progressively, using more and
more alien derivations, $\De_{\pm4\pi\I}$ and
$\De_{\pm2\pi\I}\circ\De_{\pm2\pi\I}$ to reach a set~$\cR\ss3$, etc.
\eopf

\subsection{Sectorial normalisations (Fatou coordinates) and nonlinear Stokes
phenomenon (horn maps)}

We now apply the Borel-Laplace summation process and immediately get

\begin{cor}	\label{coroll}
With the hypothesis and notations of Theorem~\ref{thmresur}, for every
$\th_0\in\left]0,\frac{\pi}{2}\right[$, there exists $\tau>0$ such that the
Borel-Laplace sums
\begin{align*}
\ph^+ = \cL^\th\hat\ph, \quad \psi^+ = \cL^\th\hat\psi,
\qquad & -\th_0 \le \th \le \th_0, \\
\ph^- = \cL^\th\hat\ph, \quad \psi^- = \cL^\th\hat\psi, 
\qquad & \pi-\th_0 \le \th \le \pi+\th_0
\end{align*}
are analytic in~$\cD^+$, \resp $\cD^-$, where
$$
\cD^+ = \bigcup_{-\th_0 \le \th \le \th_0} \ao \RE(z\,\ee^{\I\th}) > \tau \af,
\quad
\cD^- = \bigcup_{\pi-\th_0 \le \th \le \pi+\th_0} \ao \RE(z\,\ee^{\I\th}) > \tau \af,
$$
and define transformations $u^\pm=\ID+\ph^\pm$ and $v^\pm=\ID+\psi^\pm$ which
satisfy 
\begin{gather*}
v^+\circ f = v^+ + 1 \quad\text{and}\quad
u^+\circ v^+ = v^+\circ u^+ = \ID \quad\text{on $\cD^+$},\\
v^-\circ f = v^- + 1 \quad\text{and}\quad
u^-\circ v^- = v^-\circ u^- = \ID \quad\text{on $\cD^-$}.
\end{gather*}
\end{cor}

\sm

\begin{figure}

\begin{center}

\epsfig{file=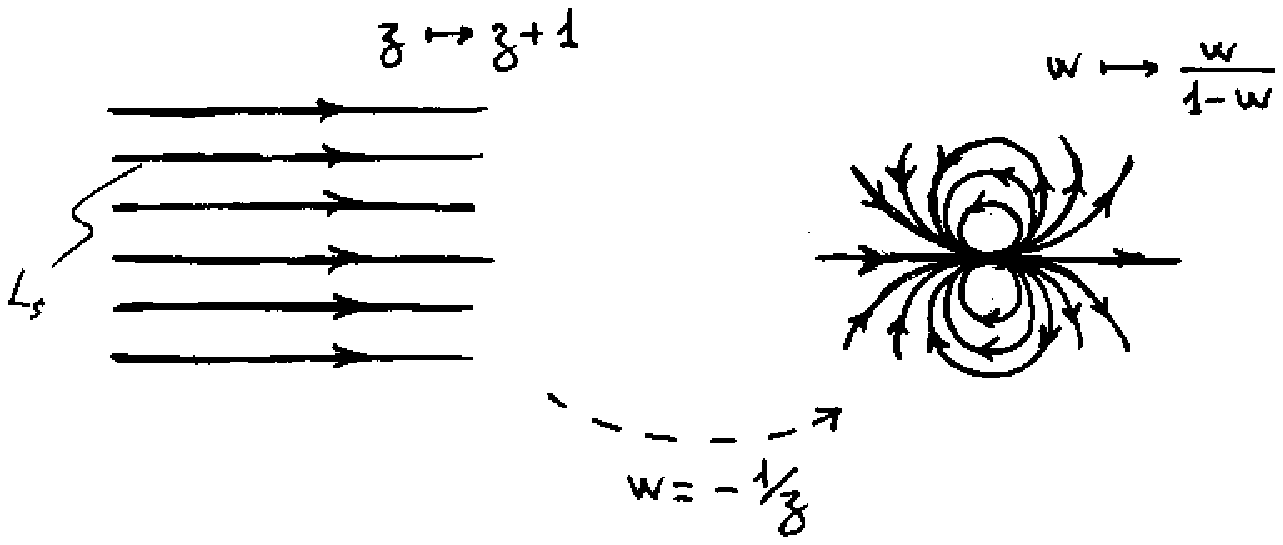,height=5.2cm,angle = 0}

\end{center}

\vspace{-.5cm}

\caption{\label{figcircles} Dynamics of $w\mapsto F_0(w)=\frac{w}{1-w}$.}


\vspace{.5cm}

\begin{center}

\epsfig{file=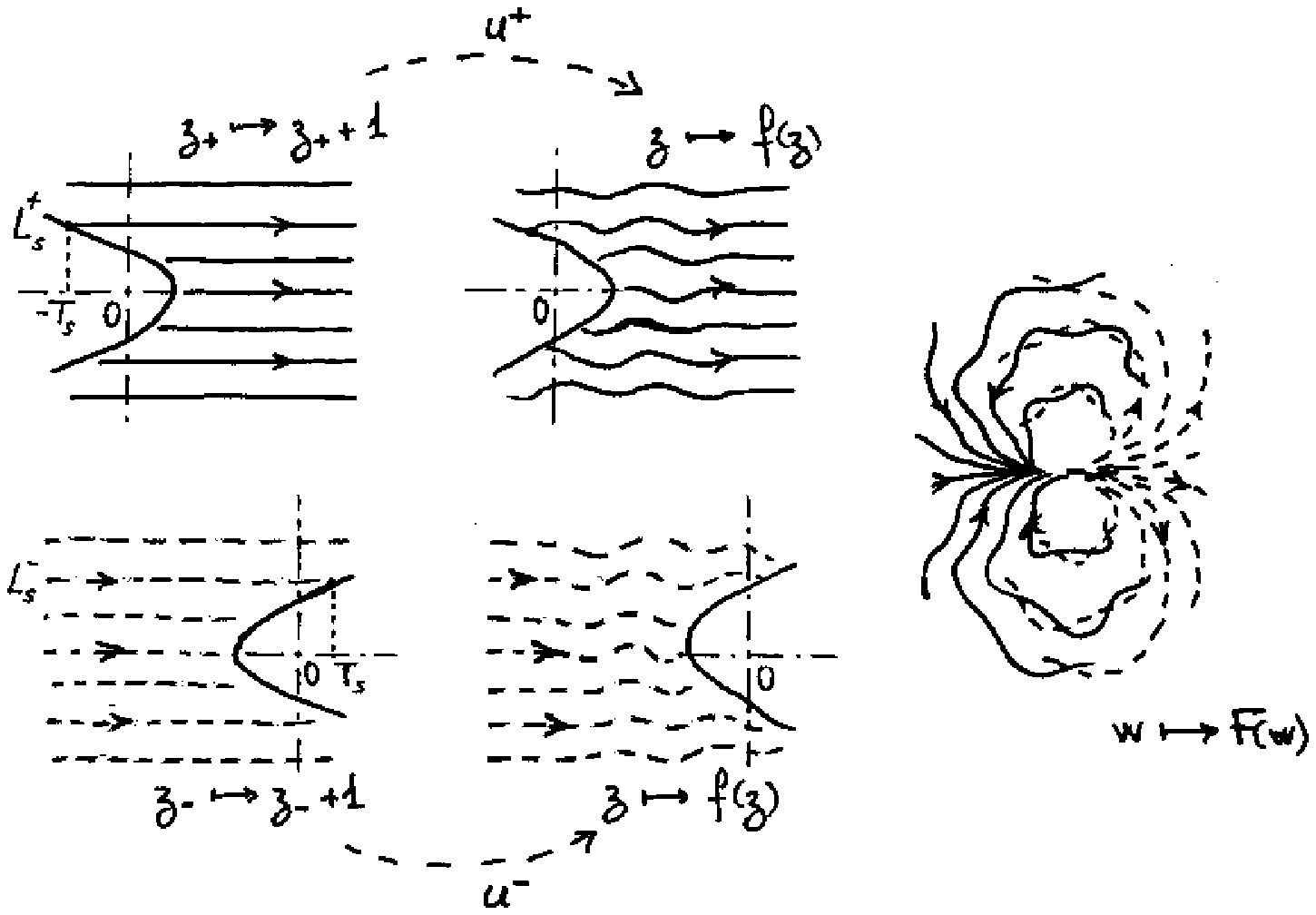,height=9.7cm,angle = 0}

\end{center}

\vspace{-.5cm}

\caption{\label{figsplitfoli} Stable and unstable foliations for~$f$ and~$F$.}

\end{figure}

One can consider $z^+ = v^+(z) = z+\psi^+(z)$ and $z^- = v^-(z) = z+\psi^-(z)$
as normalising coordinates for~$F$; they are sometimes called ``Fatou
coordinates''.\footnote{
Observe that, when the parabolic germ at the origin $F(w)\in w\C\{w\}$ extends to
an entire function, the function
$U^-(z) = -1/u^-(z)$ also extends to an entire function (because the domain of
analyticity~$\cD^-$ contains the half-plane $\RE z< -\tau$ and the relation
$U^-(z+1) = F\bigl(U^-(z)\bigr)$ allows one to define the analytic continuation
of~$U^-$ by $U^-(z) = F^n\bigl(U^-(z-n)\bigr)$, with $n\ge1$ large enough for a
given~$z$), 
which admits $-1/\ti u(z) = -z\ii\bigl(1+z\ii\ti\ph(z)\bigr)\ii \in
z\ii\C[[z\ii]]$ as asymptotic expansion in~$\cD^-$.
In this case, the formal series~$\ti\ph(z)$ must be divergent (if not, $-1/\ti
u(z)$ would be convergent, $U^-$ would be its sum and this entire function
would have to be constant), as well as~$\ti\psi(z)$, and the Fatou
coordinates~$v^+$ and~$v^-$ cannot be the analytic continuation one of the
other. 
We have a similar situation when $F\ii(w)\in w\C\{w\}$ extends to an entire
function, with $U^+(z) = -1/u^+(z)$ entire.
}
When expressed in these coordinates, the germ~$F$ simply reads $z^\pm\mapsto
z^\pm+1$ (see Figure~\ref{figsplitfoli}), the complexity of the dynamics being
hidden in the fact that neither~$v^+$ nor~$v^-$ is defined in a whole
neighbourhood of infinity
and that these transformations do not coincide on the two connected
components of~$\cD^+\cap\cD^-$
(except of course if~$F$ and~$F_0$ are analytically conjugated)---see the note
at the end of this section for a more ``dynamical'' and quicker construction
of~$v^\pm$.

\sm

This complexity can be analysed through the change of chart $v^+\circ u^- = \ID
+ \chi$, which is a priori defined on $\cE = \cD^-\cap (u^-)\ii(\cD^+)$; this
set has an ``upper'' and a ``lower'' connected components, $\cE^{\text{\rm up}}$
and~$\cE^{\text{\rm low}}$ (because $(u^-)\ii(\cD^+)$ is a sectorial neighbourhood
of infinity of the same kind as~$\cD^+$), and we thus get two analytic
functions~$\chi^{\text{\rm up}}$ and~$\chi^{\text{\rm low}}$
(this situation is reminiscent of the one described in Section~\ref{secdiffeq}).
The conjugacy equations satisfied by~$u^-$ and~$v^+$ yield $\chi(z+1)=\chi(z)$,
hence both~$\chi^{\text{\rm up}}$ and~$\chi^{\text{\rm low}}$ are $1$-periodic;
moreover, we know that these functions tend to~$0$ as $\IM z \to\pm\infty$
(faster than any power of~$z\ii$, by composition of asymptotic expansions).
We thus get two Fourier series 
\begin{alignat}{2}
\label{eqFourierLOW}
\chi^{\text{\rm low}}(z) &= v^+\circ u^-(z) - z\, &= \sum_{m\ge1} B_m \, \ee^{-2\pi\I m z}, 
\qquad & \IM z < -\tau_0, \\
\label{eqFourierUP}
\chi^{\text{\rm up}}(z) &= v^+\circ u^-(z) - z\, &= \sum_{m\le-1} B_m \, \ee^{-2\pi\I m z}, 
\quad\ens & \IM z > \tau_0,
\end{alignat}
which are convergent for $\tau_0>0$ large enough.
It turns out that the classification problem can be solved this way: {\em two
non-degenerate parabolic germs with vanishing resiter are analytically conjugate
if and only if they define the same pair of Fourier series $(\chi^{\text{\rm
up}},\chi^{\text{\rm low}})$ up to a change of variable $z\mapsto z+c$;
morevover, any pair of Fourier series of the
type~\refeq{eqFourierLOW}--\refeq{eqFourierUP} can be obtained this
way.}\footnote{		
\label{footanalyticinv}
For the first statement, consider~$f_1$ and~$f_2$ satisfying
$\chi_2^{\text{\rm up}}(z) = \chi_1^{\text{\rm up}}(z+c)$
and $\chi_2^{\text{\rm low}}(z) = \chi_1^{\text{\rm low}}(z+c)$
with $c\in\C$,
thus $v_2^+\circ u_2^- = \tau\ii\circ v_1^+\circ u_1^-\circ\tau$
in~$\cE^{\text{\rm up}}$ and~$\cE^{\text{\rm low}}$,
with $\tau(z) = z+c$.
Using $(\ti u_1\circ\tau,\tau\ii\circ\ti v_1)$
instead of $(\ti u_1,\ti v_1)$, we see that a formal conjugacy between~$f_1$
and~$f_2$ is given by $\ti u_2\circ\tau\ii\circ\ti v_1$; its Borel-Laplace sums $u_2^+\circ
\tau\ii\circ v_1^+$ and $u_2^-\circ\tau\ii\circ v_1^-$ can be glued together and give
rise to an analytic conjugacy, since
$u_2^- = u_2^+\circ\tau\ii\circ v_1^+\circ u_1^-\circ\tau$.
Conversely, if there exists $h\in\ID+\C\{z\ii\}$ such that $f_2\circ h = h\circ f_1$,
we see that $h\circ\ti u_1$ establishes a formal conjugacy between~$f_2$ and
$z\mapsto z+1$, Proposition~\ref{propformalconj} thus implies the existence of
$c\in\C$ such that $\ti u_2 = h\circ\ti u_1\circ\tau$
and $\ti v_2=\tau\ii\circ\ti v_1\circ h\ii$, with $\tau(z) = z+c$, whence
$u_2^\pm = h\circ u_1^\pm\circ\tau$ and $v_2^\pm=\tau\ii\circ v_1^\pm\circ h\ii$,
and $v_2^+\circ u_2^- = \tau\ii\circ v_1^+\circ u_1^- \circ\tau$, as desired.
The proof of the second statement is beyond the scope of the present text.
}
The numbers~$B_m$ are said to be ``analytic invariants'' for the germ~$f$
or~$F$.
The functions~$\ID+\chi^{\text{\rm low}}$ and~$\ID+\chi^{\text{\rm up}}$ themselves are
called ``horn maps''.\footnote{
In fact, this name (which is of A.~Douady's coinage) usually refers to the maps
$\ID+\chi^{\text{\rm low}}$ expressed in the coordinate $w_-=e^{-2\pi\I z}$, \ie
$w_-\mapsto w_-\,\exp\left(-2\pi\I \sum_{m\ge1}B_m w_-^m\right)$, and
$\ID+\chi^{\text{\rm up}}$ expressed in the coordinate $w_+=e^{2\pi\I z}$, \ie
$w_+\mapsto w_+\,\exp\left(2\pi\I \sum_{m\ge1}B_{-m} w_+^m\right)$, which are
holormophic for $|w_\pm|<\ee^{-2\pi\tau_0}$ and can be viewed as return
maps;
the variables~$w_\pm$ are natural coordinates at both ends of ``\'Ecalle's
cylinder''.
See \cite{MarRam}, \cite{DouadyHubb}, \cite{Milnor}, \cite{Shi}, \cite{Shi00},
\cite{Loray}.
}

\sm

Now, consider the relation~\refeq{eqFourierLOW} for instance. In the half-plane
$\Pi^{\text{\rm low}} = \ao\IM z < -\tau_0\af$, $|\chi^{\text{\rm low}}(z)|$ is
uniformly bounded by a constant which can be made arbitrarily small by
choosing~$\tau_0$ large enough.
We can use this to extend analytically~$\ph^-$ beyond~$\cD^-$,
in~$\Pi^{\text{\rm low}}$: 
we can indeed write 
$u^-(z)=u^+\bigl(z+\chi^{\text{\rm low}}(z)\bigr)$
if $z\in\cE^{\text{\rm low}}$, $\RE z>0$ and $\IM z < -\tau_0$, and the \rhs\ in this identity is
holomorphic on $\Pi^{\text{\rm low}}\cap\ao \RE z>0 \af$ (because the image of
this domain by $\ID+\chi^{\text{\rm low}}$ is included in~$\cD^+$).
However, this implies $\ph^-(z) = \chi^{\text{\rm low}}(z) +
\ph^+\bigl(z+\chi^{\text{\rm low}}(z)\bigr) = \chi^{\text{\rm low}}(z) + \cO(1/z)$ for $z$
tending to infinity in $\Pi^{\text{\rm low}}\cap\ao \RE z>0 \af$,
thus~$\ph^-(z)$ is no longer asymptotic to~$\ti\ph$ there, an oscillating term
shows up when~$z$ moves along any horizontal half-line $-\I\,s+\R^+$, $s>\tau_0$.
Similarly, $\ph^-(z) = \chi^{\text{\rm up}}(z) +
\ph^+\bigl(z+\chi^{\text{\rm up}}(z)\bigr)$ extends analytically to the
half-plane $\ao\IM>\tau_0\af$ if~$\tau_0$ is large enough, with an oscillating
asymptotic behaviour determined by~$\chi^{\text{\rm up}}$.
This can be considered as a nonlinear analogue of the classical Stokes
phenomenon (well-known in the case of linear ODEs).

\sm

In the next sections, we shall outline the resurgent approach, which consists in
extracting information from the singularities of~$\hat\ph$ or~$\hat\psi$ in
order to construct a set of analytic invariants $\{ A_{2\pi\I m}, \; m\in\Z^*
\}$, and its relation with $\{ B_m, \; m\in\Z^* \}$.
Before explaining this, let us mention an interpretation of the
non-coincidence of~$v^+$ and~$v^-$ as a ``splitting problem'', following~\cite{Gel}.

\subsub{Splitting of the invariant foliation}

The dynamical behaviour of~$F_0$ is easily visualized: the invariant foliation
by horizontal lines for~$f_0$ gives rise to an invariant foliation by circles
for~$F_0$, as shown on Figure~\ref{figcircles}
(notice that, for a global understanding of the dynamics of~$F_0$, one should
let~$w$ vary on the Riemann sphere, including the point at infinity, which is
the image of~$1$).
For the given parabolic germ~$f$ or~$F$, we can use the Fatou coordinates to
define ``stable'' and ``unstable'' foliations: for each $s\in\R$, the line $L_s
= \{ t+\I s,\; t\in\R \}$ intersects~$\cD^+$ and~$\cD^-$ along half-lines $L_s^+
= \{ t+\I s,\; t>-T_s \}$ and $L_s^- = \{ t+\I s,\; t<T_s \}$ and we may set
$Z^+_s(t)=u^+(t+\I s)$, $Z^-_s(t)=u^-(t+\I s)$
and $W_s^\pm(t) = -1/Z_s^\pm(t)$ for $t>-T_s$, \resp $t<T_s$.
We have
$$
f^n\bigl(Z^\pm_s(t)\bigr) = Z^\pm_s(t+n) = t + n + \I s + O\bigl((t+n)\ii\bigr), \qquad
F^n\bigl(W^\pm_s(t)\bigr) = W^\pm_s(t+n)
   \xrightarrow[n\to\pm\infty]{} 0.
$$
The invariant foliation $\{ L_s \}_{s\in\R}$ of~$F_0$ is so to say split, giving
rise to two foliations $\{ u^+(L_s^+) \}_{s\in\R}$ and $\{ u^-(L_s^-)
\}_{s\in\R}$ which in general do not coincide but can be compared for
$|s|>\tau_0$ large enough (because $T_s$ is then positive---see
Figure~\ref{figsplitfoli}; one can also use the analytic continuation of~$u^\pm$
to~$\Pi^{\text{\rm low}}$ and~$\Pi^{\text{\rm up}}$ and consider $u^\pm(L_s)$
for $|s|$ large enough).
It is proven in~\cite{Gel} that, if $B_1 = \cdots = B_{n-1} = 0$ and $B_n\neq0$,
for $s<-\tau_0$ the curves $\{ u^+(L_s^+) \}_{s\in\R}$ and $\{ u^-(L_s^-)
\}_{s\in\R}$ intersect along $2n$ orbits of~$f$, and that, for each of these
orbits, the intersection angle (which is the same for all the points of the
orbit because $f$ is conformal) is $2\pi n|B_n|\ee^{-2\pi n|s|} +
\cO(\ee^{-2\pi(n+1)|s|})$ (there is a symmetric result for $s>\tau_0$ in terms of
$B_{-1}, B_{-2},\ldots$).

\subsub{Note on another construction of the Fatou coordinates}

The result expressed in Corollary~\ref{coroll} can be obtained through formulas
which are reminiscent of the Note at the end of Section~\ref{seclineardiffeq},
being analogous to the explicit formulas available for the linear case.
Indeed, observe that, since $f=\ID+1+a$ with $\lim_{|z|\to\infty}|a(z)|=0$, $f(\cD^+)
\subset \cD^+$ and $f\ii(\cD^-) \subset \cD^-$ (provided these sets are defined
using a large enough constant~$\tau$); one can thus iterate Abel's equation
forward or backward in these domains: setting
$v^\pm=\ID+\psi^\pm$, equation~\refeq{eqnAbelpsi} yields 
$$
\psi^+\circ f^{k+1} - \psi^+\circ f^{k}  = - a\circ f^{k} \ens\text{in $\cD^+$}, 
\qquad
\psi^-\circ f^{-(k+1)} - \psi^-\circ f^{-k}  = a\circ f^{-(k+1)} \ens\text{in $\cD^-$}, 
\qquad k\ge0,
$$
where $f^n$ denotes the $n^\text{th}$ power of iteration of~$f$ for any
$n\in\Z$.
The assumption $a(z)\in z^{-2}\C\{z\ii\}$ implies the existence and unicity of
solutions which tend to~$0$ at infinity, which can be expressed as uniformly
convergent series
\begin{equation}	\label{eqseriepsi}
\psi^+ = \sum_{k\ge0}{a\circ f^k}, \qquad
\psi^- = - \sum_{k\ge1}{a\circ f^{-k}}.
\end{equation}
(One can check that $f^n = \ID + n + \sum_{k=0}^{n-1} a\circ f^k$ by induction
on~$n\ge1$, and similarly, that $f\ii = \ID - 1 - a\circ f\ii$, hence
$f^{-n} = \ID - n - \sum_{k=1}^{n} a\circ f^{-k}$;
formulas~\refeq{eqseriepsi} are thus consistent with the formulas
$$
v^\pm(z) = \lim_{n\to\pm\infty} \bigl( f^n(z)-n \bigr), \qquad
u^\pm(z) = \lim_{n\to\pm\infty} f^{-n}(z+n)
$$
which can be found in~\cite[Vol.~2, p.~322]{Eca81} or~\cite[p.~39]{Loray}.)

As a consequence, if we introduce the $f$-invariant functions~$\be^{\text{\rm
up}}$ and~$\be^{\text{\rm low}}$ defined by
$$
\be = \sum_{k\in\Z} a\circ f^k
$$
in the two connected components of~$\cD^+\cap\cD^-$ (they are first integrals of
the dynamics: $\be^{\text{\rm low}}=\be^{\text{\rm low}}\circ f$,
$\be^{\text{\rm up}}=\be^{\text{\rm up}}\circ f$),
we get $v^+-v^- = \psi^+-\psi^- = \be$.
Since the horn maps~$\ID+\chi^{\text{\rm low}}$ and~$\ID+\chi^{\text{\rm up}}$
are defined by $v^+ = (\ID+\chi)\circ v^-$, 
writing $v^- = \ID + \psi^-$ and $v^+ = \ID + \psi^- + \be$ we finally get
$$
\chi^{\text{\rm up}} = \be^{\text{\rm up}} \circ (\ID+\psi^-)\ii, \qquad 
\chi^{\text{\rm low}} = \be^{\text{\rm low}} \circ (\ID+\psi^-)\ii.
$$


\subsection{Alien calculus for simple resurgent fuctions}	\label{secsimpl}


In Section~\ref{secAbEqtgtId}, we saw how solving a difference equation could lead to
Borel transforms~$\hat\ph(\ze)$ or~$\hat\psi(\ze)$ which are holomorphic
in~$\cR$ and likely to develop singularities at the integer multiples
of~$2\pi\I$.
The purpose of ``alien calculus'' is to give an efficient way of encoding these
singularities and of obtaining information on them.
We shall describe the general formalism of singularities and give the definition
of alien derivations in Section~\ref{secSingul}, but we begin here with a class of
resurgent functions for which the definitions are simpler and which is
sufficient to deal with Abel's equation.

\subsub{Simple resurgent functions}

\begin{Def} 	\label{defsimplesing}
Let $\om\in\C$.
We say that a function~$\hat\ph$, which is holomorphic in an open disc
$D\subset\C$ to which~$\om$ is adherent, ``has a simple singularity at~$\om$''
if there exist $C\in\C$ and $\hat\Phi(\ze),\mathrm{reg}(\ze)\in\C\{\ze\}$ such
that 
\begin{equation}	\label{eqsimplesing}
\hat\ph(\ze) =
\frac{C}{2\pi\I(\ze-\om)} + \frac{1}{2\pi\I}\hat\Phi(\ze-\om)\log(\ze-\om) 
+ \mathrm{reg}(\ze-\om)
\end{equation}
for all $\ze\in D$ with $|\ze-\om|$ small enough.
We then use the notation
$$
\Sing_\om \hat\ph = C\,\de + \hat\Phi \in \C\,\de \oplus \C\{\ze\}.
$$
\end{Def}

Obviously, a change of branch of the logarithm in~\refeq{eqsimplesing}
only results in the replacement of~$\mathrm{reg}(\ze)$ by another regular germ;
the interesting part of the formula is the singularity encoded by the
``residuum''~$C$ and the ``variation''~$\hat\Phi$ which are
unambiguously determined.
For instance, the variation (also called ``the minor of the singularity
of~$\hat\ph$ at~$\om$'') can be written
$$
\hat\Phi(\ze) = \hat\ph(\om+\ze) - \hat\ph(\om+\ze\,\ee^{-2\pi\I}),
$$
where it is understood that considering $\om+\ze\,\ee^{-2\pi\I}$ means following
the analytic continuation of~$\hat\ph$ along the circular path 
$t\in[0,1]\mapsto \om+\ze\,\ee^{-2\pi\I t}$ 
(the analytic continuation exists when~$|\ze|$ is small enough, since $\hat\Phi$
and~$\mathrm{reg}$ are regular near the origin).

Any analytic function~$\chn\Phi(\ze)$ which differ from~$\hat\ph(\om+\ze)$ by a
regular germ is called a ``major'' of the singularity $\Sing_\om\hat\ph$. Any
major thus satisfies
$$
\Sing_0\chn\Phi = \Sing_\om\hat\ph, \qquad 
\hat\Phi(\ze) = \chn\Phi(\ze) - \chn\Phi(\ze\,\ee^{-2\pi\I})
$$
(the minor is the variation of any major).
In fact, the singularity can be identified with an equivalence class of majors
modulo~$\C\{\ze\}$. This will be used as a way of generalising the previous
definition to deal with more complicated singularities in
Section~\ref{secSingul}.

\sm

For any path~$\ga$ issuing from~$0$ and lying in $\C\setminus 2\pi\I\,\Z$, and
for any $\hat\ph\in\wh\cH(\cR)$, we denote by~$\cont_\ga\hat\ph$
the branch of~$\hat\ph$ obtained by following the analytic continuation
of~$\hat\ph$ along~$\ga$, which is a function holomorphic at least in any open
disc containing the extremity~$\ga(1)$ of~$\ga$ and avoiding $2\pi\I\,\Z$.
In particular, if $\om\in2\pi\I\,\Z$ satisfies $|\om-\ga(1)|<\pi$, there is a
disc~$D$ avoiding $2\pi\I\,\Z$ which contains~$\ga(1)$ and to which~$\om$ is adherent.

\begin{Def}	\label{defsimpleresfcn}
A ``simple resurgent function'' is any  
$\hat\chi=c\,\de + \hat\ph(\ze) \in \C\,\de\oplus\wh\cH(\cR)$ 
such that, for each $\om\in2\pi\I\,\Z$ and 
for each path~$\ga$ which starts from~$0$, lies in $\C\setminus 2\pi\I\,\Z$ and
has its extremity in the disc of radius~$\pi$ centred at~$\om$,
the branch $\cont_\ga\hat\ph$ has a simple singularity
at~$\om$.
\end{Def}
One can check that the corresponding subspace of~$\C\,\de\oplus\wh\cH(\cR)$ is
stable by convolution:

\begin{prop}	\label{propsimplestab}
The subspace~$\Rsimp$ consisting of all simple resurgent functions is a
subalgebra of the convolution algebra~$\C\,\de\oplus\wh\cH(\cR)$.

\end{prop}
(This can be done with the help of the symmetrically contractile paths of
Lemma~\ref{lemRsymcontr}; see the arguments given in the proof of
Lemma~\ref{lemDeomplus} below.)
%

\sm

As a consequence, {\em the Borel transform~$\hat\ph$ of the solution of
equation~\refeq{eqniter}, which belongs to~$\wh\cH(\cR)$ according to
Theorem~\ref{thmresur}, must be a simple resurgent function.}
Indeed, as indicated in the proof of this theorem, $\hat\ph$ can be expressed as
a uniformly convergent series $\sum_{n\ge1} \hat\ph_n$, where the functions
$\hat\ph_n\in\wh\cH(\cR)$ are defined inductively by~\refeq{eqseriesfcns}.
It is essentially sufficient to check that each~$\hat\ph_n$ belongs
to~$\Rsimp$, and this is easily done by induction 
%
%
($\hat\ph_1$ is meromorphic with simple poles because~$E$ is;
since $\hat a_r\in\ze^2\C\{\ze\}$ is entire, one can write $\hat a_r=1*1*\hat
a_r''$ with $\hat a_r''\in\Rsimp$, hence $\hat\ph_n = E\,(1*1*\hat A_n)$ with
$\hat A_n\in\Rsimp$ by the inductive hypothesis; 
one concludes by observing that the singularities of~$1*1*\hat A_n$ have no
residuum and that their variations have valuation at least~$1$ at the origin,
while $E(\ze) = -\frac{1}{\ze-\om}+\mathrm{reg}(\ze-\om)$).

\sm

{\em The Borel transform~$\hat\psi$ of the solution of equation~\refeq{eqnAbelpsi} is
also a simple resurgent function}, because the space~$\Rsimp$ enjoys stability
properties similar to the properties of $\C\,\de\oplus\wh\cH(\cR)$ indicated in
Propositions~\ref{propcompos} and~\ref{propInversion} of
Section~\ref{secformconvmod}.

\sm

What we just defined is the ``convolutive model of the algebra of simple
resurgent functions''.
The ``formal model of the algebra of simple resurgent functions'' $\wtRsimp$ is
defined as the subalgebra of~$\C\oplus\ti\cH$ obtained by pulling back~$\Rsimp$
by~$\cB$.
The formal solution $\ti v=\ID+\ti\psi$ of Abel's equation can thus be viewed as
a simple resurgent change of variable which normalises~$f$, with $\ti
u=\ID+\ti\ph\in\ID+\wtRsimp$ as inverse transformation.

\subsub{Alien derivations}

Let~$\om$ and~$\ga$ be as in Definition~\ref{defsimpleresfcn}. 
For $c\,\de + \hat\ph\in\Rsimp$, let~$C_\ga$ and~$\hat\Phi_\ga$ denote the
residuum and the minor (variation) of the singularity of $\cont_\ga\hat\ph$ at~$\om$.
It is clear that~$\hat\Phi_\ga$, which is holomorphic at the origin by assumption,
must extend holomorphically to~$\cR$ and possess itself only simple
singularities. 
The formula
$$
c\,\de + \hat\ph \ens \mapsto \ens C_\ga\,\de + \hat\Phi_\ga
= \Sing_\om \left( \cont_\ga \hat\ph \right) 
$$
thus defines a linear operator of~$\Rsimp$ to itself.

\sm

\begin{figure}

\begin{center}

\epsfig{file=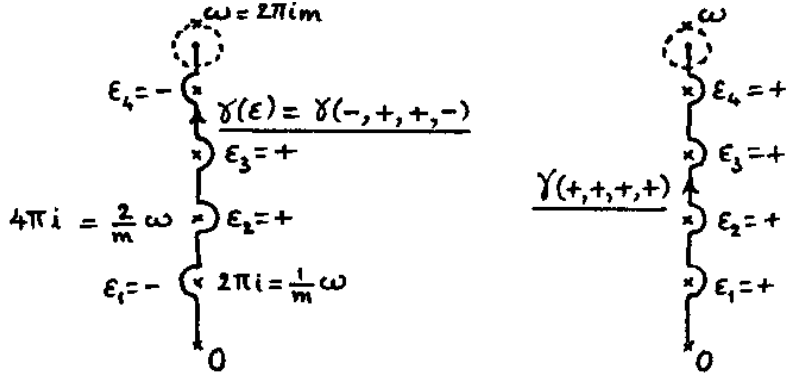,height=5.1cm,angle = 0}

\end{center}

\vspace{-.5cm}

\caption{\label{figalienpath} Paths for the definition of~$\De_\om\hat\chi$ or~$\De^+_\om\hat\chi$.}

\end{figure}

\label{secAlienDer}
\begin{Def}	\label{defDeom}
For each $\om=\pm2\pi\I\, m$, $m\in\N^*$, we define a linear operator~$\De_\om$
from~$\Rsimp$ to itself by using $2^{m-1}$ particular paths~$\ga$:
\begin{equation}    \label{eqdefDeomgen}
\De_\om(c\,\de + \hat\ph) = \sum_{\eps_1,\ldots,\eps_{m-1}\in\{+,-\}}
                \frac{p(\eps)!q(\eps)!}{m!} \left(C_{\ga(\eps)} \de +
                  \hat\Phi_{\ga(\eps)} \right),
\end{equation}
where
$p(\eps)$ and $q(\eps)=m-1-p(\eps)$ denote the numbers of signs~`$+$' and of signs~`$-$'
in the sequence~$\eps$,
and the path~${\ga(\eps)}$ connects~$\left]0,\frac{1}{m}\om\right[$ and
$\left]\frac{m-1}{m}\om,\om\right[$, following the segment~$\left]0,\om\right[$
but circumventing the intermediary singular points~$\frac{r}{m}\om$ to the right
if~$\eps_r=+$ and to the left if~$\eps_r=-$
(see Figure~\ref{figalienpath}).
We also define a linear operator~$\De^+_\om$ from~$\Rsimp$ to itself by setting
\begin{equation}    \label{eqdefDeomplus}
\De^+_\om(c\,\de + \hat\ph) = C_{\ga_\om} \de +
                  \hat\Phi_{\ga_\om},
\qquad \ga_\om = \ga(+,\ldots,+).
\end{equation}
Via~$\cB$, the operators~$\De_\om$ or~$\De^+_\om$ of the convolutive model give
rise to operators of the formal model, which we denote by the same symbols:
$$
\xymatrix{
{\rule[-.6ex]{0ex}{0ex} \makebox[5em]{$\Rsimp$}} \ar[r]^*+{\De_\om} 
& {\rule[-.6ex]{0ex}{0ex} \makebox[5em]{$\Rsimp$}} \\
{\rule[-1ex]{0ex}{3ex} \makebox[5em]{$\wtRsimp$}} \ar[r]^*+{\De_\om} \ar[u]^*+{\cB}
& {\rule[-1ex]{0ex}{3ex} \makebox[5em]{$\wtRsimp$}} \ar[u]_*+{\cB} 
} \qquad
\xymatrix{
{\rule[-.6ex]{0ex}{0ex} \makebox[5em]{$\Rsimp$}} \ar[r]^*+{\De^+_\om} 
& {\rule[-.6ex]{0ex}{0ex} \makebox[5em]{$\Rsimp$}} \\
{\rule[-1ex]{0ex}{3ex} \makebox[5em]{$\wtRsimp$}} \ar[r]^*+{\De^+_\om} \ar[u]^*+{\cB}
& {\rule[-1ex]{0ex}{3ex} \makebox[5em]{$\wtRsimp$}} \ar[u]_*+{\cB} 
}
$$
\end{Def}

\begin{prop}	\label{propsimpleAlienD}
The operator $\De_\om$ is a derivation, \ie the Leibniz rule holds in the
convolutive model:
\begin{equation}	\label{eqLeibniz}
\De_\om (\hat\chi_1*\hat\chi_2) = (\De_\om\hat\chi_1)*\hat\chi_2 +
\hat\chi_1*(\De_\om\hat\chi_2),
\qquad \hat\chi_1,\hat\chi_2\in\Rsimp.
\end{equation}
Equivalently, we have in the formal model
\begin{equation}	\label{eqLeibnizmult}
\De_\om (\ti\chi_1\ti\chi_2) = (\De_\om\ti\chi_1)\ti\chi_2 +
\ti\chi_1(\De_\om\ti\chi_2),
\qquad \ti\chi_1,\ti\chi_2\in\wtRsimp.
\end{equation}
Moreover, for any $\ti\chi,\ti\chi_1,\ti\chi_2\in\wtRsimp$, 
\begin{gather}
\label{eqcommutnatur}
\De_\om\pa\ti\chi =
\pa\De_\om\ti\chi - \om \De_\om\ti\chi,\\[1ex]
\label{eqaliencompos}
\De_\om\bigl(\ti\chi_1\circ(\ID + \ti\chi_2)\bigr) = 
\ee^{-\om\ti\chi_2} \bigl( (\De_\om\ti\chi_1)\circ(\ID+\ti\chi_2) \bigr)
+ \bigl( (\pa\ti\chi_1)\circ(\ID+\ti\chi_2) \bigr) \De_\om\ti\chi_2.
\end{gather}
In particular $\De_\om$ commutes with $\ti\chi(z)\mapsto\ti\chi(z+1)$
for each $\om\in2\pi\I\,\Z^*$.
\end{prop}

\sm

Because of~\refeq{eqLeibniz} or~\refeq{eqLeibnizmult},
the operators~$\De_\om$ were called ``alien derivations'' by \'Ecalle, by
contrast with the natural derivation $\pa=\frac{\dd}{\dd z}$.
Some formulas get simpler when introducing the ``dotted alien derivations''
$\dDeom: \ti\chi(z) \mapsto \ee^{-\om z}\,\De_\om\ti\chi(z)$ (where $\ee^{-\om
z}$ is understood as a symbol external to~$\wtRsimp$, obeying the usual rules
\wrt\ multiplication and differentiation):
the dotted alien derivations commute with~$\pa$ and satisfy
\begin{equation}	\label{eqdotaliencompos}
\dDeom\bigl(\ti\chi_1\circ(\ID + \ti\chi_2)\bigr) = 
(\dDeom\ti\chi_1)\circ(\ID+\ti\chi_2) 
+ \bigl( (\pa\ti\chi_1)\circ(\ID+\ti\chi_2) \bigr)
\dDeom\ti\chi_2.
\end{equation}

\sm

There is no relation between the operators~$\De_\om$: they generate a free Lie
algebra.\footnote{
We mean that, for any~$N\ge1$ and for any collection of non-zero simple
resurgent functions $(\ti\chi_{\om_1 \cdots \om_r})$ indexed by finitely many
words $\om_1 \cdots
\om_r$ of any length $r\in\{0,\ldots,N\}$, the operator
$$
\sum_{r=0}^N \sum_{\om_1 \cdots \om_r} \,
\ti\chi_{\om_1 \cdots \om_r} \, \De_{\om_r}\cdots \De_{\om_1}
$$
(with the convention $\ti\chi_{\om_1 \cdots \om_r} \, \De_{\om_r}\cdots
\De_{\om_1} = \ti\chi_{\emptyset}$ if $r=0$) is not identically zero on~$\wtRsimp$.  
We omit the proof.
}
They provide a way of encoding the whole singular behaviour of a minor
$\hat\ph\in\Rsimp$: 
given a sequence $\om_1,\ldots,\om_r\in2\pi\I\,\Z^*$, the evaluation of the
composition $\De_{\om_r}\cdots\De_{\om_1}(c\,\de+\hat\ph)$ is a combination of
singularities at~$\om_1+\cdots+\om_r$ for various branchs of~$\hat\ph$.
Conversely, any singularity of any branch of~$\hat\ph$ can be computed if
the collection of these objects for all sequences $(\om_1,\ldots,\om_r)$ is
known.\footnote{
\label{footnotlimit}
One must not limit oneself to $r=1$. For instance, $\De_\om\hat\chi=0$ does not
mean that there is no singularity at~$\om$ for any branch of the minor;
consider for example $\om=\om_1+\om_2$ and $\hat\chi=\hat\ph*\hat\psi$ with
$\hat\ph=1/(\ze-\om_1)$, $\hat\psi=1/(\ze-\om_2)$ and $\om_1\neq\om_2$:
$\De_{\om_1}\hat\ph=\De_{\om_2}\hat\psi=2\pi\I\,\de$ and $\De_\om\hat\chi=0$,
but $\De_{\om_1}\hat\chi=2\pi\I\hat\psi$
and $\De_{\om_2}\hat\chi=2\pi\I\hat\ph$
imply $\De_{\om_2}\De_{\om_1}\hat\chi = 
\De_{\om_1}\De_{\om_2}\hat\chi = -4\pi^2\,\de$,
which reveals the presence of a singularity at~$\om$ at least for some
branchs of~$\hat\chi$.
}

\sm

Before proving Proposition~\ref{propsimpleAlienD}, we turn to the
operators~$\De^+_\om$. Their definition is simpler (\cf formula~\refeq{eqdefDeomplus}),
but they are not derivations (except for $\om=\pm2\pi\I$, since we have then
$\De^+_\om=\De_\om$).
Here is the way they act on products:

\begin{lemma}	\label{lemDeomplus}
For $\om\in2\pi\I\,\Z^*$ and $\ti\chi_1,\ti\chi_2\in\Rsimp$,
\begin{equation}	\label{eqLeibnizplus}
\De^+_\om (\hat\chi_1*\hat\chi_2) = (\De^+_\om\hat\chi_1)*\hat\chi_2 \;+
\sum (\De^+_{\om_1}\hat\chi_1)*(\De^+_{\om_2}\hat\chi_2)
\ens +\; \hat\chi_1*(\De^+_\om\hat\chi_2),
\end{equation}
where the sum extends to all $(\om_1,\om_2)$ such that
${\om_1+\om_2=\om}$ and $\om_j \in \left]0,\om\right[ \cap 2\pi\I\,\Z^*$.
\end{lemma}

\sm

\noindent {\em Proof of Lemma~\ref{lemDeomplus}.}
Formula~\refeq{eqLeibnizplus} results from a kind of combinatorics of
symmetrically contractile paths.
We begin with a proof in the simplest case, when $\om=2\pi\I$ (the case
of~$-2\pi\I$ is similar), and then sketch the case of $\om=2\pi\I\, m$ (the case
of~$-2\pi\I\, m$ would be similar).

\sm

Let $\om=2\pi\I$.
By linearity, observing that~$\De^+_\om$ annihilates the multiples of~$\de$, it is
sufficient to consider simple resurgent functions of the form
$\hat\chi_1=\hat\ph$ and $\hat\chi_2=\hat\psi$ (\ie without any multiple
of~$\de$). 
We must prove
$$
\De^+_{2\pi\I} (\hat\ph*\hat\psi) = (\De^+_{2\pi\I}\hat\ph)*\hat\psi
+ \hat\ph*(\De^+_{2\pi\I}\hat\psi).
$$
Let $\De^+_\om\hat\ph = a\,\de + \hat\Phi$ and $\De^+_\om\hat\psi = b\,\de +
\hat\Psi$, and consider the formula
\begin{equation}	\label{eqconvell}
\hat\ph*\hat\psi(\xi) = \int_0^\xi \hat\ph(\xi_1)\hat\psi(\xi-\xi_1)\,\dd\xi_1,
\end{equation}
which holds for~$\xi$ close to~$\om$ provided~$\xi$ lies in the principal sheet
of~$\cR$, \ie the segment $\ell=[0,\xi]$ avoids~$2\pi\I$ and~$-2\pi\I$.
Writing $\xi=\om+\ze$, we have
\begin{equation}	\label{eqsingphpsi}
\hat\ph(\om+\ze) = \frac{1}{2\pi\I} \left(
\frac{a}{\ze} + \hat\Phi(\ze)\log\ze + \mathrm{reg}(\ze) \right),
\quad
\hat\psi(\om+\ze) = \frac{1}{2\pi\I} \left(
\frac{b}{\ze} + \hat\Psi(\ze)\log\ze + \mathrm{reg}(\ze) \right).
\end{equation}
It can be seen that the residuum of~$\hat\ph*\hat\psi$ at~$\om$ is zero:
in the path of integration of~\refeq{eqconvell}, the singularity for $\xi\to\om$
stems from the extremities,~$\xi$ (because of~$\hat\ph(\xi_1)$, which is
multiplied by a function holomorphic for~$\xi_1$ close to~$\om$)
and~$0$ (because of~$\hat\psi(\xi-\xi_1)$, which is multiplied by a function
holomorphic for~$\xi_1$ close to~$0$); the singularity is thus obtained by
{\em integrating} simple poles and logarithmic singularities.\footnote{
\label{footsimplify}
This argument can be avoided by arguing as at the end of the case $m\ge2$,
writing 
$\hat\ph*\hat\psi = ( \frac{\dd}{\dd\ze} )^2 (1*\hat\ph*1*\hat\psi)$.
}

\sm

\begin{figure}

\begin{center}

\epsfig{file=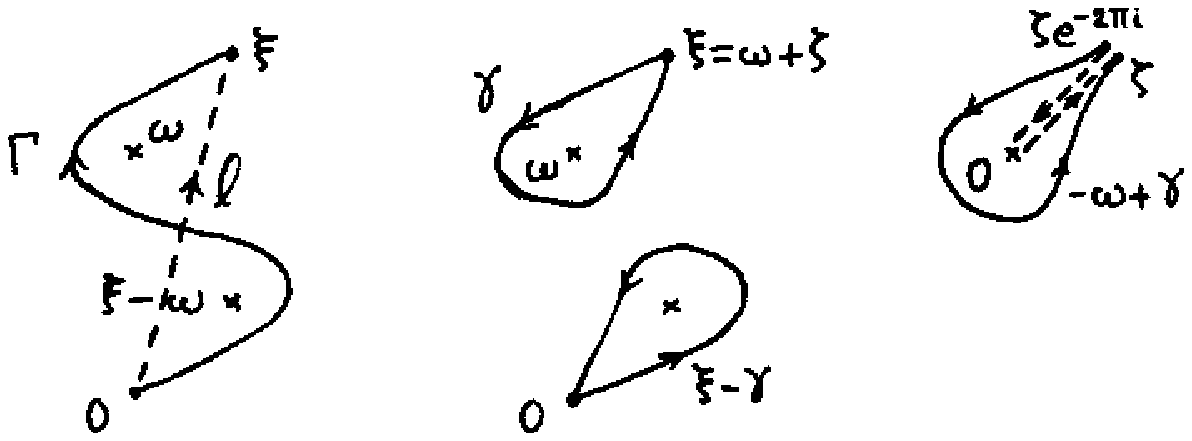,height=4cm,angle = 0}

\end{center}

\vspace{-.5cm}

\caption{\label{figsimplestconv} Derivation of formula~\refeq{eqidentLeib}.}

\vspace{.5cm}

\end{figure}

Hence, to show that 
$\De^+_\om(\hat\ph*\hat\psi) = (a\,\de + \hat\Phi)*\hat\psi +
\hat\ph*(b\,\de + \hat\Psi)$, 
we just need to check that
\begin{equation}	\label{eqidentLeib}
\hat\ph*\hat\psi(\om+\ze) - \hat\ph*\hat\psi(\om+\ze\,\ee^{-2\pi\I}) = 
a\hat\psi(\ze) + \hat\Phi*\hat\psi(\ze) 
+ b\hat\ph(\ze) + \hat\ph*\hat\Psi(\ze).
\end{equation}
For $|\ze|\le\pi$, the first term in the \lhs\
is~$\cont_\ell \hat\ph*\hat\psi (\xi)$, given by~\refeq{eqconvell} with
$\xi=\om+\ze$, while the second term is
$\cont_\Ga \hat\ph*\hat\psi (\xi) = \int_\Ga
\hat\ph(\xi_1)\hat\psi(\xi-\xi_1)\,\dd\xi_1$, 
with a symmetrically contractile path~$\Ga$ as shown on
Figure~\ref{figsimplestconv}.
The difference is thus 
$(\int_\ga-\int_{\ga'})\hat\ph(\xi_1)\hat\psi(\xi-\xi_1)\,\dd\xi_1$, where
the path~$\ga$ reaches~$\xi$, turning once anticlockwise around~$\om$, having
started from~$\xi$ (or rather from~$\om+\ze\,\ee^{-2\pi\I}$), and where
$\ga'=\ze-\ga$.
With the change of variable $\xi_1\mapsto\xi-\xi_1$ in the second integral, we
get
$$
\hat\ph*\hat\psi(\om+\ze) - \hat\ph*\hat\psi(\om+\ze\,\ee^{-2\pi\I}) = 
\int_\ga \hat\ph(\xi_1)\hat\psi(\xi-\xi_1)\,\dd\xi_1
+ \int_\ga \hat\psi(\xi_1)\hat\ph(\xi-\xi_1)\,\dd\xi_1
$$
and the identity~\refeq{eqidentLeib} will follow if we prove that
$\int_\ga \hat\ph(\xi_1)\hat\psi(\xi-\xi_1)\,\dd\xi_1 = 
a\hat\psi(\ze) + \hat\Phi*\hat\psi(\ze)$ (with a symmetric formula for the
second integral).

\sm

With the change of variable $\xi_1\mapsto \ze_1=\xi_1-\om$, this integral can be
written
$$
\frac{1}{2\pi\I} \int_{-\om+\ga}
\left(\frac{a}{\ze_1} + \hat\Phi(\ze_1)\log\ze_1 + \mathrm{reg}(\ze_1) \right)
\hat\psi(\ze-\ze_1)\,\dd\ze_1.
$$
The conclusion follows, since $\ze_1$ and~$\ze-\ze_1$ stay in a neighbourhood of
the origin where~$\hat\Phi$, $\mathrm{reg}$ and~$\hat\psi$ are holomorphic: the
residuum formula takes care of the simple pole and the Cauchy theorem cancels
the contribution of~$\mathrm{reg}(\ze_1)$, while the contribution of the
logarithm can be computed by collapsing the path~$-\om+\ga$ onto the segment
$[\ze\,\ee^{-2\pi\I},0]$ followed by~$[0,\ze]$.

\sm

\begin{figure}

\begin{center}

\epsfig{file=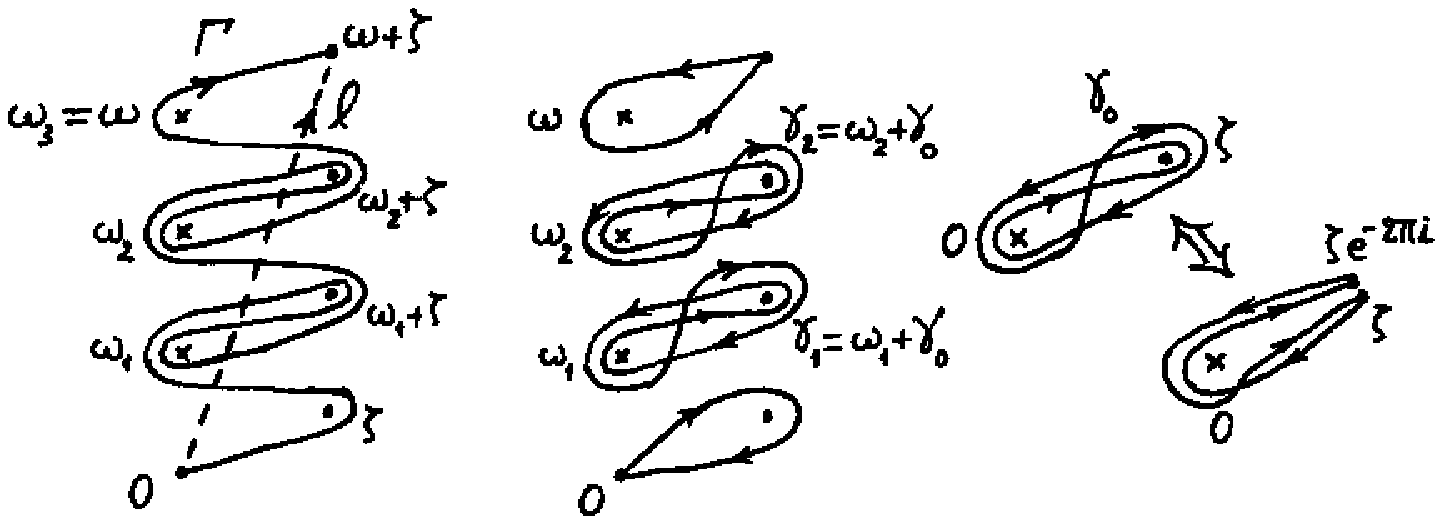,height=5.1cm,angle = 0}

\end{center}

\vspace{-.5cm}

\caption{\label{figcompliqconv} Derivation of formula~\refeq{eqidentLeibplus}.}

\end{figure}

Now let $m\ge2$ and consider $\om=2\pi\I\, m$.
Assume $\hat\chi_1=\hat\ph$ and $\hat\chi_2=\hat\psi$ with
$\De^+_{\om_j}\hat\ph = a_j\,\de + \hat\Phi_j$ and $\De^+_{\om_j}\hat\psi =
b_j\,\de + \hat\Psi_j$ for $\om_j = 2\pi\I\, j$, $j\in\{1,\ldots,m\}$.
We must prove
$$
\De^+_{2\pi\I\,m} (\hat\ph*\hat\psi) = (\De^+_{2\pi\I\,m}\hat\ph)*\hat\psi +
\sum_{m_1+m_2=m} (\De^+_{2\pi\I\,m_1}\hat\ph)*(\De^+_{2\pi\I\,m_2}\hat\psi)
+ \hat\ph*(\De^+_{2\pi\I\,m}\hat\psi).
$$
This time, to simplify the computations, we begin with the case where all the
constants~$a_j$ and~$b_j$ vanish.
This means that, considering $\xi=\om+\ze$ such that $|\ze|<\pi$ and
$\ell=[0,\xi] \subset \C \setminus 2\pi\I\,\Z^*$,
instead of~\refeq{eqidentLeib} we now must show
\begin{multline}	\label{eqidentLeibplus}
\hat\ph*\hat\psi(\om+\ze) - \hat\ph*\hat\psi(\om+\ze\,\ee^{-2\pi\I}) = 
\hat\Phi_m*\hat\psi(\ze) +
\sum_{m_1+m_2=m} \hat\Phi_{m_1}*\hat\Psi_{m_2}
\,+ \hat\ph*\hat\Psi_m(\ze),
\end{multline}
where the first term in the \lhs\ is $\cont_\ell\hat\ph*\hat\psi (\xi)$, still
given by~\refeq{eqconvell}, and the second term is
$\cont_\Ga\hat\ph*\hat\psi (\xi) = \int_\Ga
\hat\ph(\xi_1)\hat\psi(\xi-\xi_1)\,\dd\xi_1$, 
with a more complicated symmetrically contractile path~$\Ga$, as shown on
Figure~\ref{figcompliqconv}.

\sm

The difference can thus be decomposed as the sum of~$m+1$ terms, with the same
extreme terms as in the case $m=1$, which thus yield $\hat\Phi_m*\hat\psi$
and~$\hat\ph*\hat\Psi_m$ as in the first part of the proof, and with
intermediary terms
$\int_{\ga_{m_1}}\hat\ph(\xi_1)\hat\psi(\xi-\xi_1)\,\dd\xi_1$
for $m_1\in\{1,\ldots,m-1\}$, with paths $\ga_{m_1} = \om_{m_1}+\ga_0$ shown on
Figure~\ref{figcompliqconv}. 
Each intermediary term can be written
$$
\frac{1}{2\pi\I} \int_{\ga_0} \hat\ph(\om_{m_1}+\ze_1) 
\left( 
\hat\Psi_{m_2}(\ze-\ze_1)\log(\ze-\ze_1) + \mathrm{reg}(\ze-\ze_1)
\right)\,\dd\ze_1
$$
with $m_2=m-m_1$.
Collapsing the path~$\ga_0$ as indicated on Figure~\ref{figcompliqconv}, we see
that this term is nothing but
$\int_{\ti\ga_0} \hat\ph(\om_{m_1}+\ze_1) 
\hat\Psi_{m_2}(\ze-\ze_1) \, \dd\ze_1$,
where the path~$\ti\ga_0$ is identical to the path~$-\om+\ga$ of
Figure~\ref{figsimplestconv}; such an integral was already computed in the first
part of the proof: it is $\hat\Phi_{m_1}*\hat\Psi_{m_2}(\ze)$, which yields
formula~\refeq{eqidentLeibplus}. 

\sm

We end with the general case, with arbitrary constants~$a_j$ and~$b_j$. 
As alluded to in footnote~\ref{footsimplify}, it is sufficient to write
$\hat\ph*\hat\psi$ as the second derivative of $(1*\hat\ph)*(1*\hat\psi)$ and to
observe that
$(1*\hat\ph)(\om_j+\ze) = \frac{1}{2\pi\I}\Bigl(
\bigl(a_j+(1*\hat\Phi_j)(\ze)\bigr)\log\ze + \mathrm{reg}(\ze)
\Bigr)$, and similarly for $(1*\hat\psi)(\om_j+\ze)$, by
integrating~\refeq{eqsingphpsi} (we call indifferently $\mathrm{reg}(\ze)$ all
the regular germs that appear and that do not affect the final result; we used
the fact that, in the integration by parts, the primitives of
$(1*\hat\Phi_j)(\ze)\frac{\dd}{\dd\ze}(\log\ze)$ are regular).
Of course, convolution is commutative and $1*1=\ze$,
the formula for the case of vanishing residua thus yields
\begin{multline*}
\ze*\hat\ph*\hat\psi(\om+\ze) = 
\frac{1}{2\pi\I} \sum_{m_1+m_2=m} \Bigl( \bigl(
a_{m_1} b_{m_2} \ze + a_{m_1} \ze*\hat\Psi_{m_2}(\ze) \\
+ b_{m_2} \ze*\hat\Phi_{m_1}(\ze) + \ze*\hat\Phi_{m_1}*\hat\Psi_{m_2}(\ze)
\bigr) \log\ze + \mathrm{reg}(\ze) \Bigr),
\end{multline*}
with a sum extending to $(0,m)$ and~$(m,0)$, with the convention $a_0=b_0=0$.
The conclusion follows by differentiating twice (observing that
$(1*\hat A)(\ze)/\ze$ is regular for whatever regular germ~$\hat A$).
\eopf

\sm

It remains to prove Proposition~\ref{propsimpleAlienD}. For the moment, we
content ourselves with indicating a relation between the operators~$\De_\om$
and~$\De^+_\om$ (which will follow from the choice of the
weights~$\frac{p(\eps)!q(\eps)!}{m!}$), the idea being that according to
Lemma~\ref{lemDeomplus} the operators~$\De^+_{\pm2\pi\I\,m}$ are the homogeneous components
of two formal automorphisms of graded algebra and that, according to the next
lemma, the operators~$\De_{\pm2\pi\I\,m}$ are the homogeneous components of the
logarithms of these automorphims.

\begin{lemma}	\label{lemDeLogDeplus}
For each $m\in\N^*$ and $\om = \pm2\pi\I\,m$, we have
\begin{equation}	\label{eqDeLogDeplus}
\De_\om = \sum_{1\le r\le m} \tfrac{(-1)^{r-1}}{r}
\sum_{\substack{m_1,\ldots,m_r\ge1 \\ m_1+\cdots+m_r=m}}
\De^+_{\om_{m_1}}\cdots\De^+_{\om_{m_r}}
\end{equation}
with the notation $\om_j = \frac{j}{m}\om$.
\end{lemma}
We thus have
\begin{gather*}
\De_{\om_1} = \De^+_{\om_1} \\[1.5ex]
\De_{\om_2} = \De^+_{\om_2} - \tfrac{1}{2}\De^+_{\om_1}\De^+_{\om_1} \\[1.5ex]
\De_{\om_3} = \De^+_{\om_3} - \tfrac{1}{2}\left( \De^+_{\om_2}\De^+_{\om_1}
+ \De^+_{\om_1}\De^+_{\om_2} \right) 
+ \tfrac{1}{3}\De^+_{\om_1}\De^+_{\om_1}\De^+_{\om_1} \\[.5ex]
\vdots
\end{gather*}
with $\om_m=2\pi\I\,m$ for all~$m\ge1$ or $\om_m=-2\pi\I\,m$ for all~$m\ge1$.

\sm

{\em The proof of Lemma~\ref{lemDeLogDeplus} and the way it implies
formula~\refeq{eqLeibniz} of Proposition~\ref{propsimpleAlienD} are deferred to
the end of the next section.}
Formula~\refeq{eqLeibnizmult} follows by passage to the formal model.
The rest of Proposition~\ref{propsimpleAlienD} is easy:
formula~\refeq{eqcommutnatur} is best seen in the convolutive model;
formula~\refeq{eqaliencompos} follows from~\refeq{eqdotaliencompos}, which can
be checked in the convolutive model, with the help of a Taylor expansion
analogous to~\refeq{eqcomposconvol}.

\sm

We end this section by mentioning that substitution of a simple resurgent
function~$\ti\psi$ without constant term into a convergent series
$C(w)\in\C\{w\}$ gives rise to a simple resurgent function $C\circ\ti\psi$, the
alien derivatives of which are given by
\begin{equation}	\label{eqaliensubst}
\De_\om (C\circ\ti\psi) = (C'\circ\ti\psi) \De_\om\ti\psi.
\end{equation}


\subsection{Bridge equation for non-degenerate parabolic germs}	\label{secParBridge}


We can now state \'Ecalle's result for tangent-to-identity holomorphic germs,
which is at the origin of the name ``Resurgence''.

\begin{thm}	\label{thmParBridge}
Let $f$ be a non-degenerate parabolic germ at infinity with vanishing resiter
($\rho=0$), and let $\ti v(z)=z+\ti\psi(z)$ and $\ti u(z)=z+\ti\ph(z)$ be the
formal solutions of equations~\refeq{eqnAbelInv} and~\refeq{eqnAbel} with
$\ti\ph,\ti\psi\in z\ii\C[[z\ii]]$.
Then $\ti\ph,\ti\psi\in\wtRsimp$, and there exists a sequence of complex numbers
$\{A_\om\}_{\om\in2\pi\I\,\Z^*}$ such that, for each $\om\in2\pi\I\,\Z^*$,
\begin{equation}	\label{eqnBridge}
\De_\om \ti u  = A_\om \pa \ti u = A_\om \bigl(1 + \pa\ti\ph(z)\bigr),
\quad \De_\om \ti v  = -A_\om \, \ee^{-\om(\ti v(z)-z)} = -A_\om\,\ee^{-\om\ti\psi(z)}.
\end{equation}
\end{thm}

\sm

Here the notation was slightly extended \wrt\ Definition~\ref{defDeom}: 
$\De_\om\ti u$ is to be understood as
being equal to~$\De_\om\ti\ph$ (since $\ti u(z)-\ti\ph(z)=z$ is convergent:
the difference offers no singularity to be measured by any alien derivation).
In the convolutive model, this amounts to setting $\De_\om \de'=0$ (we already
had $\De_\om \de=0$).
Similarly, $\De_\om\ti v = \De_\om\ti\psi$.
The translation of~\refeq{eqnBridge} in the convolutive model is thus
\begin{equation}	\label{eqnconvBridge}
\De_\om\hat\ph = A_\om\,\de - A_\om \ze\hat\ph(\ze),
\quad
\De_\om\hat\psi = -A_\om\left(\de - \om\hat\psi
+\tfrac{1}{2!}\om^2\hat\psi^{*2} 
-\tfrac{1}{3!}\om^3\hat\psi^{*3} +\cdots \right).
\end{equation}
The existence of such relations is the resurgent phenomenon: $\hat\ph(\ze)$
or~$\hat\psi(\ze)$, which are a holomorphic germs at the origin, reappear in a
disguised form at the singularities of their analytic continuation, when
singularities are measured in an appropriate way.

\sm

Equation~\refeq{eqnBridge} is called the ``Bridge equation'', because the first
equation may be viewed as a bridge between alien calculus ($\De_\om$) and
ordinary calculus ($\pa=\frac{\dd}{\dd z}$) in the case of~$\ti u$.
Notice that it is possible to iterate these equations to compute the successive
alien derivatives $\De_{\om_r}\cdots\De_{\om_1}\ti u$ or
$\De_{\om_r}\cdots\De_{\om_1}\ti v$, since we know how the alien derivations
interact with~$\pa$ (see formula~\refeq{eqcommutnatur}) and with exponentiation
(using~\refeq{eqaliensubst}).
The computation is simpler with dotted alien derivations: we get
\begin{equation}	\label{eqiterBridgeu}
\dDe{\om_r}\cdots\dDe{\om_1} \ti u = 
D_{\om_1} \cdots D_{\om_r} \ti u, \qquad
D_\om = A_{\om}\,\ee^{-\om z}\,\pa
\end{equation}
(beware of the non-commutation of~$\pa$ and multiplication by~$\ee^{-\om z}$:
the vector fields $D_{\om_j}$ do not commute one with the other, but they do
commute with the dotted alien derivations, hence the reversal of order) and
\begin{equation}	\label{eqiterBridgev}
\dDe{\om_r}\cdots\dDe{\om_1} \ti v = 
- A_{\om_1}\cdots A_{\om_r} \,
\om_1 (\om_1+\om_2)\cdots (\om_1+\cdots+\om_{r-1})
\,\ee^{-(\om_1+\cdots+\om_{r})\ti v}.
\end{equation}
The collection of~$\{A_\om,\;\om\in2\pi\I\,\Z^*\}$ is thus sufficient to
describe the whole singular behaviour of~$\hat\ph$ and~$\hat\psi$.

\sm

\noindent {\em Proof of Theorem~\ref{thmParBridge}.}
The fact that $\hat\ph$ and~$\hat\psi$ are simple resurgent functions was
already alluded to in the previous section (after Proposition~\ref{propsimplestab}).
With our extension of~$\De_\om$ to the space $\ID+\wtRsimp$ (which contains
$f=\ID+1+a$ and~$\ti u=\ID+\ti\ph$), equation~\refeq{eqaliencompos} yields
$$
\De_\om(f\circ\ti u) = \ee^{-\om\ti\ph}\bigl( (\De_\om f)\circ\ti u \bigr)
+ ( \pa f\circ\ti u ) \De_\om\ti u.
$$
But $\De_\om f=0$ because~$a$ is a convergent power series (its Borel transform
has no singularity); on the other hand, $\De_\om(f\circ\ti u) = \De_\om (\ti
u(z+1)) = (\De_\om\ti u)(z+1)$, hence
$$
(\De_\om\ti u)(z+1) = (\pa f\circ\ti u)(z) \De_\om\ti u(z).
$$
The equation $Y(z+1) = (\pa f\circ\ti u)(z) Y(z)$ is nothing but the
linearization of the equation $u(z+1) = f\circ u(z)$ around the solution~$\ti
u$, and we know a solution of this linear difference equation, namely~$\pa\ti u$
(because~$\pa$ is also a derivation which commutes with
$\ti\chi(z)\mapsto\ti\chi(z+1)$);
moreover $\pa\ti u\in 1+z\ii\C[[z\ii]]$ admits a multiplicative inverse.
As a consequence, the formal series $A_\om = \frac{1}{\pa\ti u}\De_\om\ti u$ must be
invariant by $z\mapsto z+1$, hence constant.

\sm

The extension of equation~\refeq{eqdotaliencompos} to $\ID+\wtRsimp$ yields
$$
\dDeom(\ti u\circ\ti v) = \bigl( \dDeom\ti u \bigr) \circ\ti v 
+ (\pa\ti u\circ\ti v) \dDeom\ti v
$$
and this expression must vanish, since $\ti u\circ\ti v(z)=z$, hence
$$
\bigl( \dDeom\ti v \bigr) \circ\ti v\ii = - \frac{1}{\pa\ti u} \dDeom\ti u
= - A_\om \,\ee^{-\om z},
$$
which implies $\dDeom\ti v = - A_\om \,\ee^{-\om\ti v}$.
\eopf

\subsub{\'Ecalle's analytic invariants}

The coefficients~$A_\om$ in the Bridge equation are called ``\'Ecalle's analytic
invariants''.
Observe that $A_\om$ can be defined as the coefficient of~$\de$
in~$\De_\om\hat\ph$ (or of~$-\De_\om\hat\psi$), which is an average of the
residua of $2^{m-1}$ branchs of~$\hat\ph$ (or of~$-\hat\psi$) if $\om =
\pm 2\pi\I\,m$.
These coefficients form a complete system of analytic invariants in the following sense:
\begin{prop}
Two non-degenerate parabolic germs~$f_1$ and~$f_2$ with vanishing resiter are
analytically conjugate if and only if there exists $c\in\C$ such that $A_\om\ss2
= \ee^{-\om c} A_\om\ss1$ for all $\om\in 2\pi\I\,\Z^*$.
\end{prop}

The proof of one implication is easy. Suppose that $f_2 = h\ii\circ f_1\circ h$
with $h\in\ID+\C\{z\ii\}$; we shall prove that
$A_\om\ss2 = \ee^{-\om c} A_\om\ss1$, where~$c$ is defined by $h(z) = z + c + \cO(z\ii)$.
By the same argument as in footnote~\ref{footanalyticinv},
since $\ti v_1\circ h(z) = z + c + \cO(z\ii)$ is a formal solution of Abel's
equation for~$f_2$, Proposition~\ref{propformalconj} implies the existence of
$c'\in\C$ such that $\ti v_1\circ h = \ti v_2 + c'$; we get $c'=c$ since $\ti
v_2(z) = z + \cO(z\ii)$.
Using the chain rule for alien derivations, since $h$ is convergent, we find
$$
\dDeom \ti v_2 = \bigl( \dDeom\ti v_1 \bigr) \circ h =
- A_\om\ss1 \,\ee^{-\om \, \ti v_1\circ h} = 
- A_\om\ss1 \,\ee^{-\om c}\,\ee^{-\om\ti v_2}
$$
hence $A_\om\ss2 = A_\om\ss1 \,\ee^{-\om c}$.
(We could have expoited the relation $h\circ\ti u_2(z) = \ti u_1(z+c)$ equally).

\sm

The direct verification of the other implication requires an extra work and we
shall not give all the details. 
Suppose that $A_\om\ss 2 = \ee^{-\om c} A_\om\ss 1$ for all $\om\in
2\pi\I\,\Z^*$ and consider $\ti h = \ti u_1\circ(c+\ti v_2) \in\ID+\wtRsimp$,
which establishes a formal conjugacy between~$f_1$ and~$f_2$.
A computation similar to the above yields 
$$
\De_\om\ti h = \left( A_\om\ss1\,\ee^{-\om c} - A_\om\ss2 \right) 
\ee^{-\om\ti v_2} \, \pa\ti u_1\circ (c+\ti v_2)
= 0
$$
for all~$\om$, thus $\hat h=\cB(\ti h-\ID)$ has no singularity at all, and one can
deduce from Theorem~\ref{thmresur} that this entire function has at most
exponential growth in the non-vertical directions.
However, an extra argument is needed to make sure that~$\hat h$ has at most
exponential growth with bounded type in all directions, including the imaginary
axis, which is sufficient to conclude that~$\ti h$ is convergent.\footnote{
\label{footiter}
Here is another interesting property of the $A_\om$'s:
since the normal form $f_0(z)=z+1$ is the time-$1$ map of~$\frac{\pa}{\pa z}$,
we can define its $t^\text{th}$ power of iteration by
$f_0\sst t(z) = z+t$ for any $t\in\C$, and by formal conjugacy we retrieve the
``$t^\text{th}$ power of iteration'' of~$f$:
$\ti f\sst t = \ti u \circ f_0\sst t \circ \ti v \in \ID+\wtRsimp$, 
which is the unique $\ti f\sst t(z) \in z + t + z^{-2}\C[[z\ii]]$ such that $\ti f\ss
t\circ f = f \circ \ti f\sst t$
(in other words, we embed~$f$ in a formal one-parameter group $\{\ti f\sst t,\;
t\in\C\}$, which is generated by the formal vector field 
$\pa\ti u\circ\ti v(z)\frac{\pa}{\pa z} = \frac{1}{\pa\ti v(z)}\frac{\pa}{\pa z}$).
For a given $t\in\C$, $\ti f\sst t$ is always resurgent but usually divergent,
unless $t\in\Z$, however it may happen that some values of~$t$ give rise to a
convergent trnasformation (and this does not imply the convergence of the formal
infinitesimal generator $\frac{1}{\pa\ti v(z)}\frac{\pa}{\pa z}$).
{\em The set of all $t\in\C$ such that~$\ti f\sst t$ is convergent is the group
$\frac{1}{q}\Z$, where $q\in\N^*$ is determined by the condition
$\om\notin2\pi\I\,q\Z^* \Rightarrow A_\om=0$}, which is consistent with the
relations
$\dDeom\ti f\sst t = A_\om\left(\ee^{-\om t}-1\right)\ee^{-\om\ti v}
\frac{1}{\pa\ti v} \pa\ti f\sst t$
which follow from computations analogous to the previous ones.
This corresponds to the fact that the $1$-periodic
functions~$\chi^{\text{low,up}}$ encoding the horn maps may admit a period~$q$
larger than~$1$.
}

\sm

The extra argument that we have not given is related to the existence of
constraints on the growth of~$\hat\ph$ and~$\hat\psi$ along the lines which are
parallel to the imaginary axis, which imply constraints on the growth of
the numbers~$|A_{\pm2\pi\I m}|$ as $m\to\infty$.
This can be obtained by a fine analysis in the Borel plane. 
Another approach consists in relating \'Ecalle's analytic invariants and the horn
maps:
there is a one-to-one correspondence between
$\{ A_\om, \; \om\in 2\pi\I\,\N^* \}$ and $\{ B_m, \; m\in\N^* \}$ on the one hand,
and between
$\{ A_\om, \; \om\in -2\pi\I\,\N^* \}$ and $\{ B_{-m}, \; m\in\N^* \}$ on the other hand.
This will be the subject of the next paragraphs.
We shall see that the growth contraint on the~$|A_\om|$ amounts exactly to the
convergence of the Fourier series~\refeq{eqFourierLOW}--\refeq{eqFourierUP} with
some $\tau_0>0$.

\subsub{Relation with the horn maps}

Let us give the recipe before trying to justify it:
if we work in the graded algebras 
\begin{equation}	\label{eqdefExpRes}
\wtRsimp[[\ee^{-2\pi\I\,z}]] =
\underset{\om\in2\pi\I\,\N}{\oplus} \, \ee^{-\om z} \wtRsimp, \qquad
\text{resp.}\quad
\wtRsimp[[\ee^{2\pi\I\,z}]] =
\underset{\om\in-2\pi\I\,\N}{\oplus} \, \ee^{-\om z} \wtRsimp
\end{equation}
(so as to give a meaning to the~$\dDeom$'s as internal operators, which
commute with the multiplication by~$\ee^{-\om_0 z}$ for any~$\om_0$ and are
$\om$-homogeneous\footnote{ 
\label{footextendgraded}
We simply mean that, if $\om=2\pi\I\,m$ with $m\ge1$ for instance, there is a unique
linear operator
$$
\dDeom:\, \wtRsimp[[\ee^{-2\pi\I\,z}]] \to \wtRsimp[[\ee^{-2\pi\I\,z}]]
$$
which extends the operator~$\dDeom$ previously defined in~$\wtRsimp$ only and
which commutes with multiplication by~$\ee^{-\om_0 z}$ for any~$\om_0$.
This operator is $\om$-homogeneous in the sense that it sends the space of
$\om_0$-homogeneous elements in the set of $(\om_0+\om)$-homogeneous elements:
$$
\dDeom \bigl( \ee^{-\om_0 z} \wtRsimp \bigr) \subset
\ee^{-(\om_0+\om) z} \wtRsimp.
$$
})
and define the ``directional alien derivations'' as the two operators
$$
\De_{\I\R^+} = \sum_{\om\in 2\pi\I\,\N^*} \dDeom, \qquad
\text{resp.}\quad
\De_{-\I\R^+} = \sum_{\om\in -2\pi\I\,\N^*} \dDeom,
$$
then the Borel-Laplace summation in a direction~$\th$ slightly smaller
than~$\pm\pi/2$ is equivalent to the composition of the Borel-Laplace summation
in a direction~$\th'$ slightly larger than~$\pm\pi/2$ with the exponential of the
directional alien derivation associated with $\ee^{\pm\I\pi/2}\R^+$:
$$
\cL^\th\circ\cB \sim
\cL^{\th'}\circ\cB\circ\exp\left(\De_{\pm\I\R^+}\right)
= \cL^{\th'}\circ\cB\circ\Biggl( \sum_{r\ge0} \frac{1}{r!}
\De_{\pm\I\R^+}^r \Biggr).
$$
The symbol $\sim$ here means that the operators on both sides should
coincide when applied to formal sums (``transseries'')
$\ph = \sum \ee^{-\om z}\,\ti\ph_\om(z)$
and yield
\begin{equation}	\label{eqGradedFormalPass} 
\sum \ee^{-\om z} \, \cL^\th\hat\ph_\om = \sum \ee^{-\om z} \, \cL^{\th'}\hat\psi_\om,
\qquad
\ti\psi_\om = \sum_{r\ge0} \frac{1}{r!} \sum_{\om_0+\om_1+\cdots+\om_r=\om}
\De_{\om_r}\cdots\De_{\om_1} \ti\ph_{\om_0},
\end{equation}
whenever both sides of the first equation in~\refeq{eqGradedFormalPass} have an
analytical meaning (with all indices~$\om_j$ in~$2\pi\I\,\N^*$, or all~$\om_j$
in~$-2\pi\I\,\N^*$, with the exception of~$\om_0$ and~$\om$ which may vanish;
observe that each~$\ti\psi_\om$ is defined by a finite sum).

\sm

In our case, the recipe yields
\begin{alignat*}{3}
u^+ &= \cL^{\th'} \cB \left[
\exp\left(\De_{\I\R^+}\right)\ti u \right], \quad &
v^+ &= \cL^{\th'} \cB \left[
\exp\left(\De_{\I\R^+}\right)\ti v \right], 
\qquad && \IM z < -\tau_0, \\[1ex]
u^- &= \cL^{\pi+\th'} \cB \left[
\exp\left(\De_{-\I\R^+}\right)\ti u \right], \quad &
v^- &= \cL^{\pi+\th'} \cB \left[
\exp\left(\De_{-\I\R^+}\right)\ti v \right],
\qquad && \IM z > \tau_0,
\end{alignat*}
with $\th'$ slightly larger than~$\pi/2$.
To interpret and understand this, one can remember the computation which was made in
Section~\ref{secdiffeq} in the case of the linear difference
equation~\refeq{eqnph}: 
we had $\ti\ph(z+1)-\ti\ph(z) = a(z) \in z^{-2}\C\{z\ii\}$, hence
$\De_\om\ti\ph = A_\om$ had to be constant (indeed, $A_\om=-2\pi\I\,\hat a(\om)$),
thus $\De_{\I\R^+}\ti\ph = \sum_{\om\in2\pi\I\,\N^*} A_\om\,\ee^{-\om z}$ and 
$\De_{\I\R^+}^r\ti\ph = 0$ for $r\ge2$,
formula~\refeq{eqlinearStokes} (which was nothing but the residuum formula) can
thus be interpreted as
$$
\cL^{\th}\cB\ti\ph = 
\ph^+ = \ph^- + \cL^{\th'}\cB \, \De_{\I\R^+}\ti\ph = 
\cL^{\th'}\cB \exp\left(\De_{\I\R^+}\right)\ti\ph
\quad \text{in $\{ \IM z < -\tau \}$}.
$$
In the case of~$\ti u$ or~$\ti v$, we need to take into account the action of
$\De_{\pm\I\R^+}^r$ for $r\ge2$, but formulas~\refeq{eqiterBridgeu}
and~\refeq{eqiterBridgev} allow us to perform the calculation.

\sm

Let us begin with~$\ti u$: we have $\exp\left(\De_{\pm\I\R^+}\right)\ti u = 
\exp\left(D_{\pm\I\R^+}\right)\ti u$, with
$$
D_{\I\R^+} = \sum_{\om\in 2\pi\I\,\N^*} A_\om\,\ee^{-\om z}\frac{\dd}{\dd z}, \qquad
D_{-\I\R^+} = \sum_{\om\in -2\pi\I\,\N^*} A_\om\,\ee^{-\om z}\frac{\dd}{\dd z}.
$$
The operator $D_{\pm\I\R^+}$ is a derivation of our graded algebra, \ie a
(formal) vector field; its exponential is thus an automorphism, which can be
represented as a substitution operator:
$$
\exp\left(D_{\pm\I\R^+}\right)\ti u = \ti u \circ P_{\pm\I\R^+}, \qquad
P_{\pm\I\R^+} = \exp\left(D_{\pm\I\R^+}\right)\ID.
$$
In fact, $P_{\pm\I\R^+}$ is the time-$1$ map of~$D_{\pm\I\R^+}$.
Writing this formal vector field in the coordinate
$w_\pm = \ee^{\pm2\pi\I z}$, as
$D_\pm = \pm2\pi\I \sum A_{\pm2\pi\I m} w^{m+1} \frac{\dd}{\dd w}$,
we get a tangent-to-identity transformation, which we can write as
$$
\exp(D_\pm)\,:\ens 
w \mapsto  w\,\ee^{\pm2\pi\I \sum P_{\pm2\pi\I m} w^m}.
$$
This relation determines the coefficients~$P_\om$ in terms of the~$A_\om$'s, so
that
$$
P_{\pm\I\R^+}(z) = z + \sum_{\om\in\pm2\pi\I\,\N^*} P_\om\,\ee^{-\om z}.
$$
The final result is
\begin{equation}	\label{eqpassageu}
u^+ = u^- \circ P_{\,\I\R^+}  \ens\text{ in }\; \{ \IM z < -\tau_0 \}, \qquad
u^- = u^+ \circ P_{-\I\R^+}  \ens\text{ in }\; \{ \IM z > \tau_0 \}.
\end{equation}

\sm

The computation with~$\ti v$ is more direct:
formula~\refeq{eqiterBridgev} shows that 
$$
\exp\left(\De_{\pm\I\R^+}\right)\ti v = \ti v - \sum_{r\ge1} \frac{1}{r!}
\sum_{\om_1,\ldots,\om_r\in\pm2\pi\I\,N^*} A_{\om_1}\cdots A_{\om_r} \,
\Ga_{\om_1\cdots\om_r} \,\ee^{-(\om_1+\cdots+\om_{r})\ti v}
= Q_{\pm\I\R^+} \circ\ti v
$$
with the notation $\Ga_{\om_1}=1$, $\Ga_{\om_1\cdots\om_r}=
\om_1(\om+\om_2)\cdots(\om_1+\om_2+\cdots+\om_{r-1})$, and
\begin{equation}	\label{eqAQ}
Q_{\pm\I\R^+}(z) = z  + \sum_{\om\in\pm2\pi\I\,\N^*} Q_\om\,\ee^{-\om z},
\qquad
Q_\om = - \sum_{r\ge1} \frac{1}{r!}
\sum_{\substack{\om_1,\ldots,\om_r\in\pm2\pi\I\,N^*\\\om_1+\cdots+\om_{r}=\om}} 
\Ga_{\om_1\cdots\om_r} A_{\om_1}\cdots A_{\om_r} .
\end{equation}
The upshot is
\begin{equation}	\label{eqpassagev}
v^+ = Q_{\I\R^+} \circ v^-  \ens\text{ in }\; \{ \IM z < -\tau_0 \}, \qquad
v^- = Q_{-\I\R^+} \circ v^+  \ens\text{ in }\; \{ \IM z > \tau_0 \}.
\end{equation}

\sm

The comparison of~\refeq{eqpassageu} and~\refeq{eqpassagev} shows that 
$Q_{\pm\I\R^+} = P_{\pm\I\R^+}\ii = v^+\circ u^-$ or~$v^-\circ u^+$ (according
to the half-plane under consideration); $Q_{\pm\I\R^+}$ can thus be obtained as
the time-$1$ map of~$-D_{\pm\I\R^+}$, and
the relation relating the $-A_\om$'s to the $Q_\om$'s in~\refeq{eqAQ} can thus be
paralleled by a relation relating the $A_\om$'s to the $P_\om$'s.
We arrive at

\begin{thm}	\label{thmRelationHorn}
We have
\begin{alignat}{3}
v^+\circ u^- & = Q_{\I\R^+} & \quad\text{with inverse}\quad
v^-\circ u^+ & = P_{\,\I\R^+} 
&& \quad\text{in}\ens \{\IM z<-\tau_0\}, \\[1ex]
v^+\circ u^- & = P_{-\I\R^+} & \quad\text{with inverse}\quad
v^-\circ u^+ & = Q_{-\I\R^+} 
&& \quad\text{in}\ens \{\IM z>\tau_0\},
\end{alignat}
where
$$
Q_{\pm\I\R^+} = \ID + \sum_{\om\in \pm 2\pi\I\,\N^*} Q_\om\, \ee_\om,
\qquad
P_{\pm\I\R^+} = \ID + \sum_{\om\in \pm 2\pi\I\,\N^*} P_\om\, \ee_\om,
$$
with the notation $\ee_\om(z) = \ee^{-\om z}$, and the coefficients $Q_\om, P_\om$
depend on the coefficients of the Bridge Equation according to the formulas
\begin{multline*}
Q_\om = - \sum_{r\ge1} \frac{1}{r!}
\sum_{\substack{\om_1,\ldots,\om_r\in\pm2\pi\I\,N^*\\\om_1+\cdots+\om_{r}=\om}} 
\Ga_{\om_1\cdots\om_r} A_{\om_1}\cdots A_{\om_r},\\
P_\om = \sum_{r\ge1} \frac{(-1)^{r-1}}{r!}
\sum_{\substack{\om_1,\ldots,\om_r\in\pm2\pi\I\,N^*\\\om_1+\cdots+\om_{r}=\om}} 
\Ga_{\om_1\cdots\om_r} A_{\om_1}\cdots A_{\om_r},
\end{multline*}
using the notation $\Ga_{\om_1}=1$, $\Ga_{\om_1\cdots\om_r}=
\om_1(\om+\om_2)\cdots(\om_1+\om_2+\cdots+\om_{r-1})$.
\end{thm}

The comparison with~\refeq{eqFourierLOW}--\refeq{eqFourierUP} now shows that
$Q_{\I\R^+} = \ID + \chi^{\text{\rm low}}$ and $P_{-\I\R^+}=\ID+\chi^{\text{\rm
up}}$, hence
$$
B_m = Q_{2\pi\I m}, \quad B_{-m} = P_{-2\pi\I m}, \qquad m\ge1.
$$
The fact that the $A_\om$'s constitute a complete set of analytic invariants can
be found again this way.

\sm

Observe that~$P_{\pm\I\R^+}$ or~$Q_{\pm\I\R^+}$ was obtained as a formal Fourier
series, being the time-$1$ map of the formal vector field~$D_{\pm\I\R^+}$
or~$-D_{\pm\I\R^+}$. However, when identified with~$v^+\circ u^-$ or~$v^-\circ
u^+$, these expansions prove to be convergent.
This is the growth constraint we were alluding to: the $A_\om$'s must be such
that $Q_{\pm\I\R^+}$ defined by~\refeq{eqAQ} be convergent.

\subsub{Alien derivations as components of the logarithm of the Stokes automorphism}

Let
$$
\cS_{\I\R^+} = \ID + \sum_{\om\in  2\pi\I\,\N^*}\ee^{-\om z} \De^+_\om, \quad
\cS_{-\I\R^+} = \ID + \sum_{\om\in -2\pi\I\,\N^*}\ee^{-\om z} \De^+_\om.
$$
Lemma~\ref{lemDeomplus} can be rephrased\footnote{
Here, as in footnote~\ref{footextendgraded}, we extend these operators
from~$\wtRsimp$ to~$\wtRsimp[[\ee^{\mp2\pi\I\,z}]]$ by declaring that they
commute with multiplication by~$\ee^{-\om_0 z}$ for any~$\om_0$. 
The Borel counterparts of the relations~\refeq{eqLeibnizplus} are obtained by
projecting the automorphism property
$$
\cS_{\pm\I\R^+} \bigl(\ti\chi_1\ti\chi_2\bigr) =
\bigl(\cS_{\pm\I\R^+} \ti\chi_1 \bigr) \bigl(\cS_{\pm\I\R^+} \ti\chi_2\bigr),
\qquad \ti\chi_1,\ti\chi_2 \in \wtRsimp[[\ee^{\mp2\pi\I\,z}]]
$$
onto the spaces $\ee^{-\om_0 z}\wtRsimp$, $\om_0\in\pm2\pi\I\,\N$.
}
by saying that $\cS_{\pm\I\R^+}$ is an automorphism of the graded algebra
$\wtRsimp[[\ee^{\mp2\pi\I\,z}]]$.
We close this chapter with two things.
\begin{enumerate}
\item
We shall indicate the proof of Lemma~\ref{lemDeLogDeplus}, the content of
which can be rephrased\footnote{
Formula~\refeq{eqDeLogDeplus} of Lemma~\ref{lemDeLogDeplus} simply expresses
the fact that the $\om$-homogeneous component of~$\De_{\pm\I\R^+}$ coincides
with that of
$$
\log\cS_{\pm\I\R^+} = \sum_{r\ge1} \frac{(-1)^{r-1}}{r} 
\biggl( \sum_{m\ge1} \ee^{-2\pi\I m z} \De^+_{\pm 2\pi\I m} \biggr)^r.
$$
} 
as identities between operators of $\wtRsimp[[\ee^{\mp2\pi\I\,z}]]$: {\em the
directional alien derivations satisfy
$$
\De_{\I\R^+} = \log \cS_{\I\R^+}, \qquad
\De_{-\I\R^+} = \log \cS_{-\I\R^+}.
$$
}
This was the only step missing in the proof of
Proposition~\ref{propsimpleAlienD} (the fact that the $\De_\om$'s are
derivations follows).
\item
We shall interpret the operators~$\cS_{\I\R^+}$ and~$\cS_{-\I\R^+}$ as ``Stokes
automorphisms'' (or ``passage automorphisms'', \cite{CNP}): {\em they correspond to
composing the Laplace transform in a direction with the inverse Laplace
transform in another direction.}
This will serve as a justification of the recipe which led us to
Theorem~\ref{thmRelationHorn} in the previous section, since the exponential
of~$\De_{\pm\I\R^+}$ will then appear as the link between Borel-Laplace summations in
different directions.
\end{enumerate}

\sm

Let us focus on the singular direction~$\I\R^+$ (the case of~$-\I\R^+$ is
analogous) and introduce a graded algebra which corresponds
to~$\wtRsimp[[\ee^{-2\pi\I\,z}]]$ via formal Borel transform.
For each $\om\in2\pi\I\,\N$, we define the
translation operator
$$
\tau_\om :\; c\,\de + \hat\ph \in\Rsimp \mapsto c\,\de_\om + \hat\ph_\om,
\qquad \hat\ph_\om(\ze) = \hat\ph(\ze-\om),
$$
where $\de_\om$ is a symbol to be identified with $\cB(\ee^{-\om z})$ and
$\hat\ph_\om$ is to be thought of as a holomorphic function based at~$\om$
(well-defined on~$\left]\om-2\pi\I,\om+2\pi\I\right[$, with multivalued analytic
continuation on the rest of the singular direction~$\I\R^+$); the range
of~$\tau_\om$ will be the space~$\hat\cR^\om$ of $\om$-homogeneous elements. We
consider
$$
\hat\cR = \underset{\om\in2\pi\I\,\N}{\oplus} \hat\cR^\om, \qquad
\hat\cR^\om = \tau_\om \left(\Rsimp\right),
$$
as a graded algebra by defining the (convolution) product of two homogeneous
elements to be
$$
\tau_{\om_1}\hat\chi_1 * \tau_{\om_2}\hat\chi_2 = 
\tau_{\om_1+\om_2} (\hat\chi_1*\hat\chi_2),
\qquad \hat\chi_1,\hat\chi_2\in\Rsimp,\; \om_1,\om_2\in2\pi\I\,\N.
$$
The operators~$\De_\om^+$ extend uniquely from $\hat\cR^0=\Rsimp$ to~$\hat\cR$ by
declaring that they commute with all translations~$\tau_{\om_0}$.
The Borel counterpart of~$\cS_{\I\R^+} $ is 
$$
\De^+ = \sum_{\om\in2\pi\I\,\N} \dDep\om \,:\; \hat\cR \;\to\; \hat\cR,
\qquad \text{with}\ens \dDep0=\ID \ens 
\text{and}\ens \dDep\om = \tau_\om\De_\om^+ \ens\text{for $\om\neq0$}.
$$
Each~$\dDep\om$ is thus seen as the $\om$-homogeneous component of the operator~$\De^+$.

\sm

We now introduce elementary operators $\dot\ell_+$, $\dot\ell_-$,~$\dot A$ and~$\mu$,
which will allow us to rewrite the definition~\refeq{eqdefDeomplus} as
\begin{equation}	\label{eqNewdefDeomplus}
\De_\om^+ = \tau_\om\ii 
\left( \dot A + \dot\ell_+ - \dot\ell_- \right)\dot\ell_+^{m-1}\mu, 
\qquad \om = 2\pi\I\,m.
\end{equation}
The first two ones are the {\em lateral continuation operators} and act in
$\ChcR = {\oplus}\, \chcR\om$, where~$\chcR\om$ is the space of all the
functions~$\chb\ph_\om$ which are holomorphic on~$\left]\om,\om+2\pi\I\right[$, can be
analytically continued along any path which avoids~$2\pi\I\,\Z$ and admit at
worse simple singularities.
The functions of~$\chcR\om$ are unambiguously determined
on~$\left]\om,\om+2\pi\I\right[$, but their analytic continuation gives rise to
various branchs.
Let $\la\in\left]0,\pi\right[$, and let~$\ga_+$, resp.~$\ga_-$, be the
semi-circular path which starts from~$\om+2\pi\I-\I\la$ and ends
at~$\om+2\pi\I+\I\la$, circumventing~$\om+2\pi\I$ to the right, resp.\ to the
the left. We define~$\dot\ell_+$ and~$\dot\ell_-$ to be the operators of
analytic continuation along~$\ga_+$ and~$\ga_-$:
$$
\dot\ell_\pm :\; \chb\ph_\om \in \chcR\om \mapsto 
\cont_{\ga_\pm} \chb\ph_\om \in \chcR{\om+2\pi\I}.
$$
These operators induce $2\pi\I$-homogeneous operators of~$\ChcR$, the difference
of which sends~$\ChcR$ in~$\hat\cR$: 
for any $\chb\ph_\om\in\chcR\om$, $(\dot\ell_+ - \dot\ell_-)\chb\ph_\om \in
\hat\cR^{\om+2\pi\I}$ is simply the variation of~$\Sing_{\om+2\pi\I}\chb\ph_\om$
translated by~$\tau_{\om+2\pi\I}$.
Denoting by~$\al_{\om+2\pi\I}(\chb\ph_\om)$ the residuum of~$\Sing_{\om+2\pi\I}\chb\ph_\om$,
we set
$$
\dot A:\; \chb\ph_\om \in \chcR\om \mapsto 
\al_{\om+2\pi\I}(\chb\ph)\,\de_{\om+2\pi\I} \in \hat\cR^{\om+2\pi\I}.
$$
The whole singularity~$\Sing_{\om+2\pi\I}\chb\ph_\om$ is thus determined by
$$
\tau_{\om+2\pi\I} \Sing_{\om+2\pi\I}\chb\ph_\om = \left(\dot A + 
\dot\ell_+ - \dot\ell_- \right) \chb\ph_\om \in \hat\cR^{\om+2\pi\I}.
$$
We define the last new elementary operator of the list by
$$
\mu:\; \hat\chi_\om = c\,\de_\om + \hat\ph_\om \in \hat\cR^\om \mapsto
\hat\ph_\om \in \chcR\om,
$$
\ie we forget the multiple of~$\de_\om$, retaining only~$\hat\ph_\om$ but forgetting that
this function is regular at~$\om$: the result is considered as element
of~$\chcR\om$.
Formula~\refeq{eqNewdefDeomplus} is now an
obvious translation of Definition~\ref{defDeom}.
We thus have
\begin{equation}	\label{eqDeplus}
\De^+ = \ID + \sum_{m\ge1}
\left( \dot A + \dot\ell_+ - \dot\ell_- \right)\dot\ell_+^{m-1}\mu.
\end{equation}

\sm

In this framework, we can also rephrase the part of Definition~\ref{defDeom}
concerning~$\De_\om$, or rather the dotted version of the operator (still in the
convolutive model):
$$
\dDeom = \tau_\om\De_\om =
\sum_{\eps_1,\ldots,\eps_{m-1}\in\{+,-\}} \frac{p(\eps)!q(\eps)!}{m!}
\left( \dot A + \dot\ell_+ - \dot\ell_- \right)
\dot\ell_{\eps_{m-1}}\cdots\dot\ell_{\eps_1} \mu,
\qquad \om = 2\pi\I\,m.
$$
We now give the proof of Lemma~\ref{lemDeLogDeplus}, which amounts to the
fact that the (Borel counterpart of the) directional alien derivation
$$
\De = \sum_{\om\in2\pi\I\,\N^*} \dDeom
$$
is the logarithm of~$\De^+$. {\em We thus must show that, for each
$\om\in2\pi\I\,\N^*$, $\dDeom$ is the $\om$-homogeneous component of the
operator}
$$
\log\De^+ = \sum_{r\ge1} \frac{(-1)^{r-1}}{r} \biggl(
\sum_{m\ge0}
\left( \dot A + \dot\ell_+ - \dot\ell_- \right) \dot\ell_+^{m} \mu
\biggr)^r.
$$
Using the obvious identity $\mu\left( \dot A + \dot\ell_+ - \dot\ell_- \right)=
\dot\ell_+ - \dot\ell_-$, we can write
$$
\log\De^+ = \sum_{\substack{m_1,\ldots,m_r\ge0 \\ r\ge1}} \frac{(-1)^{r-1}}{r} 
(\dot A + \dot\ell_+ - \dot\ell_-) \dot\ell_+^{m_1} \mu
(\dot A + \dot\ell_+ - \dot\ell_-) \dot\ell_+^{m_2} \mu \cdots
(\dot A + \dot\ell_+ - \dot\ell_-) \dot\ell_+^{m_r} \mu 
$$
as $\dst\sum_{m\ge1} (\dot A + \dot\ell_+ - \dot\ell_-) B_{m-1} \mu$,
with $2\pi\I(m-1)$-homogeneous operators
$$
B_{m-1} = \sum_{\substack{m_1+\cdots+m_r+r=m \\ m_1,\ldots,m_r\ge0, \ r\ge1}} 
\frac{(-1)^{r-1}}{r} \dot\ell_+^{m_1} (\dot\ell_+ - \dot\ell_-) 
\dot\ell_+^{m_2} \cdots (\dot\ell_+ - \dot\ell_-) \dot\ell_+^{m_r}.
$$
The result follows from the following identity (which is an identity for polynomials in two
non-commutative variables):
$$
B_{m-1} = \sum_{\eps_1,\ldots,\eps_{m-1}\in\{+,-\}}
\frac{p(\eps)!q(\eps)!}{m!}
\dot\ell_{\eps_{m-1}}\cdots\dot\ell_{\eps_1},
$$
the proof of which is left to the reader.

Formula~\refeq{eqLeibniz} of Proposition~\ref{propsimpleAlienD} follows, because
the logarithm of an automorphism is a derivation, and the homogeneous components
of a derivation are also derivations.

\sm

\begin{figure}

\begin{center}

\epsfig{file=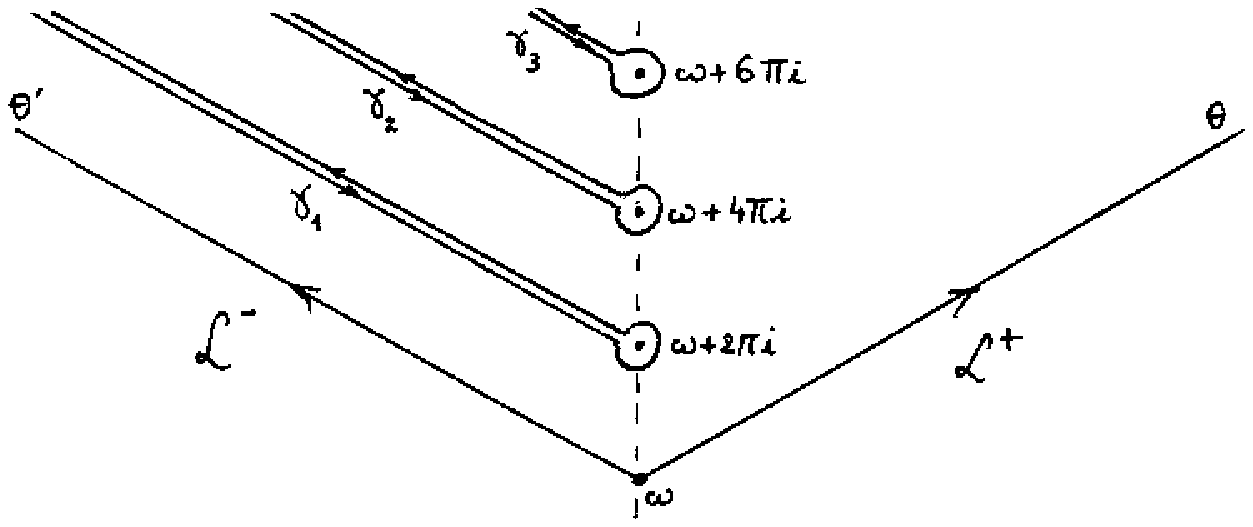,height=5.4cm,angle = 0}

\end{center}

\vspace{-.5cm}

\caption{\label{figStokes} Illustration of the formula $\cL^-\circ\De^+=\cL^+$.}

\end{figure}

As promised, we end this section with the interpretation of~$\De^+$ as Stokes
automorphism.
Let us extend $\cL^-=\cL^{\th'}$ for
$\th'\in\left]\frac{\pi}{2},\pi\right[$ to the part~$\cRexp$
of~$\hat\cR$ consisting of the formal sums $\sum(c_\om\,\de_\om+\hat\ph_\om)$
in which each~$\hat\ph_\om$ has at most exponential growth at infinity, by
setting
\begin{multline*}
\cL^-\sum_{\om\in2\pi\I\,\N} (c_\om\,\de_\om+\hat\ph_\om) = 
\sum_{\om\in2\pi\I\,\N} \bigl( c_\om\,\ee^{-\om z} +
\int_\om^{\ee^{\I\th'}\infty} \ee^{-z\ze}\hat\ph_\om(\ze)\,\dd\ze \bigr)\\[1ex]
= \sum_{\om\in2\pi\I\,\N}  \ee^{-\om z} \bigl( c_\om + 
\int_0^{\ee^{\I\th'}\infty} \ee^{-z\xi}\hat\ph_\om(\om+\xi)\,\dd\xi \bigr),
\end{multline*}
where the \rhs\ is to be considered as a formal series in~$\ee^{-2\pi\I\,z}$
with coefficients holomorphic in a domain~$\cD^-$. We define $\cL^+=\cL^{\th}$ for
$\th\in\left]0,\frac{\pi}{2}\right[$ similarly on~$\cRexp$.
We have
\begin{equation}	\label{eqStokesAut}
\cL^+\hat\chi_\om = \cL^- \De^+\hat\chi_\om,
\qquad \hat\chi_\om\in\cRexpom,
\end{equation}
at least when the series $\sum\ee^{-2\pi\I m z}\cL^-\De_{2\pi\I m z}^+\hat\chi_\om$
is convergent.
Indeed, one can decompose the contour of integration for~$\cL^+$ as follows:
$$
\cL^+\hat\chi_\om = c_\om\,\ee^{-\om z} + \left( \int_\om^{\ee^{\I\th'}\infty} +
\int_{\ga_1} + \int_{\ga_2} + \cdots \right) \ee^{-z\ze} \hat\ph_\om(\ze)\,\dd\ze
= \cL^-\hat\chi_\om + \sum_{m\ge1} \int_{\ga_m}\ee^{-z\ze} \hat\ph_\om(\ze)\,\dd\ze,
$$
where $\ga_m = \om + 2\pi\I\,m + \ga$ and~$\ga$ is the path coming
from~$\ee^{\I\th'}\infty$, turning anticlockwise around~$0$ and going back
to~$\ee^{\I\th'}\infty$ (see Figure~\ref{figStokes}), and the contribution
of~$\ga_m$ is precisely (because of the residuum formula)
\begin{multline*}
\al_{\om+2\pi\I m}(\ell_+^{m-1}\mu\hat\chi_\om) \,\, \ee^{-(\om+2\pi\I m) z} +
\int_{\om+2\pi\I m}^{\ee^{\I\th'}\infty} \ee^{-z\ze} \, (\dot\ell_+-\dot\ell_-)
\ell_+^{m-1}\mu\hat\chi_\om(\ze) \,\dd\ze = \\[1ex]
\cL^-\left(\dot A + \dot\ell_+-\dot\ell_-\right) \ell_+^{m-1}\mu\hat\chi_\om 
= \cL^-\dDep{2\pi\I m}\hat\chi_\om.
\end{multline*}

\sm

Formula~\refeq{eqStokesAut} may serve as a heuristic explanation of the fact
that $\De^+= \text{``$(\cL^-)\ii\circ\cL^+$''}$ is an automorphism, since
both~$\cL^+$ and~$\cL^-$ transform convolution into multiplication.
This formula is the expression of an ``abstract Stokes phenomenon'', without
reference to any particular equation, which manifests itself in the Stokes phenomenon when
specialized to the resurgent solution of an equation like~\refeq{eqnAbelInv}
or~\refeq{eqnAbel}. 
It was used in the form $\cL^+ = \cL^-\circ(\exp\De)$ in the previous section.


\section{Formalism of singularities, general resurgent functions and alien
derivations}	\label{secSingul} 


Let us return to the convolution algebra~$\wh\cH(\cR)$, consisting of all
holomorphic germs at the origin which extend to the Riemann surface~$\cR$.
In Section~\ref{secsimpl} we have focused on simple singularities and this led
to the definition of~$\Rsimp$, but what about more complicated singularities
than simple ones?

Even in the elementary situation described in Section~\ref{seclineardiffeq} with
the Borel transform~$\hat\psi(\ze)$ of the solution of the second linear
equation~\refeq{eqnpsi}, poles of order higher than~$1$ appear.

Or, as alluded to in Section~\ref{secAbEqtgtId} after
Proposition~\ref{propformalconj}, in the case of a nonzero resiter~$\rho$ Abel's
equation has a solution of the form
$\ID + \rho\log z + \ti\psi(z)$ with $\ti\psi\in z\ii\C[[z\ii]]$.
One can prove that~$\ti\psi\in\ti\cH$, \ie the Borel transform~$\hat\psi$ is
in~$\wh\cH(\cR)$, but the singularities one finds in the Borel plane can be of
the form $\ze^\al\hat\Phi(\ze)+\text{reg}(\ze)$ with $\al\in\C$ and
$\hat\Phi(\ze)$, $\text{reg}(\ze)$ regular at the origin (see~\cite[Vol.~2]{Eca81}).

It may also happen that an equation give rise to formal solutions involving
non-integer powers of~$z$, or monomials of the form $z^{-n}(\log z)^m$.

All these issues are addressed satisfactorily by the formalism of singularities,
as developed in~\cite[Vol.~3]{Eca81}, \cite{dulac} or~\cite{Eca93} (see
also~\cite{CNP} and~\cite{OSS}).

\subsection{General singularities. Majors and minors. Integrable singularities}

Let~$\Clog$ denote the Riemann surface of the logarithm, \ie
$$
\Clog = \widetilde{(\C\setminus\{0\},1)} =
\{ \ze = r\,\ee^{\I\th} \mid r>0, \ \th\in\R \}
$$
(\cf footnote~\ref{footClog} in Section~\ref{secAnContConv}).
We denote by $\ze \in \Clog \mapsto \dze \in \C\setminus\{0\}$ the
canonical projection (covering map).\footnote{
As a Riemann surface, $\Clog$ is isomorphic to~$\C$ (with a biholomorphism
$\ze\in\Clog \mapsto \log\ze\in\C$), but it is more significant for us to
consider it as a universal cover with a ``multivalued coordinate''~$\ze$.
}
We are interested in analytic functions which are potentially singular
at the origin of~$\C$, possibly with multivalued analytic continuation around
the origin.
We thus define~$\ANA$ to be the space of the germs of functions analytic in a
``spiralling neighbourhood of the origin'', \ie analytic in a domain
of the form
$\cV=\{ r\,\ee^{\I\th} \mid 0<r<h(\th), \ \th\in\R \} \subset \Clog$,
with a continuous function $h:\R \to\left]0,+\infty\right[$.
The space~$\C\{\ze\}$ of regular germs is obviously a subspace of~$\ANA$.

\begin{Def}
Let $\SING = \ANA/\C\{\ze\}$. The elements of this space are called
``singularities''.\footnote{	\label{footgensing}
A more general definition is given in~\cite{OSS}, with sectorial neighbourhoods
of the form $\cV = \{ r\,\ee^{\I\th} \mid \th_0-\al-2\pi < \th < \th_0+\al, \
0<r<h \}$ instead of spiralling neighbourhoods, giving rise to more general
germs $\chb\ph \in \ANA_{\th_0,\al}$ and singularities $\trb\ph \in
\SING_{\th_0,\al}=\ANA_{\th_0,\al}/\C\{\ze\}$. 
In practice, this is useful as an intermediary step to prove that the solution
of a nonlinear equation is resurgent, but here we simplify the exposition.
}
The canonical projection is denoted~$\sing_0$ and we use the notation
$$
\sing_0\,:\; \left\{ \begin{aligned}
\ANA &\to \SING \\[1ex]
\chb\ph \hspace{.65em} &\mapsto \: \trb\ph = \sing_0\bigl(\chb\ph(\ze)\bigr).
\end{aligned} \right.
$$
Any representative~$\chb\ph$ of a singularity~$\trb\ph$ is called a ``major'' of
this singularity.

The map induced by the variation map 
$\chb\ph(\ze) \mapsto \chb\ph(\ze) - \chb\ph(\ze\,\ee^{-2\pi\I})$
is denoted
$$
\var\,:\; \left\{ \begin{aligned}
\SING \hspace{.65em} &\to \hspace{.65em} \ANA \\[1ex]
\trb\ph=\sing_0\bigl(\chb\ph\bigr) &\mapsto 
\htb\ph(\ze) = \chb\ph(\ze) - \chb\ph(\ze\,\ee^{-2\pi\I}).
\end{aligned} \right.
$$
The germ $\htb\ph = \var\trb\ph$ is called the ``minor'' of the
singularity~$\trb\ph$. 
\end{Def}

Observe that the kernel of $\var:\SING\to\ANA$ is isomorphic to the space of
entire functions of~$\frac{1}{\ze}$ without constant term.
It turns out that the map~$\var$ is surjective (we omit the proof).

The simplest examples of singularities are poles
$$
\de = \sing_0\left( \frac{1}{2\pi\I\ze} \right), \qquad
\de\ss n = \sing_0\left( \frac{(-1)^n n!}{2\pi\I\ze^{n+1}} \right), \qquad
n\ge0
$$
(observe that $\var\de\ss n = 0$),
and logarithmic singularities with regular variation, for which we use the notation
\begin{equation}	\label{eqbemreg}
\bem\htb\ph = \sing_0\left( \frac{1}{2\pi\I} \htb\ph(\ze) \log\ze \right),
\qquad \htb\ph(\ze)\in\C\{\ze\}
\end{equation}
($\log\ze$ is well-defined since we work in~$\Clog$; anyway, another
branch of the logarithm would define the same singularity since the
difference of majors would be a multiple of~$\htb\ph(\ze)$ which is regular).
The last example is a particular case of ``integrable singularity''.

\begin{Def}	\label{definteg}
An ``integrable minor'' is a germ $\htb\ph\in\ANA$ which is uniformly integrable
at the origin in any sector $\th_1 \le \arg\ze \le \th_2$, in the sense that
for any $\th_1<\th_2$ there exists a Lebesgue integrable function
$f:\left]0,r^*\right]\to\left]0,+\infty\right[$ such that 
$$
|\htb\ph(\ze)| \le f(|\ze|), \qquad \ze\in S,
$$
where $S=\ao \ze\in\Clog \mid \th_1\le\arg\ze\le\th_2, \ |\ze|\le r^* \af$ and
$r^*>0$ is small enough so as to ensure that~$S$ be contained in the domain of
analyticity of~$\htb\ph$.
The corresponding subspace of~$\ANA$ is denoted by~$\ANA\integ$.

\sm

An ``integrable singularity'' is a singularity $\trb\ph\in\SING$ which admits a
major~$\chb\ph$ such that $\ze\chb\ph(\ze)\to0$ as $\ze\to0$ uniformly in any
sector $\th_1 \le \arg\ze \le \th_2$ and for which $\var\trb\ph\in\ANA\integ$.
The corresponding subspace of~$\SING$ is denoted by~$\SING\integ$.
\end{Def}

For example, the formulas
\begin{equation}	\label{eqdefIsig}
\trn I_\sig = \sing_0\bigl(\chn I_\sig\bigr), \qquad
\chn I_\sig(\ze) = \frac{\ze^{\sig-1}}{(1-\ee^{-2\pi\I\sig})\Ga(\sig)}, 
\qquad \sig\in\C\setminus\N^*
\end{equation}
define a family of singularities\footnote{
\label{footIsig}
In view of the poles of the Gamma function, $\trn I_\sig$ is well-defined for
$\sig=-n\in -\N$, and $\trn I_{-n} = \de\ss n$. 
Besides, the reflection formula yields
$$
\chn I_\sig(\ze) = \tfrac{1}{2\pi\I} \ee^{\pi\I\sig} \Ga(1-\sig) \ze^{\sig-1}.
$$
This family of singularities admits a non-trivial analytic continuation
\wrt~$\sig$ at positive integers (\cite[Vol.~1, pp.~47--51]{Eca81}): for $n\in\N^*$, one
may consider another major of~$\trn I_\sig$, which is analytic at $\sig=n$, and
define
$$
\chn I_n(\ze) = \lim_{\sig\to n} 
\frac{\ze^{\sig-1}-\ze^{n-1}}{(1-\ee^{-2\pi\I\sig})\Ga(\sig)}
= \frac{\ze^{n-1}\log\ze}{2\pi\I\,\Ga(n)},
$$
which yields $\trn I_n  = \bem\bigl( \frac{\ze^{n-1}}{(n-1)!} \bigr)$.
The formula for the minors
$$
\htn I_\sig(\ze) = \frac{\ze^{\sig-1}}{\Ga(\sig)}
$$
is thus valid for any $\sig\in\C\setminus(-\N)$
(while $\htn I_{-n}=0$ if $n\in\N$).
}
among which the integrable ones correspond to $\RE\sig>0$.
Another example is provided by polynomials of~$\log\ze$, which can be viewed as integrable
minors, and also as majors of integrable singularities. 

Observe that, when a singularity~$\trb\ph$ is integrable, any of its majors
satisfies the condition which is stated in Definition~\ref{definteg} (since the
difference between two majors, being regular, is $o(1/|\ze|)$), and that its
minor is by assumption an integrable minor.

\begin{figure}

\begin{center}

\epsfig{file=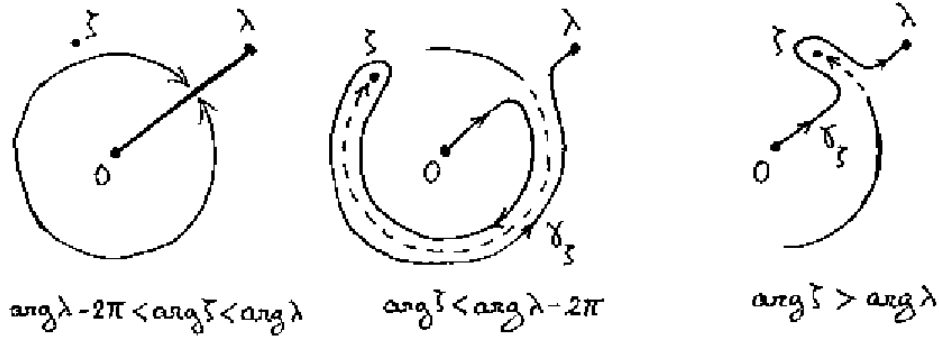,height=4.6cm,angle = 0}

\end{center}

\vspace{-.5cm}

\caption{\label{figMajMinInteg} Cauchy integral for a major associated with an
integrable minor.}

\end{figure}

\begin{lemma} 
By restriction, the variation map $\var:\SING\to\ANA$ induces a linear
isomorphism $\SING\integ \to
\ANA\integ$. 
\end{lemma}

\proof
In view of the definition of~$\SING\integ$, the variation map induces a linear
map $\var\integ :\, \SING\integ \to \ANA\integ$.
The injectivity of~$\var\integ$ is obvious: $\sing_0\bigl(\chb\ph\bigr)$ belongs
to the kernel of~$\var$ if and only if~$\chb\ph(\ze)$ is the sum of an entire
function of~$1/\ze$ and of a regular germ, and the condition
$\chb\ph(\ze)=o(1/|\ze|)$ leaves room for the regular germ only.

\sm

For the surjectivity, we suppose $\htb\ph\in\ANA\integ$ and we only need to
exhibit a germ $\chb\ph\in\ANA$ with variation~$\htb\ph$ and with the property
of being~$o(1/|\ze|)$ uniformly near the origin.
If~$\htb\ph$ is regular, we can content ourselves with setting 
$\chb\ph(\ze) = \frac{1}{2\pi\I} \htb\ph(\ze) \log\ze$, 
but for the general case we resort to a Cauchy integral.

Let us fix an auxiliary point~$\la$ in the domain of analyticity $S\subset\Clog$
of~$\htb\ph$.
The integrability of~$\htb\ph$ allows us to define a holomorphic function by the
formula
\begin{equation}	\label{eqmajintegr}
\chb\ph(\ze) = -\frac{1}{2\pi\I} \int_0^\la \frac{\htb\ph(\ze_1)}{\ze_1-\ze}\,\dd\ze_1,
\qquad \arg\la-2\pi < \arg\ze < \arg\la.
\end{equation}
This function admits an analytic continuation to $S\setminus
\la\left[1,+\infty\right[$, as is easily seen by deforming the path of
integration (see Figure~\ref{figMajMinInteg}): 
$$
\chb\ph(\ze) = -\frac{1}{2\pi\I} \int_{\ga_\ze}
\frac{\htb\ph(\ze_1)}{\ze_1-\ze} \,\dd\ze_1
$$
with a path~$\ga_\ze$ connecting~$0$ and~$\la$ inside~$S$ and
circumventing~$\ze$ to the right if $\arg\ze<\arg\la-2\pi$ or to the left if
$\arg\ze>\arg\la$ (turning around the origin as many times as necessary to reach
the sheet of~$\Clog$ where~$\ze$ lies before turning back to sheet of~$\la$).
For~$\arg\ze$ slightly larger than~$\arg\la$ and~$|\ze|$ small enough, the
residuum formula yields 
$$
\chb\ph(\ze) - \chb\ph(\ze\,\ee^{-2\pi\I}) = 
- \frac{1}{2\pi\I} \left(\int_{\ga_\ze} - \int_0^\la \right)
\frac{\htb\ph(\ze_1)}{\ze_1-\ze} \,\dd\ze_1
= \htb\ph(\ze)
$$
(because this difference of paths is just a circle around~$\ze$ with clockwise
orientation). 

\sm

For the uniform $o(1/|\ze|)$ property, it is sufficient to check that 
$\ze\chb\ph(\ze)\to0$ as $\ze\to0$ uniformly in
$\Sig_{\la,\al} = \ao \ze\in\Clog \mid \arg\la - 2\pi + \al \le \arg\ze \le
\arg\la-\al \af$, with arbitrarily small $\al>0$,
since changing~$\la$ only amounts to adding a regular germ to~$\chb\ph$.
We write 
$\ze\chb\ph(\ze) = -\frac{\la}{2\pi\I} \int_0^1 \htb\ph(t\la) 
\frac{\ze}{t\la-\ze}\dd t$
and observe that, for $\ze\in\Sig_{\la,\al}$,
$\cos\bigl(\arg\frac{\la}{\ze}\bigr)\le\cos\al$ hence
\begin{multline*}
\tst \bigl|t\frac{\la}{\ze}-1\bigr|^2 =
t^2 \bigl|\frac{\la}{\ze}\bigr|^2 + 1 
- 2t \bigl|\frac{\la}{\ze}\bigr| \cos\bigl(\arg\frac{\la}{\ze}\bigr) 
\ge F(t,|\ze|), \\[1ex]
\tst F(t,r) = t^2 \frac{|\la|^2}{r^2} + 1 
- 2t \frac{|\la|}{r} \cos\al
\ge \sin^2\al.
\end{multline*}
We can conclude by Lebesgue's dominated convergence theorem:
$$
|\ze\chb\ph(\ze)| \le \frac{|\la|}{2\pi} \int_0^1 g(t,|\ze|)\,\dd t, \qquad
g(t,r) = \frac{1}{\sqrt{F(t,r)}} |\htb\ph(t\la)| 
= \frac{r}{\bigl| t|\la|\ee^{\I\al}-r \bigr|} |\htb\ph(t\la)|,
$$
with $g(t,r)\le\frac{1}{\sin\al}|\htb\ph(t\la)|$ integrable
and $g(t,r)\to0$ as $r\to0$ for each~$t>0$,
hence $|\ze\chb\ph(\ze)| \le \eps_{\la,\al}(|\ze|)$ with 
$\eps_{\la,\al}(r)\xrightarrow[r\to0]{}0$.
\eopf

\noindent{\bf Notation:}
The inverse map will be denoted 
$$
\htb\ph\in\ANA\integ \mapsto \bem\htb\ph \in \SING\integ,
$$ 
a notation which is consistent with~\refeq{eqbemreg}.
In the spirit of Definition~\ref{defsimplesing}, we can define the ``simple
singularities'' (at the origin) as those singularities of the form
$c\,\de + \bem\htb\ph$, with $c\in\C$ and $\htb\ph(\ze) \in\C\{\ze\}$; 
we denote by $\SING\simp = \C\,\de \oplus \bem\bigl(\C\{\ze\}\bigr)$
the subspace of~$\SING$ that they form.

\subsection{The convolution algebra $\SING$}

Starting with Section~\ref{secAnContConv}, we have dealt with convolution of
regular germs. But the space~$\C\{\ze\}$ of regular germs is contained in the
space~$\ANA\integ$ of integrable minors, and we can extend the convolution law:
$$
\htb\ph_1,\htb\ph_2\in\ANA\integ \;\mapsto 
\htb\ph_3=\htb\ph_1*\htb\ph_2\in\ANA\integ,
\qquad \htb\ph_3(\ze) = \int_0^\ze \htb\ph_1(\ze_1)\htb\ph_2(\ze-\ze_1).
%
$$
Indeed, in any sector $\th_1 \le \arg\ze \le \th_2$, using integrable functions
$f_1,f_2:\,\left]0,r^*\right] \to\left]0,+\infty\right[$ such that
$|\htb\ph_i(\ze)|\le f_i(|\ze|)$, we see that the formula makes sense for
$|\ze|\le r^*$ and defines a holomorphic function such that 
$|\htb\ph_3(\ze)|\le f_3(|\ze|)$ with $f_3(r)=\int_0^1 f_1(tr)f_2((1-t)r)r\,\dd t$;
the positive function $f_3=f_1*f_2$ itself is integrable, by virtue of the Fubini
theorem:
$\int_0^{r^*} f_3(r)\,\dd r = \iint_0^{r^*} f_1(r_1)f_2(r_2)\,\dd r_1\,\dd r_2 <
\infty$.

\sm

We have for instance $\htn I_{\sig_1} * \htn I_{\sig_2} = \htn I_{\sig_1+\sig_2}$
for all complex $\sig_1$, $\sig_2$ with positive real part, whether integer or
not (by the classical formula for the Beta function).
The extended convolution is called ``convolution of integrable minors''; it is
still commutative and associative.
We thus get an algebra~$\ANA\integ$ (without unit), with~$\C\{\ze\}$ as a subalgebra.

\sm

Transporting this structure of algebra by~$\var\integ$, we can
view~$\SING\integ$ as an algebra, with convolution law
\begin{equation}	\label{eqdefconvintegsing}
\trb\ph_1=\bem\htb\ph_1,\, \trb\ph_2=\bem\htb\ph_2 \in \SING\integ \;\mapsto 
\trb\ph_1*\trb\ph_2 := \bem\bigl(\htb\ph_1*\htb\ph_2\bigr) \in\SING\integ.
\end{equation}
It turns out that the convolution law for integrable singularities can be
extended to the whole space of singularities, so as to make~$\SING$ an algebra,
of which $\SING\integ=\bem\bigl(\ANA\integ\bigr)$ will appear as a subalgebra
(and there will be a unit, namely~$\de$).

\subsub{Convolution with integrable singularities}

As an introduction to the definition of the convolution of general
singularities, let us begin with a more careful study of
$\trb\ph_3=\trb\ph_1*\trb\ph_2$ in the integrable
case~\refeq{eqdefconvintegsing}: we shall indicate formulas for the minor and a
major of~$\trb\ph_3$, which do not make reference to the minor~$\htb\ph_1$ but
only to a major of~$\bem\htb\ph_1$.

\begin{figure}

\begin{center}

\epsfig{file=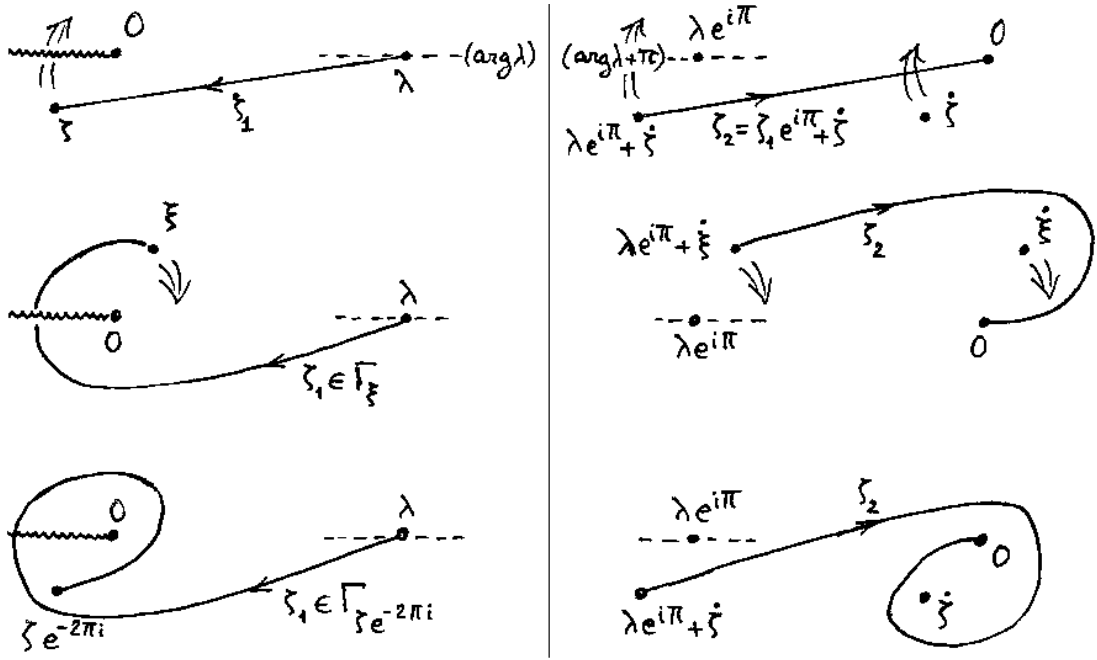,height=8.7cm,angle = 0}

\end{center}


\caption{\label{figMajMinMaj} Analytic continuation of~$\protect\chb\ph_{3,\la}$.
(Left: three examples of path~$\Ga$. Right: the corresponding paths followed by
$\ze_2=\ze_1\,\ee^{\pi\I}+{\protect\dze}$).}

\end{figure}

\begin{lemma}	\label{lemmajminM}
Let $\htb\ph_1,\htb\ph_2\in\ANA\integ$, $\htb\ph_3=\htb\ph_1*\htb\ph_2$, and
let~$\chb\ph_1$ be any major of~$\bem\htb\ph_1$.
Then, for~$\ze$ with small enough modulus, one has
\begin{equation}	\label{eqmajminmin}
\htb\ph_3(\ze) = \int_{\ga_\ze} \chb\ph_1(\ze_1) \htb\ph_2(\ze-\dze_1) \,\dd\ze_1,
\end{equation}	
where~$\ga_\ze$ is any path in~$\Clog$ which starts at~$\ze\,\ee^{-2\pi\I}$,
turns around the origin anticlockwise and ends at~$\ze$, \eg the circular path 
$t\in[0,1]\;\mapsto\; \ze\,\ee^{-2\pi\I(1-t)}$.
In the above formula, $\ze_2=\ze-\dze_1$ denotes the lift of~$\dze-\dze_1$ which
lies in the same sheet of~$\Clog$ as~$\ze$; this point thus starts from and
comes back to the origin after turning anticlockwise around~$\ze$ (rather
than~$\ze\,\ee^{-2\pi\I}$, or any other lift of~$\dze$ in~$\Clog$).

Let $\la\in\Clog$ belong to the domain of analyticity of~$\chb\ph_1$, with small
enough modulus so that $\la\,\ee^{\I\pi}$ belongs to the domain of analyticity
of~$\htb\ph_2$.
Then the formula
\begin{equation}	\label{eqmajminmaj}
\chb\ph_{3,\la}(\ze) = \int_\la^\ze \chb\ph_1(\ze_1) 
\htb\ph_2(\ze_1\,\ee^{\I\pi}+\dze) \,\dd\ze_1,
\qquad \arg\la-\pi < \arg\ze < \arg\la+\pi, \quad \text{$|\ze|$ small enough}
\end{equation}	
defines by analytic continuation an element of~$\ANA$ which is a major
of~$\bem\htb\ph_3$.
In formula~\refeq{eqmajminmaj}, it is understood that $\ze_2 =
\ze_1\,\ee^{\I\pi}+\dze$ moves along the segment $[\la\,\ee^{\I\pi}+\dze,0]$,
where $\la\,\ee^{\I\pi}+\dze$ denotes the lift in~$\Clog$ of~$-\dla+\dze$ which
has its argument closest to~$\arg\la+\pi$ (it is well-defined for
$|\ze|<|\la|$)---see the top of Figure~\ref{figMajMinMaj}.
\end{lemma}

\proof 
a) Let $\ze\in\Clog$ with small enough modulus.
Observe that, in formula~\refeq{eqmajminmin}, $\arg\ze_1$ takes all possible
values beteween $\arg\ze-2\pi$ and~$\arg\ze$ at least, for whatever choice
of~$\ga_\ze$, while $\arg(\ze-\dze_1)$ (with the convention indicated) can be
maintained arbitrarily close to~$\arg\ze$ by choosing~$\ga_\ze$ close enough to
the segments $[\ze\,\ee^{-2\pi\I},0]$ and~$[0,\ze]$.
Let $\eps$ denote a positive function on~$\left]0,|\ze|\right]$ such that
$\eps(r)\xrightarrow[r\to0]{}0$ and
$$
|\ze_1\chb\ph_1(\ze_1)| \le \eps(|\ze_1|), \qquad
\arg\ze-2\pi \le \arg\ze_1 \le \arg\ze.
$$
Deforming the contour, we rewrite the \rhs\ of~\refeq{eqmajminmin} as
$$
\left( \int_{\ze\,\ee^{-2\pi\I}}^{a\,\ee^{-2\pi\I}} + \int_{\ga_a} + \int_a^\ze \right)
\chb\ph_1(\ze_1) \htb\ph_2(\ze-\dze_1) \,\dd\ze_1,
$$
with any auxiliary point $a\in\left]0,\demi\ze\right]$.
The two integrals over rectilinear segments contribute $\int_a^\ze
\htb\ph_1(\ze_1) \htb\ph_2(\ze-\dze_1) \,\dd\ze_1$, which tends
to~$\htb\ph_3(\ze)$ as $a\to0$ by integrability of~$\htb\ph_1$ and~$\htb\ph_2$.
The integral over~$\ga_a$ tends to~$0$, since its modulus is not larger than
$2 \pi \, \eps(|a|) \max_D |\htb\ph_2|$, where~$D$ denotes the closed disc of
radius~$\demi|\ze|$ centred at~$\ze$ (inside the same sheet of~$\Clog$ as~$\ze$).

\sm

\noindent
b) Formula~\refeq{eqmajminmaj} defines a function holomorphic near the
origin of~$\Clog$ in the sheet (branch cut) which corresponds to the condition 
$0\notin[\dla,\dze]$ and which contains~$\la$,
\ie in $S_\la = \{ \arg\la-\pi < \arg\ze < \arg\la+\pi \}$ for small enough~$|\ze|$;
for such~$\ze$, the whole path of integration~$[\la,\ze]$ lies inside~$S_\la$,
while the corresponding $\ze_2 = \ze_1\,\ee^{\I\pi} + \dze$ move along the
segment $[\la\,\ee^{\I\pi} + \dze,0] \subset S_{\la\,\ee^{\I\pi}}$
(here we use again the integrability of~$\htb\ph_2$).
The analytic continuation is obtained by deforming continuously the contour of
integration as~$\ze$ moves to different sheets of~$\Clog$:
$$
\chb\ph_{3,\la}(\ze) = \int_{\Ga_\ze} \chb\ph_1(\ze_1) 
\htb\ph_2(\ze_1\,\ee^{\I\pi}+\dze) \,\dd\ze_1,
$$
with a path~$\Ga_\ze$ which connects~$\la$ and~$\ze$ without touching the origin
and without intersecting itself, it being understood that the symmetric path
followed by $\ze_2 = \ze_1\,\ee^{\I\pi}+\dze$ always starts on the same sheet
of~$\Clog$---see Figure~\ref{figMajMinMaj}.

In particular, following the analytic continuation from~$\ze$
to~$\ze\,\ee^{-2\pi\I}$, we obtain a path~$\Ga_{\ze\,\ee^{-2\pi\I}}$ which can be
decomposed as the original path~$[\la,\ze]$ followed by~$\ti\ga_\ze$,
where~$\ti\ga_\ze$ is the same as the above~$\ga_\ze$ but with reverse
orientation. 
Moreover, our convention for~$\arg\ze_2$ agrees\footnote{
See Figure~\ref{figMajMinMaj} or draw your own picture in the simpler case
when $\arg\ze$ is slightly smaller than $\arg\la+\pi$.
Or start directly with the limit case $\arg\ze=\arg\la+\pi$, with $\Ga_\ze$
and~$\Ga_{\ze\,\ee^{-2\pi\I}}$ starting at~$\la$ and following $[\dla,\dze]$ but
circumventing~$0$ to the right or to the left, while the symmetric path starts at
$\la\,\ee^{\I\pi}+\dze$ and circumvents~$\dze$ to the right or to the left; the
point is that the symmetric path thus lies in the same sheet as~$\ze$ (rather
than~$\ze\,\ee^{-2\pi\I}$, or any other lift of~$\dze$ in~$\Clog$), this is the
origin of our convention on $\arg\ze_2$.
}
with that of formula~\refeq{eqmajminmin}, hence
$\chb\ph_{3,\la}(\ze) - \chb\ph_{3,\la}(\ze\,\ee^{-2\pi\I}) = \htb\ph_3(\ze)$.

\sm

\noindent
c) Let
$$
\Sig'_{\la,\al} = \lao |\arg\ze-\arg\la| \le \pi-\al\laf, \qquad
S_{\la,\al} = \Sig'_{\la,\al} \cap \lao |\ze|\le\dem|\la| \laf,
\qquad \al\in\left]0,\tfrac{\pi}{2}\right[.
$$
We observe that, if $0<\arg\la'-\arg\la<2\pi$,
$(\chb\ph_{3,\la}-\chb\ph_{3,\la'})(\ze)$ 
(which, for $\arg\la'-\pi \le \arg\ze < \arg\la+\pi$, is given by an integral
involving the values of~$\htb\ph_1$ close to $[\la,\la']$ and those
of~$\htb\ph_2$ close to $[\la\,\ee^{\I\pi}+\dze,\la'\,\ee^{\I\pi}+\dze]$)
extends to a holomorphic function which is regular near $\ze=0$.
Thus, to prove that $\bem\htb\ph_3=\sing_0\bigl(\chb\ph_{3,\la}\bigr)$, it only
remains to check that $\ze\chb\ph_{3,\la}(\ze)\to0$ as $\ze\to0$, uniformly in
any sector~$\Sig'_{\la,\al}$.

Let $r\mapsto\eps(r)$ denote a positive function defined
on~$\left]0,|\la|\right]$ such that $\eps(r)\xrightarrow[r\to0]{}0$ and
$|\ze_1\htb\ph_1(\ze_1)|\le\eps(|\ze_1|)$ whenever $\ze_1\in\Sig'_{\la,\al}$ and
$|\ze_1|\le|\la|$.
In fact, we shall only use the fact that one can take a {\em bounded}
function~$\eps$: we only suppose $\eps(r)\le \eps^*$
(this is related to footnote~\ref{footdeplusintegr}).
Let $r\mapsto f(r)$ denote a positive integrable function defined
on~$\bigl]0,\frac{3|\la|}{2}\bigr]$ such that $|\htb\ph_2(\ze_2)| \le f(|\ze_2|)$ 
whenever $|\ze_2|\le\frac{3}{2}|\la|$ and
$\arg\la+\frac{\pi}{2} \le \arg\ze_2 \le \arg\la+\frac{3\pi}{2}$.
We can write
$$
\chb\ph_{3,\la}(\ze) = \int_0^{\la\,\ee^{\I\pi}+\sdze}
\chb\ph_1(\ze-\dze_2)\htb\ph_2(\ze_2)\,\dd\ze_2, 
\qquad \ze \in S_{\la,\al}.
$$
For any $\ze\in S_{\la,\al}$, letting $u_\ze =
\frac{\la\,\ee^{\I\pi}+\sdze}{|\la-\ze|}$, we thus have
$$
|\ze\chb\ph_{3,\la}(\ze)| \le \int_0^{|\la-\ze|} 
\left|\tfrac{\ze}{\ze-\xi u_\ze}\right| \eps^* f(\xi) \,\dd\xi.
$$
Elementary geometry shows that
$$
\ze\in S_{\la,\al} \;\Rightarrow\; 
\al' \le \arg\bigl( \tfrac{u_\ze}{\ze} \bigr) 
= \arg\bigl( \tfrac{\la\,\ee^{\I\pi}}{\ze} + 1 \bigr) \le 2\pi - \al',
\qquad \al' = \arg\left(2\,\ee^{\I\al} + 1 \right) \in \left]0,\al\right[.
$$
Hence
\begin{equation*}
\tst \bigl|\frac{\ze-\xi u_\ze}{\ze}\bigr|^2 =
1 + \bigl|\frac{\xi}{\ze}\bigr|^2 
- 2 \bigl|\frac{\xi}{\ze}\bigr| \cos\bigl(\arg\frac{u_\ze}{\ze}\bigr) 
\ge F(\xi,|\ze|), \qquad
\tst F(\xi,r) = 1 + \frac{\xi^2}{r^2} 
- 2 \frac{|\xi|}{r} \cos(\al').
\end{equation*}
But $g=1/\sqrt{F}$ is continuous and $\le 1/\sin(\al')$ on
$[0,\frac{3|\la|}{2}]\times\left]0,|\la|\right]$, 
with $g(\xi,r)\xrightarrow[r\to0]{}0$ for each $\xi>0$, thus
$$
|\ze\chb\ph_{3,\la}(\ze)| \le \eps^* \int_0^{\frac{3|\la|}{2}} 
g(\xi,|\ze|) f(\xi) \,\dd\xi =\eps'_{\la,\al}(|\ze|),
$$
with $\eps'_{\la,\al}(r)\xrightarrow[r\to0]{}0$ by Lebesgue's dominated
convergence theorem.
\eopf


We now see how we can define the convolution product of a general
singularity~$\trb\ph_1$ with an integrable one~$\bem\htb\ph_2$:
for any major $\chb\ph_1\in\ANA$ and any~$\la$, formula~\refeq{eqmajminmaj}
still defines an element~$\chb\ph_{3,\la}$ of~$\ANA$ by analytic continuation
(the integrability of~$\trb\ph_1$ is not required for this), and we can set
\begin{equation}	\label{eqdefconvolgeninteg}
\trb\ph_1*\bem\htb\ph_2 = \sing_0(\chb\ph_{3,\la}).
\end{equation}
The choice of the major~$\chb\ph_1$ does not matter (by linearity, adding
to~$\chb\ph_1$ a regular germ in~\refeq{eqmajminmaj} will add
to~$\chb\ph_{3,\la}$ a function which is regular, being a major of the null
singularity, in view of the above), nor does the choice of~$\la$ because
$\chb\ph_{3,\la}-\chb\ph_{3,\la'}$ is regular (as was mentioned in the above
proof).
The minor of~$\trb\ph_1*\bem\htb\ph_2$ is still given by
formula~\refeq{eqmajminmin}, but one must realize that $\trb\ph_1*\bem\htb\ph_2$
has no reason to be an integrable singularity when~$\trb\ph_1$ is not.\footnote{
\label{footdeplusintegr}
Notice however that, with this definition, $\de*\bem\htb\ph_2=\bem\htb\ph_2$
(compare $\chb\ph_{3,\la}(\ze)$ when $\chb\ph_1(\ze_1)=1/2\pi\I\ze_1$ and the
major $\chb\ph_2(\ze) = \int_0^{\la\ee^{\I\pi}}
\frac{\shtb\ph_2(\ze_2)}{2\pi\I(\ze-\ze_2)}\,\dd\ze_2$ 
given by~\refeq{eqmajintegr}: the difference is regular at the origin).
Hence, if $\trb\ph_1 = c\,\de + \bem\htb\ph_1$, $\trb\ph_1*\bem\htb\ph_2$ is still
an integrable singularity (namely $\bem(c\htb\ph_2+\htb\ph_1*\htb\ph_2)$).
On the other hand, the minor of $\de\ss{n+1}*\bem\htb\ph_2$ is
$\bigl(\frac{\dd}{\dd\ze}\bigr)^{n+1}\htb\ph_2$, but this singularity is usually
not integrable; for instance, if $\htb\ph_2\in\C\{\ze\}$, one can check that
$\de\ss{n+1}*\bem\htb\ph_2 = \htb\ph_2(0) \, \de\ss{n} + 
\frac{\dd\shtb\ph_2}{\dd\ze}(0) \, \de\ss{n-1} + \cdots + 
\frac{\dd^{n}\shtb\ph_2}{\dd\ze^{n}}(0) \, \de + 
\bem\bigl( \bigl( \frac{\dd}{\dd\ze}\bigr)^{n+1}\htb\ph_2 \bigr)$.
}

\subsub{Convolution of general singularities. The convolution algebra~$\SING$}

We just saw how the convolution of integrable minors gave rise to a convolution
of integrable singularities
$\SING\integ\times\SING\integ \to \SING\integ$
which could be extended to a convolution
$\SING\times\SING\integ \to \SING$.
We proceed with a further extension so as to view the space~$\SING$ as an
algebra, of which~$\SING\integ$ will appear as a subalgebra.

\sm

To this end, it is sufficient to imitate the arguments leading to
Lemma~\ref{lemmajminM} and to express
a major~$\chbph{*}_{3,\la}$ of the convolution product of two integrable
singularities by a formula similar 
to~\refeq{eqmajminmaj}, but referring to a major~$\chb\ph_2$ rather than to the
minor~$\htb\ph_2$ of the second singularity.
The new formula will then be taken as a definition of the convolution product
when the singularities are no longer assumed to be integrable.

\sm

\begin{figure}

\begin{center}

\epsfig{file=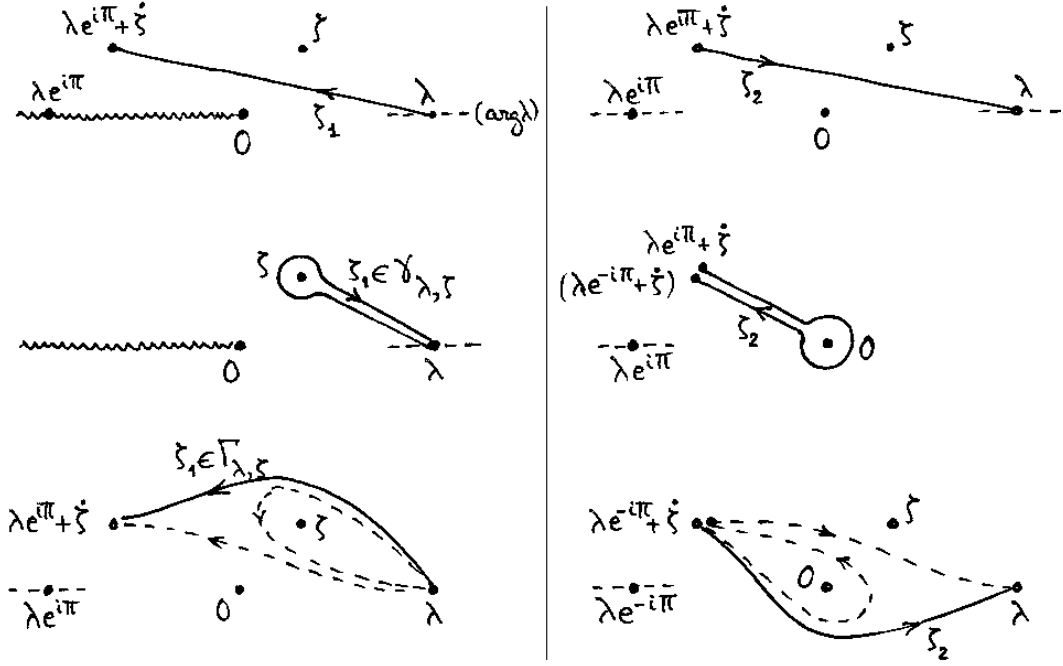,height=9cm,angle = 0}

\end{center}


\caption{\label{figConvolMaj} Top: integration path
for $\protect\chbph{*}_{3,\la}(\ze)$.
Middle: path for $\protect\chbph{**}_{3,\la}(\ze)$.
Bottom: difference.
(Left: paths followed by~$\ze_1$. Right: corresponding paths followed by~$\ze_2$.)}

\end{figure}

\begin{lemma}	\label{lemconvolmaj}
Let $\trb\ph_1,\trb\ph_2\in\SING$ have majors~$\chb\ph_1$ and~$\chb\ph_2$.
Let $\la\in\Clog$ belong to the intersection of the domains of analyticity
of~$\chb\ph_1$ and ~$\chb\ph_2$, with small enough modulus so that
$\la\,\ee^{\I\pi}$ also belongs to this intersection.
Then an element of~$\ANA$ can be defined by analytic continuation from the formula
\begin{equation}	\label{eqconvolmaj}
\chbph{*}_{3,\la}(\ze) = \int_\la^{\la\,\ee^{\I\pi}+\sdze} \chb\ph_1(\ze_1) 
\chb\ph_2(\ze_2) \,\dd\ze_1,
\qquad \ze\in H_\la, \quad \text{$|\ze|$ small enough},
\end{equation}	
where $H_\la = \ao \ze\in\Clog \mid \arg\la < \arg\ze < \arg\la+\pi \af$, and
where it is understood that $\la\,\ee^{\I\pi}+\dze$ is the lift in~$\Clog$
of~$-\dla+\dze$ which lies in~$H_\la$ and~$\ze_2$ is the lift of $\dze-\dze_1$
which is also in~$H_\la$ (and thus moves backwards along the same segment
$[\la,\la\,\ee^{\I\pi}+\dze]$)---see the top of Figure~\ref{figConvolMaj}.
This germ gives rise to a singularity
\begin{equation}	\label{eqdefconvolsing}
\trb\ph_3 := \sing_0\bigl( \chbph{*}_{3,\la} \bigr)
\end{equation}
which does not depend on~$\la$, nor on the choice of the majors~$\chb\ph_1$
and~$\chb\ph_2$, but only on the singularities~$\trb\ph_1$ and~$\trb\ph_2$.

Moreover, when $\trb\ph_2 \in \SING\integ$, $\trb\ph_3$ coincides with the
singularity~$\trb\ph_1*\trb\ph_2$ defined by formula~\refeq{eqdefconvolgeninteg}
in the previous section.
(In particular, when both~$\trb\ph_1$ and~$\trb\ph_2$ are integrable
singularities, we recover $\bem\htb\ph_1 * \bem\htb\ph_2 = \trb\ph_3$.)
\end{lemma}

\proof
Formula~\refeq{eqconvolmaj} defines an analytic function, since the segment
$[\la,\la\,\ee^{\I\pi}+\dze]$ is contained in the domains of analyticity
of~$\chb\ph_1$ and~$\chb\ph_2$.
This is not the case when $\arg\ze$ crosses $\arg\la$ or $\arg\la+\pi$, but the
analytic continuation is then obtained by deforming the path so that the origin
be avoided by~$\ze_1$ and~$\ze_2$.
Another way of obtaining the analytic continuation is to observe that, whenever
$|\arg\la'-\arg\la| < \pi$, the difference $(\chbph{*}_{3,\la'}-\chbph{*}_{3,\la})(\ze)$ 
(which, for $\ze\in H_{\la}\cap H_{\la'}$, is given by an integral involving the
values of~$\chb\ph_1$ and~$\chb\ph_2$ close to $[\la,\la']$ and
$[\la\,\ee^{\I\pi},\la'\,\ee^{\I\pi}]$)
extends analytically to a full neighbourhood of the origin.
Thus $\chbph{*}_{3,\la}\in\ANA$ and $\trb\ph_3$ does not depend on~$\la$.

\sm

The fact that~$\trb\ph_3$ does not depend on the choice of the major~$\chb\ph_2$
follows by linearity from the last statement (a regular~$\chb\ph_2$ can be
viewed as the major of an integrable singularity, namely the null singularity, the
convolution product of which with any~$\trb\ph_1$ is the null singularity).
The fact that it does not depend on the choice of the major~$\chb\ph_1$ then follows
from the commutativity of formula~\refeq{eqconvolmaj}.

\sm

Consider
\begin{equation}	\label{eqchbphsst}
\chbph{**}_{3,\la}(\ze) = 
\int_{\ga_{\la,\ze}} \chb\ph_1(\ze_1) \chb\ph_2(\ze_2) \,\dd\ze_1,
\qquad \arg\la-\pi < \arg\ze < \arg\la+\pi, \quad \text{$|\ze|$ small enough},
\end{equation}
where the path~$\ga_{\la,\ze}$ starts at~$\la$, goes towards~$\ze$, follows a
circle of small radius around~$\ze$ with clockwise orientation, and comes back
to~$\la$, and where $\ze_2$ is the lift of~$\dze-\dze_1$ which starts at
$\la\,\ee^{\I\pi}+\dze$ and ends at $\la\,\ee^{-\I\pi}+\dze$ after having turned
clockwise around the origin (see the middle of Figure~\ref{figConvolMaj}).
We observe that, for $\ze\in H_\la$ with small enough modulus,
$$
- \chbph{*}_{3,\la}(\ze) + \chbph{**}_{3,\la}(\ze) = 
\int_{\Ga_{\la,\ze}} \chb\ph_1(\ze_1) \chb\ph_2(\ze_2) \,\dd\ze_1,
$$
where the path~$\Ga_{\la,\ze}$ starts at~$\la$, circumvents both~$0$ and~$\ze$
to the right and ends at $\la\,\ee^{\I\pi}+\dze$, while the lift~$\ze_2$ of
$\dze-\dze_1$ starts at $\la\,\ee^{-\I\pi}+\dze$ and ends at~$\la$ after having
circumvented both~$0$ and~$\ze$ to the right (see the bottom of
Figure~\ref{figConvolMaj}),
and the function defined by the last integral is regular at the origin (because
$\Ga_{\la,\ze}$ can keep off the origin even if~$\ze$ varies in a full
neighbourhood of the origin).
Hence
$$
\sing_0\bigl( \chbph{**}_{3,\la} \bigr) 
= \sing_0\bigl( \chbph{*}_{3,\la} \bigr)
= \trb\ph_3.
$$
Suppose now that $\trb\ph_2\in\SING\integ$ and define $\trb\ph_1*\trb\ph_2$ by
formula~\refeq{eqdefconvolgeninteg}. When letting an auxiliary point~$a$ tend
to~$\ze$ along $[\la,\ze]$ and using a path
$\ga_{\la,\ze} = [\la,a] \cup \ga_{a,\ze} \cup [a,\la]$,
the alternative major~$\chbph{**}_{3,\la}$ of~$\trb\ph_3$ appears to be nothing
but the major~$\chb\ph_{3,\la}$ of $\trb\ph_1*\trb\ph_2$ delivered
by~\refeq{eqmajminmaj}.
This ends the proof.
\eopf


\begin{Def}
We define the convolution product
$$
\trb\ph_3 = \trb\ph_1 * \trb\ph_2
$$
of any two singularities $\trb\ph_1,\trb\ph_2\in\SING$ by
formulas~\refeq{eqconvolmaj} and~\refeq{eqdefconvolsing}.
\end{Def}

\begin{prop}
The convolution law just defined on the space~$\SING$ is commutative and associative; it
turns it into a commutative algebra, with unit $\de = \sing_0\Bigl(
\frac{1}{2\pi\I\,\ze} \Bigr)$.
\end{prop}

\proof
The commutativity is obvious. 
The relation $\trb\ph_1 * \de = \trb\ph_1$ is immediate when using the
alternative major~$\chbph{**}_{3,\la}$ of formula~\refeq{eqchbphsst} with
$\chb\ph_2(\ze_2) = 1/2\pi\I\ze_2$, since the residuum formula gives
$\chbph{**}_{3,\la}=\chb\ph_1$.
\sm

For the associativity, the quickest proof consists in extending it
from~$\SING\integ$ (in restriction to which it is a mere consequence of the
associativity of the convolution of integrable minors) to~$\SING$ by {\em
continuity} and {\em density}.
Indeed, we may call a sequence $\bigl(\trb\ph_n\bigr)_{n\ge0}$ of~$\SING$
convergent if these singularities admit majors~$\chb\ph_n$ analytic in the same
spiralling neighbourhood of the origin~$\cV$ and if there exists~$\chb\ph$
analytic in~$\cV$ such that $\bigl(\chb\ph_n\bigr)_{n\ge0}$ converges uniformly
towards~$\chb\ph$ in every compact subset of~$\cV$; the singularity
$\sing_0\bigl(\chb\ph\bigr)$ is then unique and is called the limit of the
sequence.
It is easy to check that 
$$
\trb\ph_n\to\trb\ph \ens\text{and}\ens \tr\psi_n\to\tr\psi \quad\Rightarrow\quad
\trb\ph_n * \tr\psi_n \to \trb\ph*\tr\psi.
$$
On the other hand, any singularity is the limit of a sequence of integrable
singularities, majors of which can be chosen to be polynomials in~$\log\ze$
(this essentially amounts to the Weierstrass theorem in the
variable~$\log\ze$). 
We thus obtain $\trb\ph*\bigl(\tr\psi*\trb\chi\bigr) =
\bigl(\trb\ph*\tr\psi\bigr)*\trb\chi$ by passing to the limit in the
corresponding identity for integrable singularities.
\eopf

\sm

Observe that we have two subalgebras without unit $\bem\bigl(\C\{\ze\}\bigr)
\subset \SING\integ \subset \SING$, and that simple singularities form a
subalgebra
$\SING\simp = \C\,\de \oplus \bem\bigl(\C\{\ze\}\bigr) \subset \SING$.
Here are a few properties of the algebra~$\SING$:

\begin{enumerate}

\item
The family of singularities $\bigl(\trn I_\sig\bigr)_{\sig\in\C}$ defined
by~\refeq{eqdefIsig} and
$$
\trn I_{-n} = \de\ss n, \qquad \trn I_{n+1} = \bem\bigl(\tfrac{\ze^n}{n!}\bigr),
\qquad n\in\N
$$
satisfies $\trn I_{\sig_1} * \trn I_{\sig_2} = \trn I_{\sig_1+\sig_2}$ for all
$\sig_1,\sig_2\in\C$, as can be checked from the integrable case by analytic
continuation\footnote{	\label{footanafam}
There is a notion of singularity~$\trb\ph_{\!s}$ depending analytically on a
parameter $s\in S$, where $S$ is an open subset of~$\C$: following \cite[Vol.~1,
p.~48]{Eca81}, we assume that for each $s_0\in S$, there exist $r>0$, an open
subset~$V$ of~$S$ containing~$s_0$ and a holomorphic function
$\chb\ph(s,\ze)=\chb\ph_{\!s}(\ze)$ on~$V\times D_r$, where
$D_r = \ao \ze\in\Clog \,;\; |\ze|<r \af$,
such that $\trb\ph_{\!s} = \sing_0\bigl(\chb\ph_{\!s}\bigr)$ for each $s\in V$.
According to footnote~\ref{footIsig}, the family 
$\bigl(\trn I_\sig\bigr)_{\sig\in\C}$ satisfies this with $s=\sig\in S=\C$, and
this example shows that there may be no major~$\chb\ph_{\!s}(\ze)$ which
satisfies the above with $V=S$.

With this definition, the uniqueness of the continuation of analytic identities
is guaranteed by the following fact: 
if $S$ is connected and $Z=\ao s\in S\mid \trb\ph_{\!s} = 0 \af$ has an accumulation
point, then $Z=S$.
(Let $Z'$ denote the interior of~$Z$. We first observe that if a
non-stationary sequence of points $s_n\in Z$ converges to $s_\infty\in S$, then
$s_\infty\in Z'$.
Indeed, let $V$ be an open connected neighbourhood of~$s_\infty$ in~$S$, $r>0$ and
$\chb\ph_{\!s}(\ze)$ be a major of~$\trb\ph_{\!s}$ which is analytic in $(s,\ze)\in
V\times D_r$. For $n\ge N$ large enough, the functions $\chb\ph_{\!s_n}$ are regular at the
origin; for any fixed $\ze\in D_r$, the analytic identity 
$\chb\ph_{\!s}(\ze)=\chb\ph_{\!s}(\ze\,\ee^{-2\pi\I})$ holds for $s=s_n$, $n\ge N$, thus
for all $s\in V$; now each analytic function
$s\in V \mapsto c_k(s) = \int \ze^k\chb\ph_{\!s}(\ze)\,\dd\ze$, where $k\in\N$ and
the integral is taken over a circle centred at the origin, vanishes for $s=s_n$,
$n\ge N$, thus for all $s\in V$; hence every $\chb\ph_{\!s}$, $s\in V$, is regular
at the origin, \ie $\trb\ph_{\!s}=0$, and $s_\infty\in Z'$.
The open subset~$Z'$ of~$S$ is thus non-empty and closed, and the conclusion
follows from the connectedness of~$S$.)

Moreover, if we are given two families of singularities $\trb\ph_{\!s}$
and~$\tr\psi_{\!s}$ which depend analytically on $s\in S$,
formula~\refeq{eqconvolmaj} with $\la\in D_{r/2}$ and its analytic continuation for
$\ze\in D_{|\la|}$ show that $\trb\ph_{\!s} * \tr\psi_{\!s}$ also depends analytically
on $s\in S$.
One can use these facts to continue analytically the identity 
$\trn I_{\sig_1} * \trn I_{\sig_2} - \trn I_{\sig_1+\sig_2} = 0$
from $\RE\sig_1>0$ to arbitrary $\sig_1\in\C$ with a fixed $\sig_2$ of positive
real part, and then from $\RE\sig_2>0$ to arbitrary $\sig_2\in\C$ with any
fixed~$\sig_1$.
}
in~$(\sig_1,\sig_2)$.

\item
If $\al(\ze)$ is a regular germ, the multiplication of majors by~$\al$ obviously
passes to the quotient:
$$
\trb\ph = \sing_0\bigl(\chb\ph\bigr) \in \SING \;\mapsto\;
\al \trb\ph := \sing_0\bigl(\al \chb\ph\bigr) \in \SING.
$$
This turns $\SING$ into a $\C\{\ze\}$-module.

\item
In particular, we have a linear operator of $\SING$
$$
\pa:\,\trb\ph \mapsto -\ze\trb\ph,
$$ 
which turns out to be a derivation (multiply by $\dze=\dze_1+\dze_2$ in
formula~\refeq{eqconvolmaj}), the kernel of which is $\C\,\de$.

\item
Differentiation of majors passes to the quotient and defines a linear operator 
$\frac{\dd}{\dd\ze}$ which coincides with convolution by $\de\ss1 =
\frac{\dd}{\dd\ze} \de$ 
(differentiate the relation $\chb\ph_1(\ze) = \int_{\ga_{\la,\ze}}
\chb\ph_1(\ze_1) \chb\ph_2(\dze-\dze_1)\,\dd\ze_1$, where $\chb\ph_2(\ze_2) =
1/2\pi\I\ze_2$) 
and is invertible, with inverse
$\bigl(\frac{\dd}{\dd\ze}\bigr)\ii\trb\ph = \de\ss{-1}*\trb\ph$.
More generally
$$
\bigl(\tfrac{\dd}{\dd\ze}\bigr)^n\trb\ph = \de\ss{n}*\trb\ph, 
\qquad \trb\ph\in\SING, \quad n\in\Z.
$$
Notice that $\frac{\dd}{\dd\ze}$ is not a derivation; its action on a
convolution product is given by
$$
\bigl(\tfrac{\dd}{\dd\ze}\bigr)^n \bigl(\trb\ph_1*\trb\ph_2\bigr) =
\bigl(\bigl(\tfrac{\dd}{\dd\ze}\bigr)^n \trb\ph_1\bigr)*\trb\ph_2 =
\trb\ph_1*\bigl(\bigl(\tfrac{\dd}{\dd\ze}\bigr)^n \trb\ph_2\bigr).
$$
We have for instance
$\bigl(\tfrac{\dd}{\dd\ze}\bigr)^n \trn I_\sig = \trn I_{\sig-n}$
for all $\sig\in\C$ and $n\in\N$.

\end{enumerate}

\subsub{Extensions of the formal Borel transform}	\label{secextBT}

The formal Borel transform that we have used so far was defined on the space of
Gevrey-$1$ formal series~$\C[[z\ii]]_1$.
%
In the language of singularities, this means that we have an algebra isomorphism
$$
\cB\,: \ti\ph = \sum_{n\ge0} c_n z^{-n} \in \C[[z\ii]]_1 \;\mapsto\;
\trb\ph = c_0\,\de + \bem\htb\ph \in \SING\simp,
\qquad \htb\ph(\ze) = \sum_{n\ge1} c_n \tfrac{\ze^{n-1}}{(n-1)!}.
$$
The field of fractions of~$\C[[z\ii]]_1$ is $\fracC = \C[[z\ii]]_1[z]$,
the space of sums of a polynomial in~$z$ and a Gevrey-$1$ series in~$z\ii$
(because, when $c_0\neq0$, the above~$\ti\ph$ admits a multiplicative inverse
in~$\C[[z\ii]]_1$), and it is natural to extend the formal Borel transform by
setting $\cB(z^n) = \de\ss n$: we get an algebra isomorphism
$$
\cB\,: \ti\ph = \sum_{n\ge-N} c_n z^{-n} \in \fracC \;\mapsto\;
\trb\ph = \sum_{k=0}^N c_{-k}\,\de\ss k + \bem\htb\ph \in \SING\sram,
$$
with $N$ depending on~$\ti\ph$ and the same $\htb\ph(\ze)$ as above, and where
$\SING\sram\subset\SING$ is the subalgebra of ``simply ramified singularities'',
consisting of those singularities which admit a major of the form
$P(1/\ze)+\frac{1}{2\pi\I}\htb\ph(\ze)\log\ze$ with $P$ polynomial and
$\htb\ph(\ze)\in\C\{\ze\}$
(notice that this subalgebra is smaller than $\var\ii\bigl(\C\{\ze\}\bigr)$,
which consists of those singularities which admit a major of the same form but
with $P$ entire function).

\sm

In~$\fracC$, we have well-defined difference operators 
$\ti\ph(z) \mapsto \ti\ph(z+1)-\ti\ph(z)$
and 
$\ti\ph(z) \mapsto \ti\ph(z+1)-2\ti\ph(z)+\ti\ph(z-1)$,
the counterpart of which are 
$\trb\ph \mapsto \bigl(\ee^{-\ze}-1\bigr)\trb\ph$ and
$\trb\ph \mapsto 4\sinh^2\bigl(\frac{\ze}{2}\bigr)\trb\ph$.

\sm

The inverse formal Borel transform is not defined in the space of all
singularities, but further extensions are possible beyond $\SING\sram$. For
instance, setting
$$
\cB(z^{-\sig}) = \trn I_\sig, \qquad 
\cB\bigl((-1)^m z^{-\sig}(\log z)^m\bigr) = \trn J_{\sig,m} := 
\bigl(\tfrac{\pa}{\pa\sig}\bigr)^m \trn I_\sig, 
\qquad \sig\in\C,\; m\in\N,
$$
allows one to deal with formal expansions involving non-integer powers
of~$z$ and integer powers of~$\log z$
(\cf footnote~\ref{footanafam} for the differentiation \wrt\ a parameter in an
analytic family of singularities).
In practice, when studying the formal solutions of a problem, one chooses a
suitable subset of~$\SING$ according to one's needs.
This choice can be dictated by the shape of the formal solutions one finds and
of their formal Borel transforms, and also by the nature of the singularities of
the analytic continuation of the minors of these Borel transforms.

\subsub{Laplace transform of majors}

Let us denote by~$R_\th$, for any $\th\in\R$, the ray
$\left]0,\ee^{\I\th} \infty\right[$ in~$\Clog$. 
If the minor~$\htb\ph$ of a singularity $\trb\ph =
\sing_0\bigl(\chb\ph\bigr)$ extends analytically along~$R_\th$ 
%
and 
has at most exponential growth at infinity in this direction, one can define the
Laplace transform by
$$
\bigl(\cL^\th\,\trb\ph\bigr)(z) =
\int_{a\,\ee^{\I(\th-2\pi)}}^{a\,\ee^{\I\th}} 
\ee^{-z\ze}\,\chb\ph(\ze)\,\dd\ze 
+ \int_{a\,\ee^{\I\th}}^{\ee^{\I\th} \infty}
\ee^{-z\ze}\,\htb\ph(\ze)\,\dd\ze,
$$
where $a>0$ is chosen small enough and the first integral is taken over a circle
centred at the origin; the result then does not depend on the choice of~$a$ nor
on the chosen major, it is a function holomorphic in a half-plane of the form
$\RE(z\,\ee^{\I\th})>\tau$.
Observe that if $\trb\ph$ is an integrable singularity, one can let~$a$ tend
to~$0$, which yields the usual formula:
$$
\cL^\th\bigl(\bem\htb\ph\bigr)(z) =
\int_{0}^{\ee^{\I\th} \infty}
\ee^{-z\ze}\,\htb\ph(\ze)\,\dd\ze.
$$
The idea is in fact very similar to the one which was used to define the
convolution of general singularities: integration of the minor up to the origin,
being possible only in the integrable case, must be replaced by an integration
of a major around the origin.
Thus extended, the Laplace transform is compatible with the convolution of
general singularities: 
$\cL^\th\bigl(\trb\ph*\tr\psi\bigr) =
\bigl(\cL^\th\,\trb\ph\bigr)\bigl(\cL^\th\,\tr\psi\bigr)$.
If the singularity admits a major which has at most exponential growth
along~$R_\th$ and~$R_{\th-2\pi}$, one can then resort to the Laplace transform
of majors:
$$
\bigl(\cL^\th\,\trb\ph\bigr)(z) =
\left(
- \int_{a\,\ee^{\I(\th-2\pi)}}^{\ee^{\I(\th-2\pi)}\infty} 
+ \int_{a\,\ee^{\I(\th-2\pi)}}^{a\,\ee^{\I\th}} 
+ \int_{a\,\ee^{\I\th}}^{\ee^{\I\th} \infty}
\right) \ee^{-z\ze}\,\chb\ph(\ze)\,\dd\ze
$$
(with the second integral taken over the same circle as above, \ie the sum of
the three parts amounts to a Hankel's contour integral).

\sm

Laplace transforms in nearby directions $\th_1<\th_2$ can be glued together and
yield a function~$\cL^{\left]\th_1,\th_2\right[}\,\trb\ph$ analytic in a sectorial
neighbourhood of infinity when the minor has no singularity in the sector
$\th_1\le\arg\ze\le\th_2$; it is then more convenient to consider~$z$ as element
of~$\Clog$, with $-\th_2 - \frac{\pi}{2} < \arg z < -\th_1 + \frac{\pi}{2}$ 
and~$|z|$ large enough.
The difference with what we saw at the beginning (Section~\ref{secFBT}) is in
the possible asymptotic expansions of~$\cL^\th\,\trb\ph$.

\sm

In the simply ramified case, when a Laplace
transform~$\cL^{\left]\th_1,\th_2\right[}\,\trb\ph$ can be defined for
$\trb\ph\in\SING\sram$, one finds $\cB\ii\trb\ph\in\fracC$ as asymptotic
expansion.
Similarly, if the case of one of the extensions of the formal Borel transform
mentioned in the previous section,
$\cL^{\left]\th_1,\th_2\right[}\,\trb\ph(z)\sim\cB\ii\trb\ph(z)$, as a consequence
of the formula
$$
(-1)^m z^{-\sig}(\log z)^m = \cL^\th\, \trn J_{\sig,m}(z), \qquad
z\in\Clog, \quad -\th-\tfrac{\pi}{2} < \arg z < -\th+\tfrac{\pi}{2}.
$$

\subsection{General resurgent functions and alien derivations}	\label{secGenRes}

We are now ready to give definitions which are more general than in
Sections~\ref{secformconvmod} and~\ref{secsimpl}.
We shall not provide many details; the reader is referred to
\cite[Vol.~3]{Eca81}, \cite{dulac}, \cite{Eca93} or \cite{CNP}.

\med

We first define the space of ``resurgent minors'', $\whRES_{2\pi\I\Z}$, as the
set of all the germs of~$\ANA$ which extend to the universal cover
$\widetilde{(\C\setminus2\pi\I\,\Z,1)}$ (using the notation of
footnote~\ref{footClog}, meaning that the holomorphic function~$\htb\ph$
determined by the germ in a spiralling neighbourhood of the origin~$\cV$ extends
analytically along any path of~$\Clog$ which starts in~$\cV$ and avoids the lift
of~$2\pi\I\,\Z^*$ in~$\Clog$).
Some resurgent minors are integrable minors, among these some are even regular
at the origin; this gives rise to subspaces $\whRES_{2\pi\I\Z}\cap\ANA\integ$
and $\wh\cH(\cR) = \whRES_{2\pi\I\Z}\cap\,\C\{\ze\}$, which are both stable by
convolution (the former for reasons similar to what was explained in
Section~\ref{secAnContConv} for the latter).

Next, we define the ``convolutive model of resurgent functions'' as the space of
all the singularities of~$\SING$, the minors of which belong to
$\whRES_{2\pi\I\Z}$:
$$
\trRES_{2\pi\I\Z} := \var\ii\bigl( \whRES_{2\pi\I\Z} \bigr) \subset \SING.
$$
This space is stable by convolution (we omit the proof): {\em
$\trRES_{2\pi\I\Z}$ is a subalgebra of~$\SING$}, which obviously contains the
unit~$\de$.
We may call 
$$
\trRESint{2\pi\I\Z} = \bem\bigl( \whRES_{2\pi\I\Z} \cap \ANA\integ \bigr)
\subset \trRES_{2\pi\I\Z}
$$
the ``minor model''; it is a subalgebra of $\SING\integ$ (without unit), the
elements of which are determined by their minor, so that there is no loss in
information when reasoning on the minors only.
The convolution algebra $\wh\cH(\cR)$ of Section~\ref{secformconvmod}, being
isomorphic to $\bem\bigl( \whRES_{2\pi\I\Z} \cap \,\C\{\ze\} \bigr)$,
can now be considered as a subalgebra of~$\trRESint{2\pi\I\Z}$.

\med

The algebra $\Rsimp$ of Section~\ref{secsimpl} 
\begin{equation}	\label{eqdeftrRESsimp}
\Rsimp \simeq \trRESsimp{2\pi\I\Z} := \ao 
\trb\ph \in \trRES_{2\pi\I\Z} \cap \SING\simp \mid
\text{$\var\bigl(\trb\ph\bigr)$ has only simple singularities} \af
\end{equation}
corresponds to singularities which are determined by a regular minor up to the
addition of a multiple of~$\de$, such that the minor extends to~$\cR$, and {\em
with the further restriction that the analytic continuation of the minor
possesses only simple singularities.}
This restriction made it possible to define the alien derivations~$\De_\om$ in
Section~\ref{secsimpl} as internal operators of~$\Rsimp$. 
Relaxing the conditions of regularity at the origin and on the shape of the
singularities, we are now in a position to define the alien
derivations~$\De_\om$ in a somewhat enlarged framework, as internal operators
of~$\trRES_{2\pi\I\Z}$ and not only of~$\trRESsimp{2\pi\I\Z}$.

\sm

Because of the possible ramification of the minor at the origin, the alien
derivations will now be indexed by all $\om\in\Clog$ such that $\dom \in
2\pi\I\,\Z^*$. 
Here is the generalisation of Definition~\ref{defDeom}, with notations similar
to those of Section~\ref{secsimpl}:
\begin{quote}{\em
For $\trb\ph\in\trRES_{2\pi\I\Z}$ and $\om = 2\pi m\,\ee^{\I\th}$ with
$m\in\N^*$ and $\th\in\frac{\pi}{2}+\pi\Z$, 
\begin{equation}	\label{eqdefDeomsing}
\De_\om\trb\ph = \sum_{\eps_1,\ldots,\eps_{m-1}\in\{+,-\}}
\frac{p(\eps)!q(\eps)!}{m!} \sing_0 \bigl( \chn\Phi_{\ga(\eps)} \bigr),
\end{equation}
where the path~$\ga(\eps)$ connects $\left]0,\frac{1}{m}\om\right[$ and
$\left]\frac{m-1}{m}\om,\om\right[$ and circumvents 
$2\pi r\,\ee^{\I\th} = \frac{r}{m}\om$
to the right if $\eps_r=+$ and to the left if $\eps=-$, and where the analytic
continuation of the minor~$\htb\ph$ of~$\trb\ph$ determines the major
$$
\chn\Phi_{\ga(\eps)}(\ze) = \bigl(\cont_{\ga(\eps)} \htb\ph\bigr)(\om+\dze),
\qquad \arg\om-2\pi < \arg\ze < \arg\om, \ens |\ze| < 2\pi.
$$
}\end{quote}

One can check that the operators~$\De_\om$ are {\em derivations}
of~$\trRES_{2\pi\I\Z}$, which satisfy the rules of alien
calculus that we have indicated in the case of simple resurgent functions.
Notice that if $\dom_1=\dom_2$, the restrictions of~$\De_{\om_1}$
and~$\De_{\om_2}$ to $\var\ii\bigl(\wh\cH(\cR)\bigr)$, and a fortiori
to~$\trRESsimp{2\pi\I\Z}$, coincide.

\med

Up to now we have restricted ourselves to minors which extend analytically
provided one avoids always the same fixed set of potentially singular points,
namely $2\pi\I\,\Z$ (or its lift in~$\Clog$).
But one can consider other lattices $\Om\subset\C$ of singular points and
define accordingly the space $\whRES_\Om$ of germs which extend to
$\widetilde{(\C\setminus\Om,1)}$, the algebra $\trRES_\Om = \var\ii\bigl(
\whRES_\Om \bigr)$, and the alien derivations~$\De_\om$ with $\dom\in\Om$
($\Om$ must be an additive semi-group of~$\C$ to ensure stability by
convolution, and a group to ensure stability by alien derivations\footnote{
If $\dom\in\Om$ and $\trb\ph\in\trRES_\Om$, some branchs of the minor
of~$\De_{2\om}\trb\ph$ will usually be singular at~$-\dom$, which should thus be
included in~$\Om$
(for instance, with the minor $\htb\ph(\ze) = \frac{\om}{\ze-\om}\log\bigl( 1 -
\frac{\ze}{2\om} \bigr)$, which is regular at the origin, $\De_{2\om}\trb\ph =
\bem\bigl( 1 + \frac{\ze}{\om} \bigr)$).
}).
The subalgebra of simple resurgent functions with singular support in~$\Om$ will be
defined by
$$
\trRESsimp{\Om} = \ao \trb\ph \in \trRES_\Om \cap \SING\simp \mid
\forall r\ge1, \forall \om_1,\ldots, \om_r, \; 
\De_{\om_r}\cdots\De_{\om_1}\trb\ph \in \SING\simp \af
$$
(which is consistent with~\refeq{eqdeftrRESsimp}---\cf
footnote~\ref{footnotlimit}).
Similarly, replacing~$\SING\simp$ by~$\SING\sram$ in the above formula, one can
define the larger subalgebra $\trRsramOm$ of simply ramified resurgent functions
with singular support in~$\Om$.

\sm

Finally, the most general algebra we can construct at this level is the
space~$\trRES$ of all the singularities, the minors of which are {\em endlessly
continuable along broken lines}: a singularity~$\trb\ph$ is said to belong
to~$\trRES$ if, on any broken line~$L$ of finite length drawn on~$\Clog$ and
starting in the domain of analyticity of the minor~$\htb\ph$, there exists a
finite set~$\Om_L$ (depending on~$\htb\ph$) such that~$\htb\ph$ admits an
analytic continuation along the paths which follow~$L$ but circumvent the points
of~$\Om_L$ to the right or to the left.

This means that we do not impose the location of possibly singular points in
advance, nor any constraint on the shape of the possible singularities.
The alien derivations~$\De_\om$ acting in~$\trRES$ are thus indexed by all
$\om\in\Clog$; they are defined by the same formula as~\refeq{eqdefDeomsing}, but
with $m-1 = \card\Om_L$, for a segment $L = [\tau\om, (1-\tau)\om]$ with $\tau>0$
small enough so that the points of $\Om_L=\{\om_1,\ldots,\om_{m-1}\}$ be the
only singular points encountered in the analytic continuation of the minor
along~$\left]0,\om\right[$ (instead of
$\bigl\{\frac{1}{m}\om,\ldots,\frac{m-1}{m}\om\bigr\}$), and with $\chn\Phi_{\ga(\eps)}$
denoting the $2^{m-1}$ corresponding branchs of the minor near~$\om$
(of course, if none of them is singular at~$\om$, then $\De_\om\trb\ph=0$).
One can check that the modified formula defines an operator~$\De_\om$ which is a
derivation of the convolution algebra~$\trRES$.

\subsub{Bridge equation for non-degenerate parabolic germs in the case $\rho\neq0$}

As an illustration of this enlarged formalism with more general resurgent
functions than simple ones, let us return to non-degenerate parabolic germs with
arbitrary resiter, as defined in Section~\ref{secAbEqtgtId}:
the holomorphic germ~$F$ at the origin gives rise to a germ at infinity
$$
f(z) = z + 1 + a(z), \qquad
a(z) = - \rho z\ii + \cO(z^{-2}) \in \C\{z\ii\},
$$
with $\rho\in\C$.
As was mentioned after Proposition~\ref{propformalconj}, Abel's equation $v\circ
f = v + 1$ admits a formal solution
$$
\ti v(z) = z + \ti\psi(z), \qquad
\ti\psi = \rho\log z + \sum_{n\ge0} c_n z^{-n},
$$
which is unique if we impose $c_0=0$. This ``iterator'' is formally invertible, with
inverse
$$
\ti u(z) = z + \ti\ph(z), \qquad
\ti\ph(z) = -\rho\log z + \sum_{\substack{n,m\ge0\\ n+m\ge1}} C_{n,m} z^{-n} (z\ii\log z)^m.
$$
The inverse iterator allows one to describe all the solutions of the
difference equation $f\circ u(z) = u(z+1)$ in the set
$\{ z - \rho\log z + \chi(z), \; \chi(z) \in \C[[z\ii,z\ii\log z]] \}$:
they are the series of the form $\ti u(z+c)$ with arbitrary $c\in\C$.

\sm

The space $\C[[z\ii,z\ii\log z]]$ is one of those spaces to which the formal
Borel transform can be extended (\cf Section~\ref{secextBT}).
For the Borel transform of $J_{\sig,m}(z) = (-1)^m z^{-\sig}(\log z)^m$, we have
the formula
$$
\trn J_{\sig,m} = \bem\htn J_{\sig,m}, \qquad
\htn J_{\sig,m}(\ze) = \sum_{k=0}^m  
\bino{m}{k} \Bigl(\frac{1}{\Ga}\Bigr)^{(k)}\!\!(\sig)\, \ze^{\sig-1}(\log\ze)^{m-k},
\qquad \RE\sig>0,\quad m\in\N.
$$
The analogue of Theorem~\ref{thmresur} is that the series 
$$
\hta\psi_1(\ze) = \sum c_n \tfrac{\ze^{n-1}}{(n-1)!}, \qquad
\htb\ph_1(\ze) = \sum (-1)^m C_{n,m} \htn J_{n+m,m}(\ze)
$$ 
converge for $|\ze|$ small enough and define germs of~$\ANA$, which belong
to~$\whRES_{2\pi\I\Z}$, with at most exponential growth at infinity in the
non-vertical directions.
Thus
$$
\trx v = \de' - \rho\trn J_{0,1} + 
\bem\hta\psi_1
\in \trRES_{2\pi\I\Z}, \qquad
\trx u = \de' + \rho\trn J_{0,1} + 
\bem\htb\ph_1
\in \trRES_{2\pi\I\Z}.
$$

The analogue of Theorem~\ref{thmParBridge} is the existence of complex
numbers~$A_\om$, the ``analytic invariants'' of~$f$, such that, when pulled back in the
$z$-variable, the action of the alien derivations is given by
$$
\De_\om \ti u = A_\om \,\frac{\dd\ti u}{\dd z}, \qquad
\De_\om \ti v = -A_\om \,\ee^{-\om(\ti v(z)-z)}.
$$
Equivalently, we can say that $\ti\ph,\ti\psi$ have formal Borel tansforms
$\trb\ph,\tr\psi\in\trRES_{2\pi\I\Z}$ and satisfy
$$
\De_\om \ti\ph = A_\om \bigl(1 + \frac{\dd\ti\ph}{\dd z} \bigr), \qquad
\De_\om \ti\psi = -A_\om \,z^{-\rho\om}\, \ee^{-\om\ti\psi_1}.
$$
The successive alien derivatives can also be computed, giving rise to formulas
analogous to~\refeq{eqiterBridgeu} and~\refeq{eqiterBridgev}.
This means in particular that, near~$\om$, any branch of~$\htb\ph(\ze)$
is of the form $\frac{B}{2\pi\I(\ze-\om)} + \chb\chi(\ze-\om) +
\text{reg}(\ze-\om)$, with a complex number~$B$ and an integrable singularity
$\sing_0\bigl(\chb\chi(\ze)\bigr)$ which can be computed from~$\ti\ph$ and from
the invariants for each chosen branch),
and any branch of~$\hta\psi(\ze)$ is of the form 
$ -B \bigl( 
\chn I_{\rho\om}(\ze) + b_1\chn I_{\rho\om+1}(\ze) + b_2\chn I_{\rho\om+2}(\ze) 
+ \cdots \bigr) + \text{reg}(\ze-\om)$ with computable complex
numbers~$b_1,b_2,\ldots$
When $\rho\neq0$, neither $\trb\ph,\tr\psi$ nor their alien derivatives\footnote{
Except the ones corresponding to $\om$ such that $A_\om=0$, in case certain
invariants vanish.
}
are simple singularities
(and $\tr\psi_1=\bem\hta\psi_1$ is a simple singularity, but
$\De_\om \tr\psi_1 = \De_\om \tr\psi$ is not).

\sm

In the above, $\om$ is any point of~$\Clog$ with $\dom\in 2\pi\I\Z^*$, but it
turns out that the numbers~$A_\om$ depend on~$\dom$ only, since 
$\htb\psi = -\rho\htn J_{0,1} + \htb\psi_1$ where $\htb\psi_1(\ze)\in\C\{\ze\}$
and $\htn J_{0,1}(\ze) = 1/\ze$ carries no singularity outside the
origin.\footnote{
Using the formula for~$\chn I_\sig$ indicated in footnote~\ref{footIsig}, one
finds 
$\trn J_{0,1} = \sing_0\bigl( \frac{\log\ze + \ga + \I\,\pi}{2 \pi\I \ze} \bigr)$,
where $\ga$ is Euler's constant.
}
This fact can also be deduced from the existence of a resurgent $\ti U(z)\in \ID +
z\ii\C[[z\ii]]$ such that $\ti u(z) = \ti U\bigl(z+\chi(z\ii,z\ii\log z)\bigr)$ 
with $\chi(u,v)\in\C\{u,v\}$ (the alien derivatives $\De_\om\ti U$ depend on the
projections~$\dom$ only, since $\htn U(\ze)\in\C\{\ze\}$, and the ``alien chain
rule'' yields 
$\De_\om\ti u = \De_\om \bigl(\ti U\circ(\ID+\chi)\bigr) =
\ee^{-\dom\chi} \, (\De_\om\ti U)\circ(\ID+\chi)$).
The series~$\ti U$ and~$\chi$ can be found by using the operator of ``formal
monodromy'', \ie the substitution $z\mapsto z\,\ee^{2\pi\I}$ in the solution~$\ti
u(z)$ defined by $z^{-n}(\log z)^m\mapsto z^{-n}(\log z+2\pi\I)^m$, which leads
to another solution $\ti u(z\,\ee^{2\pi\I}) = z - \rho\log z - 2\pi\I\rho +
o(1)$, whence $\ti u(z\,\ee^{2\pi\I}) = \ti u(z-2\pi\I\rho)$ (because all the
solutions are known to be of the form $\ti u(z+c)$);
one then observes that the (convergent) transformation $z_*\mapsto
z=z_*+\rho\log z_*$, $\ti u(z)=\ti U(z_*)$, transforms this monodromy relation
into the trivial one $\ti U(z_*\,\ee^{2\pi\I}) = \ti U(z_*)$ and has a convergent
inverse of the form $z_*=z+\chi(z\ii,z\ii\log z)$.

\sm

The reader is referred to~\cite[Vol.~2]{Eca81} for the results of this section.


\section{Splitting problems}


\subsection{Second-order difference equations and complex splitting problems}	\label{secHenon}

We now wish to present the results of the article~\cite{sauzin01}, and hint at
some of the more general results to be found in the work in
progress~\cite{GelfSauz}, in the context of $2$-dimensional holomorphic
transformations with a parabolic fixed point.

\sm

The aim is to understand a part of the local dynamics for a germ
\begin{equation}	\label{henonGen}
F : \bino{x}{y} \longmapsto
\bino{x + y + f(x,y)}{y + f(x,y)},
\end{equation}
where $f(x,y)\in\C\{x,y\}$ is of the form
\begin{equation*}	
f(x,y) = -x^2 -\ga xy + \cO_3(x,y), \qquad \ga\in\C.
\end{equation*}
It is shown in~\cite{GelfSauz} that any germ of holomorphic map of~$(\C^2,0)$
with double eigenvalue~$1$ but non-identity differential at the origin can be
reduced, by a local analytic change of coordinates, to the form~\refeq{henonGen}
with $f(x,y) = b_{20} x^2 + b_{11} xy + \cO_3(x,y)$; in the non-degenerate case,
\ie when $b_{20}\neq0$, it is easy to normalise further the map so as to have
$b_{20}=-1$ (observe that one can always remove the $y^2$-term, but it is not so
for the $xy$-term: the complex number~$\ga$ is in fact a {\em formal invariant}
of the non-degenerate parabolic germ under consideration).

\sm

The article~\cite{sauzin01} is devoted to the particular map one obtains when
$f(x,y)=-x^2$, which is an instance of the H\'enon map (for a special choice of
parameters). It is an invertible quadratic map $\C^2\to\C^2$, the only fixed
point of which is the origin, and which is symplectic for the standard symplectic
structure $\dd x\wedge \dd y$.

\subsub{Formal separatrix}

Our main subject of investigation will be the {\em formal separatrix} of the
map~\refeq{henonGen}, which is a pair of formal series $\ti p(z) = \bigl(\ti
x_0(z),\ti y_0(z)\bigr)$ satisfying $\ti p(z+1) = F\bigl(\ti p(z)\bigr)$ and formally
asymptotic to the origin, and its Borel sums~$p^+(z)$ and~$p^-(z)$. 
A slightly more geometric way of introducing~$\ti p(z)$ is to consider first the
{\em formal infinitesimal generator} of~$F$.

\sm

Indeed, it turns out that there exists a unique formal vector field 
$$
\ti X = \bigl( y + \ti A(x,y) \bigr) \frac{\pa}{\pa x} 
+ \ti B(x,y) \frac{\pa}{\pa y}, \qquad
\ti A(x,y),\ti B(x,y) = \cO_2(x,y) \in \C[[x,y]],
$$
the formal flow of which $\ti\Phi(t,\,\cdot\,,\,\cdot\,) = \exp_t\ti X$
satisfies $\exp_1\ti X = F$
(the flow of a formal vector field like~$\ti X$ is determined as the unique 
$\ti\Phi\in\bigl(\C[t][[x,y]]\bigr)^2$ such that $\ti\Phi_{|t=0}=\ID$ and
$\pa_t\ti\Phi = \ti X\circ\ti\Phi$).
This is the $2$-dimensional analogue of the infinitesimal generator mentioned in
footnote~\ref{footiter}.
Now, for any $2$-dimensional vector field with a singularity at the origin,
there exists at least one ``separatrix'', \ie a solution of the vector field
which is asymptotic to the origin. This is the celebrated Camacho-Sad theorem
for analytic vector fields; it is thus not a surprise that for our formal vector
field there exists a formal solution~$\ti p(z)$ which is formally asymptotic to
the origin.
In practice, one finds that both components of~$\ti p(z)$ belong to the space
$z^{-2}\C[[z\ii]]$ when $\ga=0$, and to the space $z^{-2}\C[[z\ii,z\ii\log z]]$
in the general case; moreover this solution is unique (under our non-degeneracy
hypothesis) up to a time-shift $z\mapsto z+a$.

\sm

In some sense the dominant part in~$\ti X$ is
$y \frac{\pa}{\pa x} - x^2 \frac{\pa}{\pa y}$, which is the Hamiltonian vector
field generated by $h(x,y)=\demi y^2 + \frac{1}{3}x^3$ for the standard symplectic
structure (however the whole vector field~$\ti X$ itself is Hamiltonian only
when~$F$ is symplectic, \ie when the function~$f(x,y)$ depends on its first
argument only).
The separatrix for this dominant part is given by the cusp $\{h(x,y)=0\}$ (zero
energy level), the time-parametrisation of which is $(-6z^{-2},12z^{-3})$.
This is the leading term of the formal separatrix~$\ti p(z)$.

\sm

It turns out that the formal separatrix $\ti p(z)=\bigl(\ti x_0(z),\ti
y_0(z)\bigr)$ can be found directly from~$F$, without any reference to the
formal vector field~$\ti X$, \ie without solving the formal differential
equation $\pa_z\ti p(z) = \ti X\bigl(\ti p(z)\bigr)$. 
One just needs to consider the difference equation
\begin{equation}	\label{eqHenonSys}
\ti p(z+1) = F\bigl(\ti p(z)\bigr)
\end{equation}
(in fact the equations $\ti p(z+t)=\exp_t\ti X\bigl(\ti p(z)\bigr)$ with $t\neq0$
are all equivalent).
This vector equation, in turn, is equivalent to the scalar equations 
$$
P \ti x_0 = f(\ti x_0, D \ti x_0), \qquad \ti y_0 = D \ti x_0,
$$
with difference operators~$P$ and~$D$ defined by
\begin{equation}	\label{eqdefDP}
P\ph(z) = \ph(z+1) - 2\ph(z) + \ph(z-1),
\qquad D\ph(z) = \ph(z) - \ph(z-1)
\end{equation}
(thanks to the special form~\refeq{henonGen} that we gave to the map~$F$).
One can thus eliminate~$\ti y_0$ and work with the equation $P \ti x_0 = f(\ti
x_0, D \ti x_0)$ alone, which is nothing but the nonlinear second-order difference
equation~\refeq{eqnsecond} of Section~\ref{secdiffeq} (up to a slight change of
notation).
The equation which corresponds to the H\'enon map and is studied
in~\cite{sauzin01} is simply $P\ti x_0 = - \ti x_0^2$.

\sm

When one is given an analytic vector field, the corresponding separatrix is
convergent: the series~$\ti x_0(z)$ and~$\ti y_0(z)$ converge for~$|z|$ large
enough. One should not expect convergence for a general map~$F$. 
The case of an entire map, like the H\'enon map, is of particular interest, as
one can prove divergence in this case.
We shall see that the components of the separatrix are always resurgent and
generically divergent.

\sm

We speak of ``separatrix splitting'' because the separatrix, which was
convergent in the case of an analytic vector field, breaks (becomes formal) when
one passes to maps and gives rise to two distinct curves (two Borel sums of~$\ti
p(z)$, none of which is the analytic continuation of the other).

\sm

We shall now proceed and describe the results concerning~$\ti x_0(z)$ and its
resurgent structure which are given in~\cite{sauzin01} for the H\'enon map.
We shall see that one of the novelties of the $2$-dimensional case \wrt\
Section~\ref{secAlCalAbEq} is the necessity of considering a formal solution
more general than the formal series~$\ti x_0(z)$, namely the ``formal
integral''~$\ti x(z,b)$ which depends on a further formal variable and allows
one to write a Bridge equation.
The terminology comes from~\cite[Vol.~3]{Eca81}, as well as the ideas for the
resurgent approach (although the case of parabolic maps like~\refeq{henonGen} is
not covered by this reference).

\subsub{First resurgence relations}

From now on we thus set $f(x,y)=-x^2$ and we keep using the
notation~\refeq{eqdefDP} for the difference operators~$D$ and~$P$.

\begin{thm}	\label{thmResurx}
The nonzero solutions $x(z)\in\C[[z\ii]]$ of the equation
\begin{equation}	\label{eqHenon}
P x = - x^2
\end{equation}
are the formal series $x(z) = \ti x_0(z+a)$, where $a\in\C$ is arbitrary
and~$\ti x_0(z)$ is the unique nonzero even formal solution:
$$
\ti x_0(z) = -6z^{-2}+\frac{15}{2}z^{-4}-\frac{663}{40}z^{-6}+\cdots
$$
They are resurgent and Borel summable: the formal Borel transform $\htx x_0 =
\cB\ti x_0$ has positive radius of convergence and extends to~$\cR$
($\htx x_0(\ze) \in \wh\cH(\cR) = \C\{\ze\}\cap \whRES_{2\pi\I\Z}$),
with at most exponential growth at infinity along non-vertical directions.
\end{thm}

\idproof
The formal part of the statement can be obtained by substitution of 
$x(z) = \sum_{n\ge n_0} a_n z^{-n}$ with $a_{n_0}\neq0$ into~\refeq{eqHenon}:
one finds that necessarily $n_0=2$ and $a_{2}=-6$, then $a_3$ is free whereas
all the successive coefficients are uniquely determined. Choosing $a_3=0$ yields
the even solution~$\ti x_0(z)$, while the general solution must coincide with
$\ti x_0(z+\frac{a_3}{12})$.

\sm

The Borel transforms of the solutions are studied through the equation they
satisfy: the counterpart of~$P$ is multiplication by $\al(\ze) =
\ee^{-\ze}-2+\ee^\ze$, hence
\begin{equation}	\label{eqHenonBT}
\al \, \htx x = - \htx x*\htx x, \qquad
\al(\ze) = 4 \sinh^2\frac{\ze}{2}.
\end{equation}
We know in advance that this equation has a unique formal solution of the form
$\htx x_0(\ze) = -6\ze + \htx v(\ze)$ with $\htx v(\ze) \in \ze^3\C[[\ze]]$.
The corresponding equation for~$\htx v$ is
$$
\al \, \htx v - 12\ze*\htx v = 6\bigl(\ze\al(\ze)-\ze^3\bigr) - \htx v*\htx v.
$$
As in the proof of Theorem~\ref{thmresur} in Section~\ref{secAbEqtgtId}, one can
devise a method of majorants to prove that~$\htx v(\ze)$ has positive radius of
convergence and extends analytically to the sets~$\cR_c\ss0$ which were defined
there, with at most exponential growth at infinity (and, in fact, exponential
decay). The method can also be adapted to reach the union~$\cR\ss1$ of the
half-sheets which are contiguous to the principal sheet; analyticity is then
propagated to the rest of~$\cR$ through the resurgence relations to be shown
below.
See~\cite{sauzin01} for the details (or~\cite{OSS} for an analogous proof).
\eopf

Observe that the only source of singularities in the Borel plane is the division
by~$\al(\ze)$ when solving equation~\refeq{eqHenonBT}, this is why the only
possible singular points are the points of~$2\pi\I\,\Z$.

\sm

We shall see that the first singularities of~$\htx x_0(\ze)$, \ie the
singularities at $\ze=\pm2\pi\I$, are not apparent ones, thus~$\ti x_0(z)$ is
divergent.
The coefficients of~$\ti x_0(z)$ are real numbers and~$\htx x_0(\ze)$
is thus real-analytic;
therefore, the singularity at~$-2\pi\I$ can be deduced by symmetry from the
singularity at~$2\pi\I$.
We now use alien calculus to analyse the singularity at~$2\pi\I$.

\sm

Let $\trx x_0 = \bem\htx x_0$; we thus have $\trx x_0 \in
\bem\bigl(\C\{\ze\}\cap \whRES_{2\pi\I\Z}\bigr)$, solution of
\begin{equation}	\label{eqHenonSG}
\al \,\trx x_0 = - \trx x_0*\trx x_0.
\end{equation}
A major of the singularity $\trb\chi = \De_{2\pi\I}\trx x_0$ can be defined by
$\chb\chi(\ze) = \htx x_0(2\pi\I+\ze)$, $-\frac{3\pi}{2}<\arg\ze<\frac{\pi}{2}$,
$|\ze|<2\pi$.
We shall show that~$\trb\chi$ can be expressed as the linear combinations of two
elementary singularities deduced from~$\trx x_0$.
The key point is the possibility of ``alien-differentiating''
equation~\refeq{eqHenonSG}: for any $\om\in2\pi\I\,\Z$, the singularity $\trb\ph
= \De_\om\trx x_0$ must satisfy $\al \, \trb\ph = -2\trx x_0*\trb\ph$ (because
$\De_\om$ is a derivation which commutes with the multiplication
by~$\al$).\footnote{
We do not present the arguments in full rigour here. For instance, the method of
majorants we alluded to for the proof of Theorem~\ref{thmResurx} yields the
analyticity of~$\htx x_0$ in~$\cR\ss1$, hence the above~$\chb\chi$, at this
level, is only a major of ``sectorial'' singularity:
$\chb\chi\in\ANA_{\frac{\pi}{2},\de}$ for $0<\de<\frac{\pi}{2}$ with the
notation of footnote~\ref{footgensing}. The subsequent arguments should thus be
rephrased in the corresponding space~$\SING_{\frac{\pi}{2},\de}$ rather
than~$\SING$, in order to establish the relations~\refeq{eqfirstresrel}.
One would then use these relations to propagate the analyticity of~$\htx x_0$
in~$\cR\ss2$, and argue similarly to reach farther and farther half-sheets
of~$\cR$, using gradually all the resurgent relations expressed by the Bridge
Equation of Theorem~\ref{thmBridgeHen} below.
}

\begin{prop}	\label{propHenLin}
The linear difference equation
\begin{equation}	\label{eqHenLin}
P\ph = -2 \ti x_0 \ph
\end{equation}
admits a unique even solution $\ti\ph_2(z)\in\C[[z\ii]][z]$ of the form 
$\frac{1}{84}z^4\bigl(1 + \cO(z\ii) \bigr)$.
It belongs in fact to the space $\fracC = \C[[z\ii]]_1[z]$ and has a formal
Borel transform of the form
$$
\trb\ph_2 = \frac{1}{84}\de\ss4 + \frac{17}{840}\de\ss2 -\frac{17}{2240}\de +
\bem\htb\ph_2
$$ 
with $\htb\ph_2(\ze)\in\wh\cH(\cR)$.
Moreover, the solutions of the linear equation
\begin{equation}	\label{eqHenLinSG}
\al \, \trb\ph = -2\trx x_0*\trb\ph, \qquad \trb\ph\in\SING
\end{equation}
are the linear combinations (with constant coefficients) of $\trb\ph_1=\pa\trx
x_0$ and~$\trb\ph_2$.
In particular, there exist $\mu\in\R$ and $\Th\in\I\,\R$ such that 
\begin{equation}	\label{eqfirstresrel}
\De_{2\pi\I}\trx x_0 = \mu\,\trb\ph_1 + \Th\,\trb\ph_2, \qquad
\De_{-2\pi\I}\trx x_0 = -\mu\,\trb\ph_1 + \Th\,\trb\ph_2.
\end{equation}
\end{prop}

\proof
Let us first consider equation~\refeq{eqHenLin} in the
space of formal series~$\C[[z\ii]][z]$. This equation being the linearization of
equation~\refeq{eqHenon}, the ordinary derivative $\ti\ph_1 = \pa\ti x_0 = 12
z^{-3} \bigl(1 + \cO(z^{-2})\bigr)$ is obviously a particular solution (thus its
formal Borel transform $\trb\ph_1 = \pa\trx x_0$, which we know belongs to
$\bem\bigl(\C\{\ze\}\bigr) \subset \SING$, satisfies~\refeq{eqHenLinSG}).

\sm

Standard tools of the theory of second-order linear difference equations
(see~\cite{GL}) allow us to find an independent solution~$\ti\ph_2$:
$\ti\ph_2 = \ti\psi\ti\ph_1$ is solution as soon as 
$D\ti\psi=\ti\chi$, where 
$\ti\chi(z)=\frac{1}{\ti\ph_1(z)\ti\ph_1(z-1)}=\frac{1}{144}z^6\bigl(1+\cO(z\ii)\bigr)$.
The latter equation determines~$\ti\psi(z)$ up to an additive constant, since it
can be rewritten $\pa\,\be(\pa)\ti\psi = \ti\chi$, with an invertible power series
$\be(X) = \frac{1-\ee^{-X}}{X} = 1 + \cO(X)$.
We just need to choose a primitive $\pa\ii\ti\chi$ and we get the corresponding
$\ti\psi = \ga(\pa)\,\pa\ii\ti\chi$, with $\ga(X) = \frac{1}{\be(X)} = 1+\ga_1 X + \ga_2
X^2 + \cdots$.
It turns out that the primitives of~$\ti\chi(z)$ belong to~$\C[[z\ii]][z]$:
there is no $\log z$ term because the coefficient of~$z\ii$ in~$\ti\chi(z)$ is
zero (this is due to the special form of
$\ti\chi(z)=\ti\chi_1(z)\ti\chi_1(z-1)$, with $\ti\chi_1=1/\ti\ph_1$ odd: use the
Taylor formula and observe that $\ti\chi_1\ti\chi_1\ss k$ is even when $k$ is
even, and that it is the derivative of an element of~$\C[[z\ii]][z]$ when $k$ is
odd). 
These solutions~$\ti\psi\ti\ph_1$ are thus elements of~$\C[[z\ii]][z]$, of the
form $\frac{1}{84}z^4\bigl(1+\cO(z\ii)\bigr)$, differing one from the other by a
multiple of~$\ti\ph_1$, and it is easy to check that exactly one of them is our
even solution~$\ti\ph_2$. The coefficients can be determined inductively from
those of~$\ti x_0$; one finds $\ti\ph_2(z) = \frac{1}{84}z^4 + 
\frac{17}{840}z^2 -\frac{17}{2240}+\cO(z^{-2})$.

\sm

We have $\ti\ph_2\in\fracC$ and the minor of its formal Borel transform is
in~$\wh\cH(\cR)$, because the same is true of the above series~$\ti\chi$,
$\pa\ii\ti\chi$ and~$\ti\psi$.
Indeed, we can write $\ti\ph_1(z)\ti\ph_1(z-1) = 144 z^{-6} \bigl( 1 - \ti w(z) \bigr)$
with $\ti w(z) \in z\ii\C[[z\ii]]$ and 
$\htx w(\ze) = -\frac{1}{144} \bigl(\frac{\dd}{\dd\ze}\bigr)^6 
\bigl[ (\ze\htx x_0)*(\ze\ee^\ze\htx x_0) \bigr] \in\wh\cH(\cR)$; 
Proposition~\ref{propsubs} then ensures that $(1-\ti w)\ii$ is in~$\C[[z\ii]]_1$
(thus $\ti\chi(z) = \frac{1}{144}z^6\bigl(1-\ti w(z)\bigr)\ii$ and its
primitives lie in~$\fracC$) and that the minor of the formal Borel transform of
$(1-\ti w)\ii$ is in~$\wh\cH(\cR)$ (thus the minor~$\htb\chi$ of $\trb\chi=\cB\ti\chi$ lies
also in~$\wh\cH(\cR)$, and so does $\frac{1}{\ze}\htb\chi(\ze)$ which is the
minor corresponding to the primitives of~$\ti\chi$);
the conclusion for~$\ti\psi$ and~$\ti\ph_2$ follows easily (the operator~$\ga(\pa)$ 
amounts to the multiplication by $-\frac{\ze}{1-\ee^\ze}$ in the Borel plane).

\sm

We can now easily describe the solutions of equation~\refeq{eqHenLinSG}, or
even those of the inhomogeneous equation
\begin{equation}	\label{eqHenLinNHSG}
\al \, \trb\ph + 2 \trx x_0*\trb\ph = \tr\psi,
\end{equation}
with any given $\tr\psi\in\SING$.
It is sufficient to use the finite-difference Wronskian
$$
\cW(\ti\psi_1,\ti\psi_2)(z) := \det \begin{pmatrix}
\ti\psi_1(z-1) & \ti\psi_2(z-1) \\
\ti\psi_1(z) & \ti\psi_2(z) \end{pmatrix},
\qquad \ti\psi_1,\ti\psi_2\in\C[[z\ii]][z]
$$
or rather its Borel counterpart
$$
\cW(\tr\psi_1,\tr\psi_2) = (\ee^\ze\tr\psi_1)*\tr\psi_2 - \tr\psi_1*(\ee^\ze\tr\psi_2),
\qquad \tr\psi_1,\tr\psi_2\in\SING.
$$
One can indeed check that $\cW(\ti\ph_1,\ti\ph_2)=1$ and that
equation~\refeq{eqHenLinNHSG} is thus reduced to
\begin{equation}	\label{eqHenLinRED}
(\ee^{-\ze}-1)\trx c_1 = -\tr\psi*\trb\ph_2, \qquad 
(\ee^{-\ze}-1)\trx c_2 = \tr\psi*\trb\ph_1
\end{equation}
by the change of unknown
$$
\left\{\begin{aligned}
\trb\ph &= \trx c_1*\trb\ph_1 + \trx c_2*\trb\ph_2 \\[1ex]
\ee^\ze \trb\ph &= \trx c_1*(\ee^\ze\trb\ph_1) + \trx c_2*(\ee^\ze\trb\ph_2)
\end{aligned}\right.
\quad\Leftrightarrow\quad
\left\{\begin{aligned}
\trb c_1 &= \cW(\trb\ph,\trb\ph_2) \\[1ex]
\trb c_2 &= \cW(\trb\ph_1,\trb\ph).
\end{aligned}\right.
$$
For the homogeneous equation~\refeq{eqHenLinSG} we have $\tr\psi=0$, thus the
only solutions of~\refeq{eqHenLinRED} are 
$(\trx c_1,\trx c_2) = (C_1\,\de,C_2\,\de)$ with arbitrary $C_1,C_2\in\C$.
Indeed, the equation for~$\trx c_1$ for instance amounts to 
$\trx c_1 = \sing_0\bigl( \frac{\text{reg}(\ze)}{\ee^{-\ze}-1} \bigr)$ with an
arbitrary $\text{reg}(\ze)\in\C\{\ze\}$, of which the value at~$0$ will be
$-\frac{C_1}{2\pi\I}$.

\sm

We already saw that $\De_{2\pi\I}\trx x_0$ was solution of~\refeq{eqHenLinSG},
this yields complex numbers $C_1=\mu$ and $C_2=\Th$. 
Since~$\htx x_0$ is real-analytic and odd, it is purely imaginary on the
imaginary axis; since the coefficients of~$\ti\ph_1$ and~$\ti\ph_2$ are real,
the parity properties imply that $\mu\in\R$ and $\Th\in\I\,\R$.
The statement for~$\De_{-2\pi\I}\trx x_0$ is obtained by symmetry.
\eopf

The fact that $\De_{\pm2\pi\I}\trx x_0\in\SING\sram$ means that the first
singularities of~$\htx x_0$ are of the form \{polar part\}~$+$~\{logarithmic
singularity with regular variation\}. More precisely, with the notation
$$
\ti\ph_1(z) = \pa\ti x_0(z) = \sum_{k\ge1} b_k z^{-2k-1}, \qquad
\ti\ph_2(z) = \sum_{k\ge-2} d_k z^{-2k},
$$
the principal branch of~$\htx x_0$ satisfies
\begin{gather}	\label{eqsingx}
\htx x_0(2\pi\I+\ze) = \frac{\Th}{2\pi\I}\left(
d_{-2}\frac{4!}{\ze^5} + d_{-1}\frac{2!}{\ze^3} + d_0\frac{1}{\ze}
\right) + 
\frac{1}{2\pi\I} \bigl(\Th\htb\ph_2(\ze)+\mu\htb\ph_1(\ze)\bigr) \log\ze 
+ \text{reg}(\ze), \\
\htx x_0(-2\pi\I+\ze) = \frac{\Th}{2\pi\I}\left(
d_{-2}\frac{4!}{\ze^5} + d_{-1}\frac{2!}{\ze^3} + d_0\frac{1}{\ze}
\right) + 
\frac{1}{2\pi\I} \bigl(\Th\htb\ph_2(\ze)-\mu\htb\ph_1(\ze)\bigr) \log\ze 
+ \text{reg}(\ze).
\end{gather}
The coefficients~$b_k$ and~$d_k$ can be computed inductively, whereas the
constants~$\Th$ and~$\mu$ must be considered as transcendent: in the case of a
more general map~\refeq{henonGen}, there would be relations analogous
to~\refeq{eqfirstresrel} in which the corresponding constants do not depend on
finitely many coefficients of~$f$ only.
In the case of the H\'enon map, a positivity argument leads to
\begin{prop}	\label{propposit}
The constant $\Th\in\I\,\R$ determined in equation~\refeq{eqfirstresrel} of
Proposition~\ref{propHenLin} satisfies 
$$
\IM\Th<0.
$$
\end{prop}

\proof
We first show that $\ti x_0(z) = \sum_{k\ge1} a_k z^{-2k}$ with $(-1)^k a_k>0$.
Consider the auxiliary series 
$\ti U(t) = \ti x_0(\I t) = \sum_{k\ge1} (-1)^k a_k t^{-2k}$: it is the unique
nonzero even solution of 
$$
-\ti U(t+\I) + 2 \ti U(t) - \ti U(t-\I) = \ti U(t)^2.
$$
This equation can be rewritten 
$\pa_t^2\ti U = \Ga(\pa_t) (\ti U^2)$, with a convergent power series
$\Ga(X) = \frac{X^2}{4\sin^2\frac{X}{2}}= 1 + \Ga_1 X^2 + \Ga_2 X^4 + \cdots$
which has only non-negative coefficients (as can be seen from the decomposition
of the meromorphic function $1/\sin^2\frac{X}{2}$ as a series of second-order
poles).
One can thus write induction formulas for the coefficients of~$\ti U(t)$ which
show that they are positive.

\sm

As a consequence, the Borel transform $\hat U(\tau) = \I\,\htx x_0(\I\tau) = 
\sum_{k\ge1} (-1)^k a_k \frac{\tau^{2k-1}}{(2k-1)!}$ (which is convergent at
least in the disc of radius~$2\pi$) is positive and increasing on the
segment~$\left]0,2\pi\right[$.
But this function satisfies 
$$
\left(4\sin^2\frac{\tau}{2}\right) \hat U(\tau) = \hat U*\hat U(\tau),
$$
hence it cannot be bounded on~$\left]0,2\pi\right[$
(if it were, the \lhs\ would tend to~$0$ as $\tau\xrightarrow{<}2\pi$, whereas
the \rhs\ is positive increasing).

\sm

Now, in view of~\refeq{eqsingx}, the fact that~$\htx x_0(\ze)$ is not bounded
on~$\left]0,2\pi\I\right[$ shows that $\Th\neq0$ (because
$\htb\ph_1(\ze)=\cO(\ze^2)$).
Moreover, $0 < \hat U(\tau) = \I\,\htx x_0(\I\tau)
\sim \frac{1}{2\pi} \frac{\Th}{84} \frac{4!}{(\I\tau-2\pi\I)^5}$
for $\tau\xrightarrow{<}2\pi$ implies $\I\Th>0$.
\eopf

Numerically, one finds
$|\Th|\simeq 2.474 \cdot 10^6$, $\mu\simeq 4.909\cdot 10^3$
(much better accuracy can be achieved thanks to the precise information we
have on the form of the singularity---see~\cite{sauzin01}).

\subsub{The parabolic curves~$p^+(z)$ and~$p^-(z)$ and their splitting}

We pause here in the description of the resurgent structure of~$\ti x_0$ to give
a look at the analytic consequences we can already deduce from the above.

\sm

Borel-Laplace summation yields two analytic solutions~$x^+(z)$ and~$x^-(z)$ of
equation~\refeq{eqHenon}:
$$
x^\pm(z) = \cL^\pm \htx x_0(z) \sim \ti x_0(z), \qquad
z \in \cD^\pm
$$
(with notations analogous to those of Section~\ref{secAlCalAbEq}).
The analysis of the first singularities of~$\htx x_0(\ze)$ is sufficient to
describe the asymptotic behaviour of the difference $(x^+-x^-)(z)$ when
$z\in\cD^+\cap\cD^-$ with $\IM z>0$ or $\IM z<0$:
when $\IM z<0$ for instance,
one can argue as at the end of Section~\ref{secParBridge} and write
$$
(x^+-x^-)(z) = \int_{\ga_1}\ee^{-z\ze}\,\htx x_0(\ze)\,\dd\ze
+ \int_{\Ga}\ee^{-z\ze}\,\htx x_0(\ze)\,\dd\ze,
$$
with the same path~$\ga_1$ as on Figure~\ref{figStokes} (with $\om=0$), and with
a path~$\Ga$ coming from~$\ee^{\I\th'}\infty$, passing slightly below~$4\pi\I$
and going to~$\ee^{\I\th}\infty$.
The first integral is exactly $\ee^{-2\pi\I z}\cL^-(\De_{2\pi\I}\trx x_0)(z)$ (\cf
the section on the Laplace transform of majors), while the second integral is
$\cO(\ee^{-(4\pi-\de)|\IM z|})$ for any $\de>0$.
We thus obtain
$$
\ee^{2\pi\I z} (x^+-x^-)(z) \sim \De_{2\pi\I}\ti x_0(z) = \Th\ti\ph_2(z)+\mu\ti\ph_1(z),
\qquad z\in\cD^+\cap\cD^-, \quad \IM z<0
$$
(and similarly with $\De_{-2\pi\I}\ti x_0(z)$ for $\IM z>0$).
Since $\Th\neq0$, the leading term for this exponentially small discrepancy is
$$
x^+(z)-x^-(z) \sim \tfrac{\Th}{84}z^4\,\ee^{\pm2\pi\I z},
\qquad z\in\cD^+\cap\cD^-, \quad \pm\IM z>0.
$$
Observe that both~$x^+$ and~$x^-$ extend from~$\cD^\pm$ to entire functions, as
would any solution of the difference equation~\refeq{eqHenon} analytic in a
half-plane bounded by a non-horizontal line (for instance, the analytic
continuation of~$x^+$ is defined by iterating $x^+(z-1) = 2 x^+(z) - x^+(z+1) -
\bigl(x^+(z)\bigr)^2$). 
But the simple asymptotic behaviour of~$x^\pm$ described by~$\ti x_0$ does not
extend beyond~$\cD^\pm$: this is the Stokes phenomenon we have just described.

\sm

We finally return to the dynamics to the H\'enon map~$F$ and supplement the
solutions~$\ti x_0(z)$ and~$x^\pm(z)$ of the scalar equation~\refeq{eqHenon} with
the appropriate $y$-components, $\ti y_0=D\ti x_0$ and $y^\pm=D x^\pm$, in order
to define solutions of equation~\refeq{eqHenonSys}:

\begin{prop}	\label{propPC}
There exist two holomorphic maps $p^+:\C\to\C^2$ and $p^-:\C\to\C^2$ satisfying
equation~\refeq{eqHenonSys}:
$$
p^\pm(z+1) = F\bigl(p^\pm(z)\bigr), \qquad z\in\C,
$$
and admitting the same asymptotic expansion in different domains:
$$
p^\pm(z) \sim \ti p(z), \qquad z\in\cD^\pm,
$$
where $\ti p(z) = \bigl(-6z^{-2}+\cO(z^{-4}), 12z^{-3}+\cO(z^{-4})\bigr)$ is a
formal solution of~\refeq{eqHenonSys}.
Moreover, 
$$
\ee^{\pm2\pi\I z}\bigl(p^+(z)-p^-(z)\bigr) \sim
\Th \ti N(z) \pm \mu \frac{\dd\ti p}{\dd z}(z),
\qquad z\in\cD^+\cap\cD^-, \quad \pm\IM z<0,
$$
with the notations of Proposition~\ref{propHenLin} and 
$\ti N = ({\ti\ph_2},{D\ti\ph_2}) 
= \bigl(\frac{1}{84}z^4+\cO(z^2),\frac{1}{21}z^3+\cO(z^2)\bigr)$.
\end{prop}

Observe that the symplectic $2$-form $\dd x\wedge\dd y$ yields~$1$ when
evaluated on $\bigl( \frac{\dd\ti p}{\dd z}(z), \ti N(z) \bigr)$.
The constant~$\Th$ thus describes the normal component of the splitting,
while~$\mu$ describes the tangential component.

\sm

Because of equation~\refeq{eqHenonSys}, the curves $\cW^+ = \{ p^+(z), \;
z\in\C\}$ and $\cW^- = \{ p^-(z), \; z\in\C\}$ are invariant by~$F$ with
$$
F^n\bigl( p^\pm(z) \bigr) = p^\pm(z+n) \xrightarrow[n\to\pm\infty]{}0
$$
(in view of their common asymptotic series).
They may be called ``stable and unstable separatrices'' (by analogy with the
stable and unstable manifolds of a hyperbolic fixed point), or ``parabolic
curves'' (as is more common in the litterature on $2$-dimensional holomorphic
maps).

\subsub{Formal integral and Bridge equation}

We end with the complete description of the resurgent structure of the formal
solution~$\ti x_0(z)$ of equation~\refeq{eqHenon}.
Taking for granted the possibility of following the analytic continuation
of~$\htx x_0(\ze)$ in~$\cR$, as ascertained by Theorem~\ref{thmResurx}, and
consequently the possibility of defining the singularities
$\De_{\om_1}\cdots\De_{\om_r}\trx x_0$ for all $r\ge1$ and
$\om_1,\ldots,\om_r\in 2\pi\I\,\Z^*$, we see that the there must exist complex
numbers $A_{\om,0}$, $B_{\om,0}$ such that
$A_{\pm2\pi\I,0}=\pm\mu$, $B_{\pm2\pi\I,0}=\Th$ and
\begin{equation}	\label{eqresrel}
\De_\om\trx x_0 = A_{\om,0}\, \trb\ph_1 + B_{\om,0}\, \trb\ph_2, 
\qquad \om\in2\pi\I\,\Z^*.
\end{equation}
Indeed, the same arguments as those which led to Proposition~\ref{propHenLin}
apply and each $\De_\om\trx x_0$ is a solution of the linearized
equation~\refeq{eqHenLinSG}. 
But this is not a closed system of resurgence relations: in order to compute 
$\De_{\om_1}\De_{\om_2}\trx x_0$ for instance, we need to know the alien
derivatives of~$\trb\ph_2$ (on the other hand the alien derivatives
of~$\trb\ph_1$ can be deduced from~\refeq{eqresrel} and from the commutation
relation~\refeq{eqcommutnatur} of Proposition~\ref{propsimpleAlienD}).

\sm

This problem is solved by the notion of ``formal integral''. 
Recalling that $\ti\ph_1=\pa\ti x_0\in z^{-3}\C[[z\ii]]$ and setting $\ti
x_1=\ti\ph_2\in z^4\C[[z\ii]]$, we can consider $\ti x_0(z) + b\ti x_1(z)$,
where $b$ is a small deformation parameter, as a solution of
equation~\refeq{eqHenon} up to $\cO(b^2)$.
We can also consider this expression as the beginning of an exact solution
belonging to $\bigl(\C[[z\ii]][z]\bigr)[[b]]$:

\begin{prop}
There exist formal series $\ti x_2(z),\ti x_3(z),\ldots$ in~$\C[[z\ii]][z]$ such
that 
\begin{equation}	\label{eqformint}
\ti x(z,b) = \sum_{n\ge0} b^n \ti x_n(z)
\end{equation}
solves equation~\refeq{eqHenon}.
The $\ti x_n$'s are uniquely determined by the further requirement that they be
even and do not contain any $z^4$-term.
Moreover, 
$$
\ti x_n(z) \in z^{6n-2}\C[[z\ii]]_1, \qquad 
\trx x_n = \cB\ti x_n \in \trRES_{2\pi\I\Z}.
$$
\end{prop}

\proof
Plugging~\refeq{eqformint} into~\refeq{eqHenon} and expanding in powers of~$b$,
we recover the known equations $P\ti x_0=-\ti x_0^2$ and
$P\ti x_1+2\ti x_0\ti x_1=0$ at orders~$0$ and~$1$, and then a system of
inhomogeneous linear equations to be solved inductively:
\begin{equation}	\label{eqxn}
P\ti x_n+2\ti x_0\ti x_n = \ti\psi_n, 
\qquad \ti\psi_n = -\sum_{n_1=1}^{n-1} \ti x_{n_1}\ti x_{n-n_1},
\qquad n\ge2.
\end{equation}
In the course of the proof of Proposition~\ref{propHenLin}, we saw how to solve
such equations (admittedly their Borel counterparts, but this makes no
difference to the algebraic structure of their solution):
$\ti x_n = \ti c_1\ti\ph_1 + \ti c_2\ti\ph_2$ is solution of~\refeq{eqxn} as soon as
$$
\ti c_j(z+1)-\ti c_j(z) = \ti\chi_j(z), \qquad
\ti\chi_1 = -\ti\ph_2\ti\psi_n, \qquad
\ti\chi_2 = \ti\ph_1\ti\psi_n.
$$
Cancellations occur so that the coefficient of~$z\ii$ in~$\ti\chi_1$
and~$\ti\chi_2$ vanishes,\footnote{
See the proof of Proposition~6 in~\cite{sauzin01} (beware of the misprint in the
last formula on p.~551, which corresponds to an incorrect expansion of the (correct)
identity indicated in the footnote on that page).
Besides, these cancellations are special to the case of a symplectic
map~\refeq{henonGen}, \ie with~$f$ depending on its first argument only; the
formalism of Section~\ref{secSingul} is perfectly capable to handle the general
case with $\log z$ terms---see~\cite{GelfSauz}.
}
hence the primitives of~$\ti\chi_j$ have no $\log z$
term.
The end of the proof is just a matter of selecting appropriately the primitives
$\pa\ii\ti\chi_j$, setting $\ti c_j = \ga(-\pa)\,\pa\ii\ti\chi_j$ (with the
notations of the proof of Proposition~\ref{propHenLin}) and counting the
valuations---observe that the Borel tranforms $\trx c_j =
-\frac{\ze}{\ee^{-\ze}-1}\cB(\pa\ii\ti\chi_j)$, whence the analyticity in~$\cR$ of
the minors follows by induction.
The reader is referred to~\cite[\S 5.1]{sauzin01} for the details.
\eopf

The series $\ti x(z,b) = \sum b^n \ti x_n(z)$ is called a {\em formal integral} of
equation~\refeq{eqHenon}.
We thus have a kind of two-parameter family of solutions of~\refeq{eqHenon},
namely  $\ti x(z+a,b)$,
which is consistent with the fact that we are dealing with a second-order
equation.\footnote{
For the first-order difference equation~\refeq{eqniter}, we had the one-parameter family of
solutions~$\ti\ph(z+c)$.
}
At the level of the map~$F$, this corresponds to a parametrisation of the
``formal invariant foliation'' of the formal infinitesimal generator~$\ti X$.
Here is the implication for the resurgent structure of~$\ti x_0$:

\begin{thm}	\label{thmBridgeHen}
For each $\om\in 2\pi\I\,\Z^*$, there exist formal series 
$$
A_\om(b) = \sum_{n\ge0} A_{\om,n} b^n, \ens
B_\om(b) = \sum_{n\ge0} B_{\om,n} b^n
\,\in\, \C[[b]]
$$
such that 
\begin{equation}	\label{eqBridgeHen}
\De_\om \ti x(z,b) = \bigl( A_\om(b) \pa_z + B_\om(b)\pa_b \bigr)\ti x(z,b),
\end{equation}
this ``Bridge equation'' being understood as a compact writing of infinitely many
``resurgence relations''
\begin{equation}	\label{eqDetrxn}
\De_\om\trx x_n = \sum_{n_1+n_2=n} \bigl( A_{\om,n_1} \pa\trx x_{n_2}
+ (n_2+1) B_{\om,n_1} \trx x_{n_2+1} \bigr),
\qquad n\ge0.
\end{equation}
\end{thm}

\proof
The point is that $\pa_z\ti x(z,b)$ and~$\pa_b\ti x(z,b)$ are independent
solutions of the linearized equation $P\ti\ph + 2\ti x(z,b)\ti\ph =0$; what we
shall check amounts to the fact that their formal Borel transforms $\pa\trx x$
and~$\pa_b\trx x$, which lie in~$\SING\sram[[b]]$, span the space of solutions
of the equation
\begin{equation}	\label{eqLINbBT}
\al \, \trb\ph + 2\trx x * \trb\ph =0, \qquad \trb\ph\in\SING[[b]],
\end{equation}
a particular solution of which is $\De_\om\trx x = \sum b^n\De_\om\trx x_n$.

\sm

We prove~\refeq{eqDetrxn} by induction on~$n$ and suppose that the coefficients
of
$$
A^*(b) = \sum_{n=0}^{N-1} A_{\om,n} b^{n}, \quad
B^*(b) = \sum_{n=0}^{N-1} B_{\om,n} b^{n}
$$
were already determined so as to satisfy~\refeq{eqDetrxn} for $n=0,\ldots,N-1$.
The \rhs\ of~\refeq{eqDetrxn} with $n=N$ can be written as
$\trb\chi_N + A_{\om,N}\,\pa\trx x_0 + B_{\om,N}\,\trx x_1$, where
$\trb\chi_N\in\SING\sram$ is known in terms of the coefficients of~$\trx
x(z,b)$, $A^*(b)$ and~$B^*(b)$.
Thus, we only need to check that $\De_\om\trx x_N-\trb\chi_N$ is a
linear combination of $\pa\trx x_0=\tr\ph_1$ and $\trx x_1=\trb\ph_2$.

\sm

The singularity $\De_\om\trx x_N-\trb\chi_N$ is the coefficient of~$b^N$
in $\trb\ph = \De_\om\trx x - \bigl( A^*(b) \pa + B^*(b)\pa_b \bigr)\trx x$.
On the one hand, $\trb\ph = \cO(b^N)$ by our induction hypothesis; on the other
hand, $\trb\ph$ is solution of~\refeq{eqLINbBT} (because $\De_\om$, $\pa$
and~$\pa_b$ are derivations of~$\SING[[b]]$ that commute with the multiplication
by~$\al$).
Therefore $\De_\om\trx x_N-\trb\chi_N$ is solution of~\refeq{eqHenLinSG}, hence
of the form $A_{\om,N} \, \trb\ph_1 + B_{\om,N} \, \trb\ph_2$ according to
Proposition~\ref{propHenLin}.
\eopf

Now we can compute all the successive alien derivatives
$\De_{\om_1}\cdots\De_{\om_r}\trx x_n$ of all the components of the formal
integral.
Since the $\De_\om$'s commute with~$\pa_b$ and satisfy the commutation
rule~\refeq{eqcommutnatur} with~$\pa$, we get
\begin{multline*}
\De_{\om_1}\cdots\De_{\om_r}\trx x = 
\bigl( A_{\om_r}(\pa-\chx\om_{r-1}) + B_{\om_r}\pa_b \bigr)
\bigl( A_{\om_{r-1}}(\pa-\chx\om_{r-2}) + B_{\om_{r-1}}\pa_b \bigr)
\cdots\\
\cdots
\bigl( A_{\om_1}(\pa-\chx\om_0) + B_{\om_1}\pa_b \bigr) \trx x
\end{multline*}
with $\chx\om_j := \om_1+\cdots+\om_j$, $\chx\om_0:=0$.
By expanding in powers of~$b$, we have access to the resurgent structure of
each~$\trx x_n$. In particular, we see that all of them are simply ramified
resurgent functions:
$$
\trx x \,\in\, \trRsramZ\,[[b]].
$$


\subsection{Real splitting problems}

We now give a brief account of the work by V.~Lazutkin and V.~Gelfreich on the
exponentially small splitting for area-preserving planar maps, indicating the
connection with Section~\ref{secHenon}.

\subsub{Two examples of exponentially small splitting}

Let us consider a one-parameter family of real $2$-dimensional symplectic maps
\begin{equation*}	
\cG_\eps : \bino{X}{Y} \longmapsto
\bino{X + Y + \eps^2 g(X)}{Y + \eps^2 g(X)},
\end{equation*}
where $\eps\ge0$, and $g:\R\to\R$ is analytic with $g(X)=X+\cO(X^2)$.
More specifically, we have in mind two examples in which~$g$ is a (possibly
trigonometric) polynomial:
\begin{itemize}
\item
$g(X) = g\QU(X) = X(1-X)$ gives rise to $\cG_\eps=\cG_\eps\QU$, which is a
normal formal for non-trivial quadratic diffeomorphisms of the plane which are
symplectic and possess two fixed points (see~\cite[\S~1.2]{sauzin01} and the
references therein);
\item
$g(X) = g\SM(X) = \sin X$ gives rise to $\cG_\eps=\cG_\eps\SM$, the so-called
{\em standard map} (see~\cite{GL} for a survey and references to the
litterature); in this case we consider~$X$ as an angular variable, \ie the phase
space of~$\cG_\eps\SM$ is $(\R/2\pi\Z)\times\R$.
\end{itemize}
For $\eps=0$, we have an invariant foliation by horizontal lines, with fixed
points for $Y=0$.
In both examples, we wish to study the behaviour of~$\cG_\eps$ with $\eps>0$
small, in a region $|Y|=\cO(\eps)$ of the phase space.

\sm

Rescaling the variables by $x=X$, $y=\eps Y$, we get
\begin{equation*}	
\cF_\eps : \bino{x}{y} \longmapsto
\bino{x + \eps \bigl( y + \eps g(x) \bigr)}{y + \eps g(x)},
\end{equation*}
which has a hyperbolic fixed point at the origin.
But the hyperbolicity is weak when $\eps$ is small (the eigenvalues
are~$\ee^{\pm h}$ with $h=\eps+\cO(\eps^2)$) and this is the source of an exponentially
small phenomenon:
parts of the stable and unstable manifolds~$\cW_\eps^+$ and~$\cW_\eps^-$ of the
origin are very close one to the other (see~\cite{FS}). 
However, it turns out that they do not coincide and there is a homoclinic
point~$h_\eps$ at which they intersect transversely,\footnote{
Moreover, in both examples $\cW_\eps^-$ can be deduced from~$\cW_\eps^+$ by a
linear symmetry and~$h_\eps$ lies on the symmetry line.
}
with an exponentially small angle~$\al$: 
\begin{equation}	\label{eqasympangle}
\al\QU = \frac{\Om\QU}{\eps^7}\ee^{-\frac{2\pi^2}{\eps}} \bigl(1+\cO(\eps)\bigr),
\qquad
\al\SM = \frac{\Om\SM}{\eps^3}\ee^{-\frac{\pi^2}{\eps}} \bigl(1+\cO(\eps)\bigr),
\end{equation}
where $\Om\QU = \frac{64\pi}{9}|\Th|$ and $\Om\SM = \pi|\Th'|$, 
with the same constant~$\Th$ as the one discussed in
Propositions~\ref{propHenLin} and~\ref{propposit} in the study of the H\'enon
map of Section~\ref{secHenon}, and with an analogous constant~$\Th'$ stemming
from the study of another map without parameter.

\sm

The proof of the result by Vladimir Lazutkin and his coworkers in the case of
the standard map has a long story, which starts with his VINITI preprint in~1984
and ends with an article by Vassili Gelfreich in~1999 filling all the remaining
gaps. 
In the quadratic case (and also for other cases, corresponding to other
algebraic or trigonometric polynomials~$g(X)$ in the map~$\cG_\eps$), the result
was indicated by V.~Gelfreich without a complete proof (which should be, in
principle, a mere adaptation of the proof for the standard map).

\subsub{The map~$F$ as ``inner system''}

The strategy for proving~\refeq{eqasympangle} is to represent the invariant
manifolds~$\cW_\eps^\pm$ as parametrised curves
$P_\eps^\pm(t) = \bigl(x_\eps^\pm(t),y_\eps^\pm(t)\bigr)$, 
which are real-analytic, extend holomorphically to a half-plane $\pm\RE t\gg0$,
with $x_\eps^\pm(t) \sim \text{const}\,\ee^t$ as $\pm\RE t\to\infty$, and
satisfy $P_\eps^\pm(t+\eps) = \cF_\eps\bigl(P_\eps^\pm(t)\bigr)$.
The second component can be eliminated: $y_\eps^\pm(t) =
x_\eps^\pm(t)-x_\eps^\pm(t-\eps)$, and the first is then solution of a nonlinear
second-order difference equation with small step size:
\begin{equation}	\label{eqdiffeqnsmstp}
x_\eps^\pm(t+\eps) - 2 x_\eps^\pm(t) + x_\eps^\pm(t-\eps) = 
\eps^2 g\bigl( x_\eps^\pm(t) \bigr).
\end{equation}
Pictures in~\cite{GL} or~\cite{sauzin01} show that both curves are close to the
separatrix solution~$P_0(t)$ of the Hamiltonian vector field generated by
$H(x,y) = \demi y^2 + G(x)$, where $G'=-g$.
Indeed, the differential equation $\frac{\dd^2 x}{\dd t^2} = g(x)$,
which is the singular limit of~\refeq{eqdiffeqnsmstp} when $\eps\to 0$, 
has a particular solution~$x_0(t)$ wich tends to~$0$ both for $t\to+\infty$ and
$t\to-\infty$ (with the identification $0\equiv2\pi$ in the case of the standard
map), corresponding to the upper part of the separatrix of the pendulum in the
second example and to a homoclinic loop in the first one:
\begin{equation}	\label{eqxQUxSM}
x_0\QU(t) = \frac{3}{2\cosh^2\tfrac{t}{2}}
\quad\text{or}\quad
x_0\SM(t) = 4 \arctan \ee^t.
\end{equation}
The exponentially small phenomenon that we wish to understand can be described
as the splitting of this separatrix:
when passing from the Hamiltonian flow corresponding to~$\eps=0$ to the
map~$\cF_\eps$ with $\eps>0$, the separatrix solution~$x_0(t)$ is replaced by
two solutions~$x_\eps^+(t)$ and~$x_\eps^-(t)$ of
equation~\refeq{eqdiffeqnsmstp}, the difference of which is exponentially small
but not zero.
Proving~\refeq{eqasympangle} amounts essentially to estimating~$x_\eps^+(t)-x_\eps^-(t)$.

\sm

It is easy to see that the separatrix must split for an entire map
like~$\cF_\eps$: since $g$ is entire, equation~\refeq{eqdiffeqnsmstp} affords
analytic continuation to the whole $t$-plane of the functions~$x_\eps^+$
and~$x_\eps^-$ and Liouville's theorem prevents them from coinciding
(see~\cite{Ush} for a more general result of the same kind).

\sm

The method which was successfully developed by Vladimir Lazutkin and his
coworkers can then be described in the language of complex matching asymptotics.
The separatrix solution~$x_0(t)$ is a good approximation
of~$x_\eps^\pm(t)$ for real values of~$t$, but it is no so in the complex
domain, if only because~$x_\eps^+$ and~$x_\eps^-$ are entire functions, whereas~$x_0$ has
singularities in the complex plane.
It turns out that the asymptotics that we want to capture is governed by the
singularities of~$x_0$ which are closest to the real axis,
$$ 
t^*=\I\pi \quad\text{or}\quad \frac{\I\pi}{2}, 
\qquad \text{and}\quad \bar t^*=-t^*
$$
(the exponential in the final result~\refeq{eqasympangle} is nothing
but~$\ee^{\frac{2\pi\I t^*}{\eps}}$).
One is thus led to look for a better approximation of~$x_\eps^\pm(t)$ when $t$
lies in a domain called ``inner domain'', close to~$t^*$.

\sm

In the case of the quadratic map, one can perform the change of variables
$u=\eps^2 x-\frac{\eps^2}{2}$, $v=\eps^3 y$,
which transforms~$\cF_\eps$ into
\begin{equation*}
F_\eps : \bino{u}{v} \longmapsto
\bino{u + v -u^2 +\tfrac{\eps^4}{4}}{v -u^2 +\tfrac{\eps^4}{4}}.
\end{equation*}
We recognize here an $\eps^4$-perturbation of the H\'enon map of
Section~\ref{secHenon} (\ie the map~\refeq{henonGen} with $f(x,y)=-x^2$), and
the second-order difference equation with small step size~\refeq{eqdiffeqnsmstp}
can be treated as a perturbation of equation~\refeq{eqHenon} with $t = t^* +
\eps z$.
The H\'enon map, being obtained by forgetting the $\eps^4$-terms in~$F_\eps$
and providing good complex approximations of~$x_\eps^+(t)$ and~$x_\eps^-(t)$
through its parabolic curves, is called the {\em inner system} of the
family~$\cF_\eps$. 

\sm

The parabolic fixed point of the H\'enon map appears as the
organising centre of our perturbation problem; the detailed study of the
parabolic curves, including Proposition~\ref{propPC}, is an important step in
the proof of the estimate~\refeq{eqasympangle} for the splitting in the
quadratic case; it is here that the non-zero constant~$\Th$ appears.
In the case of the standard map, the inner system is the so-called
``semi-standard map'', which one could also study by resurgent methods, although
this is not exactly what Lazutkin and his coworkers did.

\sm

Once the separatrix splitting for the inner map has been analysed, an extra work
is required to reach the result for~$\cG_\eps$, which was not yet written in
full details for the quadratic family~$\cG_\eps\QU$, contrarily to the case of
the standard map~$\cG_\eps\SM$.

\subsub{Towards parametric resurgence}

One can suggest another approach to the proof of~\refeq{eqasympangle}.
Equation~\refeq{eqdiffeqnsmstp} admits a formal solution
\begin{equation}	\label{eqformalx} 
\ti x_\eps(t) = \sum_{n\ge0} \eps^{2n} x_n(t),
\end{equation}
where the first term is the separatrix solution~$x_0$ of the limit flow, and
where the subsequent functions can be computed inductively. 
The formal solution is unique if one imposes the condition $x_n(0)=0$ for each $n\ge1$.
One finds that each function~$x_n(t)$ tends to~$0$ when~$t$ tends to~$+\infty$
or to~$-\infty$. The formal solution is thus a candidate to represent the stable
solution~$x_\eps^+(t)$ and the unstable one~$x_\eps^-(t)$ as well.
What happens is that the formal series~\refeq{eqformalx} is divergent for
$t\neq0$ and only provides an asymptotic expansion both for~$x_\eps^+(t)$
and~$x_\eps^-(t)$.

\sm

The formal series~\refeq{eqformalx} truncated far enough is used in the work by
Lazutkin et al.\ as an approximation of~$x_\eps^\pm(t)$, but one can envisage a
more radical use of this formal solution. Unpublished computations performed in
collaboration with Stefano Marmi tend to indicate that it is resurgent in
$z=\frac{1}{\eps}$, that~$x_\eps^+(t)$ and~$x_\eps^-(t)$ can be recovered from
it by Borel-Laplace summation for $\eps>0$ as
$$
x_\eps^+(t) = x_0(t)+\int_0^{+\infty}\ee^{-\ze/\eps}\,\htx x(\ze,t)\,\dd\ze	
\quad\text{if $t>0$}, \quad
x_\eps^-(t) = x_0(t)+\int_0^{+\infty}\ee^{-\ze/\eps}\,\htx x(\ze,t)\,\dd\ze	
\quad\text{if $t<0$},
$$
where $\dst\htx x(\ze,t) = \sum_{n\ge1}\frac{\ze^{2n-1}}{(2n-1)!}x_n(t)$,
and that 
$$
x_\eps^+(t)-x_\eps^-(t) \sim \ee^{-{\om_*(t)}/{\eps}}\,\De_{\om_*(t)}\ti x
-\ee^{-{\ov{\om}_*(t)}/{\eps}}\,\De_{\ov{\om}_*(t)}\ti x,
\qquad t>0,
$$
where $\om_*(t) = 2\pi\I(t-t^*)$ and $\ov{\om}_*(t) = -2\pi\I(t+t^*)$ are the
singularities of~$\htx x(\,.\,,t)$ with the smallest positive real part.

\sm

We shall try to explain in the next section why singularities of the principal
branch of~$\htx x(\ze,t)$ should appear at the points 
$\om_{a,b}(t) = 2\pi\I a \bigl(t - (2b+1)t^*\bigr)$, $a\in\Z^*$, $b\in\Z$.
This kind of resurgence we expect is called {\em parametric resurgence}, because
the resurgent variable~$1/\eps$ appears as a parameter in
equation~\refeq{eqdiffeqnsmstp}, whereas the ``dynamical'' variable is~$t$.
This makes the analysis more complicated, since, for instance, the singular
points of the formal Borel transform \wrt~$1/\eps$ are ``moving singular
points'' (they move along with~$t$).

\sm

In the above conjectural statements, the Borel summability is the most
accessible. 
It amounts to the fact that, for fixed~$t$, the formal Borel transform~$\htx
x(\ze,t)$ has positive radius of convergence and extends holomorphically, with
at most exponential growth at infinity, to the branch cut obtained by removing
from the complex plane the moving singular half-lines
$\pm\om_{1,b}(t)\left[1,+\infty\right[$.
See~\cite[\S~5.3]{MarmiSauzin} for such a statement concerning a slightly
simpler second-order difference equation with small step size (related to the
semi-standard map).


\subsection{Parametric resurgence for a cohomological equation}

Since we just alluded to a possible phenomenon of parametric resurgence, we end
with two linear examples which are comparable to equations~\refeq{eqnph}
and~\refeq{eqnpsi} of Section~\ref{seclineardiffeq} and can help to understand
the origin of this phenomenon in difference equations with small step size.

\begin{prop}	\label{proplinpar}
Let $U$ be an open connected and simply connected subset of~$\C$ and let
$\cH(U)$ denote the space of all the functions holomorphic in~$U$.
Let $\al,\be\in\cH(U)$ and consider the two difference equations
\begin{equation}	\label{eqlinpar}
\ph(t+\eps)-\ph(t) = \eps \al(t), \qquad
\psi(t+\eps)-2\psi(t)+\psi(t-\eps) = \eps^2 \be(t).
\end{equation}
Then there exist sequences $(\ph_n)_{n\ge1}$ and $(\psi_n)_{n\ge1}$ of elements
of~$\cH(U)$ such that the solutions of the first equation in~$\cH(U)[[\eps]]$
are the formal series
$$
\ti\ph_\eps(t) = c_0(\eps) + \sum_{n\ge0} \eps^n\ph_n(t),
$$
where $c_0(\eps) \in \C[[\eps]]$ is arbitrary and $\ph_0$ is any primitive of~$\al$
in~$U$,
and the solutions of the second equation in~$\cH(U)[[\eps]]$ are the formal
series
$$
\ti\psi_\eps(t) = c_1(\eps) + c_2(\eps) t + \sum_{n\ge0} \eps^{2n}\psi_n(t),
$$
where $c_1(\eps),c_2(\eps) \in \C[[\eps]]$ are arbitrary and $\psi_0$ is any
second primitive of~$\be$ in~$U$ (\ie $\psi_0''=\be$).
Moreover the formal Borel transforms
$$
\hat\ph(\ze,t) = \sum_{n\ge1} \frac{\ze^{n-1}}{(n-1)!}\ph_n(t), \qquad
\hat\psi(\ze,t) = \sum_{n\ge1} \frac{\ze^{2n-1}}{(2n-1)!}\psi_n(t)
$$
have positive radius of convergence for any fixed $t\in U$ and define
holomorphic functions of~$(\ze,t)$, which can be expressed as 
$$
\hat\ph(\ze,t) = -\demi\al(t) - \sum_{\nu\in2\pi\I\Z^*} \frac{1}{\nu} \left(
\al\bigl(t+\tfrac{\ze}{\nu}\bigr) - \al(t) \right),
\qquad \hat\psi(\ze,t) = \sum_{\nu\in2\pi\I\Z^*} \frac{\ze}{\nu^2}\,\be\bigl(t+\tfrac{\ze}{\nu}\bigr)
$$
and extend holomorphically to 
$$
\cE = \ao (\ze,t) \in\C\times U \mid \bigl[
t-\tfrac{\ze}{2\pi\I},t+\tfrac{\ze}{2\pi\I} \bigr] \subset U \af.
$$
\end{prop}

When $\al(t)=g(\ee^t)$ with $g(z)\in z\C\{z\}$, equation~\refeqa{eqlinpar}{a} is a
particular case of the {\em cohomological equation} studied
in~\cite{MarmiSauzin}. This is the equation
$$
f(qz) - f(z) = g(z),
$$
with an unknown function $f(z)\in z\C\{z\}$ depending on the complex
parameter~$q$.
The article~\cite{MarmiSauzin} investigates the nature of the dependence of~$f$
upon~$q$ when~$q$ crosses the unit circle. Roots of unity appear as
``resonances'', at which the solution~$f(q,z)$ has sorts of simple poles:
$f(q,z) = \sum \frac{\cL_\La(z)}{q-\La}$ 
with a sum extending to all roots of unity~$\La$ and explicit
``residua''~$\cL_\La(z)$, except that these singular points are not isolated and
the unit circle is to be considered as a natural boundary of analyticity.
Still, $(q-\La)f(q,z)$ tends to~$\cL_\La(z)$ when~$q$ tends to~$\La$
non-tangentially \wrt\ the unit circle (uniformly in~$z$), and there is even a
kind of Laurent series at~$\La$: an asymptotic expansion~$\ti f_\La$ which is valid
near~$\La$, inside or outside the unit circle, and which must be divergent due
to the presence of arbitrarily close singularities.

\sm

The asymptotic series~$\ti f_\La$ can be found by setting $q = \La\,\ee^\eps$
and $\ph(t) = \eps f(z)$ (notice that $\eps\sim\frac{q-\La}{\La}$); this yields
the equation
\begin{equation}	\label{eqparlinres}
\ph(t + 2\pi\I\tfrac{N}{M} + \eps) - \ph(t) = \eps \al(t),
\end{equation}
where $\La = \exp\left(2\pi\I\frac{N}{M}\right)$.
The formal solution corresponding to~$\ti f_\La$ and its Borel transform are
described in~\cite[Chap.~4]{MarmiSauzin}, whith a statement which generalizes
Proposition~\ref{proplinpar}.

\sm

If~$\al(t)$ is a meromorphic function,\footnote{
This is the case when $\al(t)=g(\ee^t)$ with a meromorphic function~$g$;
this is worth of interest since $g_*(z)=\frac{z}{1-z}$ is the unit of the
Hadamard product, from which the case of any $g(z)\in z\C\{z\}$ can be deduced
(see~\cite{MarmiSauzin}).
}
then so is~$\hat\ph(\ze,t)$ as a function of~$\ze$ for any fixed~$t$,
and similarly with~$\hat\psi(\ze,t)$ when~$\be(t)$ is meromorphic.
This provides us with elementary examples of parametric resurgence.

\med

\noindent {\em Proof of Proposition~\ref{proplinpar}.}
As previously mentioned, equation~\refeqa{eqlinpar}{a} and the more general
equation~\refeq{eqparlinres} are dealt with in~\cite{MarmiSauzin}, so we content
ourselves with the case of equation~\refeqa{eqlinpar}{b} (which is not very
different anyway).

\sm

The equation can be written 
$\left( 4\sinh^2 \frac{\eps\pa_t}{2} \right) \psi = \eps^2\be$. 
We set
$4\sinh^2 \frac{X}{2} = \frac{X^2}{1+\Ga(X)}$
and introduce the Taylor series
$$
\Ga(X) = - 1 + \frac{X^2}{4\sinh^2 \frac{X}{2} } = \sum_{n\ge0} \Ga_n X^{2n+2}
$$
and its Borel transform
$\hat\Ga(\xi) = \sum_{n\ge0} \Ga_n \frac{\xi^{2n+1}}{(2n+1)!}$.
With these notations, 
$$
\refeqa{eqlinpar}{b} \ens\Leftrightarrow\ens
\pa_t^2\psi = \be + \Ga(\eps\pa_t)\be
\ens\Leftrightarrow\ens
\psi = \psi_0 + \sum_{n\ge0} \Ga_n\, \eps^{2n+2} \pa_t^{2n}\be
\ens\text{with}\ens \pa_t^2\psi_0 = \be.
$$
This gives all the formal solutions, with $\psi_n =
\Ga_{n-1}\,\pa_t^{2(n-1)}\be$ for $n\ge1$. Now,
\begin{equation}	\label{eqhatpsi}
\hat\psi(\ze,t) = \sum_{n\ge0} \Ga_n \frac{\ze^{2n+1}}{(2n+1)!}\pa_t^{2n}\be 
= \hat\Ga(\ze\pa_t)\,\pa_t\ii\be,
\end{equation}
where $\pa_t\ii\be$ denotes any primitive of~$\be$.
A classical identity yields
$$
\Ga(X) = \sum_{\nu\in2\pi\I\Z^*} \frac{X^2}{(X-\nu)^2} 
= \sum_{\nu\in2\pi\I\Z^*} \nu^{-2}X^2 (1-\nu\ii X)^{-2}
= \sum_{\nu\in2\pi\I\Z^*} \sum_{n\ge0} (n+1) \nu^{-n-2} X^{n+2},
$$
whence $\hat\Ga(\xi) = \sum \nu^{-2} \, \xi \,\ee^{\nu\ii\xi}$.
Inserting this into~\refeq{eqhatpsi} and using the Taylor formula in the form
$\ee^{\nu\ii\ze\,\pa_t}\,\be = \be(t+\nu\ii\ze)$ yields the conclusion.
\eopf

\med

We can now explain why we expect parametric resurgence for the formal
solution~$\ti x_\eps(t)$ of equation~\refeq{eqdiffeqnsmstp}, with singularities
of the Borel transform located at the points $\om_{a,b}(t) = 2\pi\I a \bigl(t -
(2b+1)t^*\bigr)$, $a\in\Z^*$, $b\in\Z$.
The idea is simply that equation~\refeqa{eqlinpar}{b} with $\be(t) =
g\bigl(x_0(t)\bigr)$ can be considered as a non-trivial approximation
of~\refeq{eqdiffeqnsmstp}; but $g\circ x_0 = \frac{\dd^2 x_0}{\dd t^2}$ is
meromorphic in the two cases we are interested in, in view of~\refeq{eqxQUxSM},
with $t^*+2t^*\Z$ as set of poles.
This yields a resurgent solution~$\ti\psi_\eps(t)$, with a meromorphic Borel
transform~$\hat\psi(\ze,t)$ for which the set of poles is
$2\pi\I\,\Z^*(-t+t^*+2t^*\Z)$. We expect~$\htx x(\ze,t)$ to resemble
somewhat~$\hat\psi(\ze,t)$; we do not expect it to be meromorphic, because of
the nonlinear character of equation~\refeq{eqdiffeqnsmstp}, but a majorant
method can be devised to control at least its principal branch.


\pagebreak

\end{document}